\documentclass[lettersize,journal]{IEEEtran}
\usepackage{amsmath,amssymb,amsfonts}
\usepackage{array}
\usepackage[caption=false,font=normalsize,labelfont=sf,textfont=sf]{subfig}
\usepackage{textcomp}
\usepackage{stfloats}
\usepackage{url}
\usepackage{verbatim}
\usepackage{graphicx}
\usepackage{xcolor,float}
\usepackage{bm}
\usepackage{algorithmic,algorithm}

\newtheorem{property}{Property}
\newtheorem{assumption}{Assumption}
\newtheorem{theorem}{Theorem}
\newtheorem{lemma}{Lemma}
\newtheorem{remark}{Remark}

\def\BibTeX{{\rm B\kern-.05em{\sc i\kern-.025em b}\kern-.08em
    T\kern-.1667em\lower.7ex\hbox{E}\kern-.125emX}}
\usepackage{balance}
\begin{document}
\title{Byzantine-Resilient Decentralized Online \\ Resource Allocation}
\author{Runhua Wang, Qing Ling, Hoi-To Wai and Zhi Tian
\thanks{Runhua Wang and Qing Ling are with the School of Computer Science and Engineering and the Guangdong Provincial Key Laboratory of Computational Science, Sun Yat-Sen University, Guangzhou, Guangdong 510006, China.
Hoi-To Wai is with the Department of Systems Engineering and Engineering Management, The Chinese University of Hong Kong, Hong Kong 999077, China.
Zhi Tian is with the Department of Electrical and Computer Engineering, George Mason University, Fairfax, VA 22030, USA.
Corresponding author: Qing Ling (lingqing556@mail.sysu.edu.cn).}\thanks{Qing Ling (corresponding author) is supported by National Key R\&D Program of China grant 2024YFA1014002, NSF China grant 62373388, Guangdong Basic and Applied Basic Research Foundation grant 2023B1515040025, and Guangdong Provincial Key Laboratory of Mathematical Foundations for Artificial Intelligence grant 2023B1212010001. 
A short, preliminary version of this paper has appeared in DSP 2025 \cite{b-RunhuaWang2025}.}
}


\maketitle

\begin{abstract}
In this paper, we investigate the problem of decentralized online resource allocation in the presence of Byzantine attacks. In this problem setting, some agents may be compromised due to external manipulations or internal failures, causing them to behave maliciously and disrupt the resource allocation process by sending incorrect messages to their neighbors. Given the non-consensual nature of the resource allocation problem, we formulate it under a primal-dual optimization framework, where the dual variables are aggregated among the agents, enabling the incorporation of robust aggregation mechanisms to mitigate Byzantine attacks. By leveraging the classical Byzantine attack model, we propose a class of Byzantine-resilient decentralized online resource allocation algorithms that judiciously integrate the adaptive robust clipping technique with the existing robust aggregation rules to filter out adversarial messages. We establish theoretical guarantees, showing that the proposed algorithms achieve tight linear dynamic regret and accumulative constraint violation bounds, where the constants depend on the properties of robust aggregation rules. Numerical experiments on decentralized online economic dispatch validate the effectiveness of our approach and support our theoretical results.
\end{abstract}


\section{Introduction}\label{sec 1}
Decentralized online resource allocation seeks to determine an optimal sequence of resource allocation strategies that satisfy long-term time-varying global resource constraints and local resource constraints, while minimizing the accumulative time-varying agent costs or maximizing the accumulative time-varying agent utilities over a given time horizon.
It arises in various application scenarios, such as smart grids \cite{b-Guo-Chen-2016, b-Zhenwei-Guo-2021}, cloud computing \cite{b-Chunrong-Wu-2021} and wireless communications \cite{b-Shichao-Xia-2021}.
Solving the decentralized online resource allocation problem relies on information exchange among neighboring agents. However, such exchange is not always reliable, as some agents may behave maliciously and send incorrect messages to their neighbors due to faults, communication failures, or cyber attacks. For instance, recent years have witnessed a surge in cyber incidents targeting power systems, such as the 2022 Ukraine power grid attack and 2025 India power grid attack. Such malicious behaviors can result in unfavorable online resource allocation strategies. In smart grids, this may lead to serious consequences such as large-scale blackouts. Therefore, this paper aims to study resilient decentralized online resource allocation algorithms to mitigate negative impacts caused by malicious agents.

Decentralized online resource allocation belongs to the broad class of decentralized constrained online convex optimization \cite{b-Xiuxian-Li2023}, for which the constraints can be either consensual \cite{b-Xuanyu-Cao2021, b-Parvin-Nazari2022,b-Xingrong-Dong-2024, b-Qingyang-Sheng2025} or non-consensual \cite{b-Xuanyu-Cao2022, b-Wenjing-Yan2024}.
It is typically modeled as decentralized online convex optimization subject to time-varying, coupled, non-consensual equality constraints.
The performance metrics are static/dynamic regret and accumulative constraint violation. Static regret compares the accumulated cost with the optimal cost relative to an optimal strategy in hindsight, which is constant over the entire time horizon. For dynamic regret, in contrast, the baseline becomes a series of instantaneous optimal strategies. The work of \cite{b-Jueyou-Li2022} proposes a decentralized online primal-dual algorithm with gradient feedback, and establishes its sublinear static regret and accumulative constraint violation. Similarly, the work of \cite{b-Xinlei-Yi2020} proposes a decentralized online primal-dual dynamic mirror descent algorithm. Sublinear static and dynamic regrets, as well as sublinear accumulative constraint violation, are established.
Unlike \cite{b-Jueyou-Li2022} and \cite{b-Xinlei-Yi2020} that rely on doubly stochastic mixing matrices to aggregate messages of neighboring agents,  \cite{b-Keishin-Tada2024} proposes a decentralized online primal-dual subgradient algorithm based on a row-stochastic mixing matrix,
and proves that the algorithm achieves sublinear dynamic regret and accumulative constraint violation.
All the aforementioned online algorithms rely on gradient or subgradient feedback from the cost functions. Several other studies focus on bandit feedback and propose decentralized one-point \cite{b-Xinlei-Yi2021a} and two-point \cite{b-Xinlei-Yi2023, b-Kunpeng-Zhang2024} online algorithms.
These algorithms also achieve sublinear dynamic regret and accumulative constraint violation.

When the agents are reliable, the online algorithms discuss- ed above can solve the decentralized online resource allocation problem. 
Nevertheless, some agents may behave maliciously, transmit incorrect messages to their neighboring agents, and consequently, disrupt the optimization process for decentralized online resource allocation.
We use the classical Byzantine attack model to describe the malicious behaviors of such agents, referring to them as Byzantine agents \cite{b-LeslieLamport-1982, b-ZhixiongYang-2020}. 

Resilience to Byzantine attacks has been extensively studied in decentralized multi-agent consensus optimization.
The basic idea is to aggregate messages in a robust manner. The works of \cite{b-Lili-Su-2021} and \cite{b-Zhixiong-Yang-2019} propose to implement trimmed mean ($TM$), a robust aggregation rule, to filter out erroneous scalar messages.
In $TM$, a benign agent discards the smallest and largest $b$ messages among those received from its neighbors, and then averages the remaining ones and its own message, where $b$ denotes the upper bound on the number of Byzantine agents.
When dealing with high-dimensional optimization variables, $TM$ is extended to coordinate-wise trimmed mean ($CTM$), in which $TM$ is executed in each dimension.
Another robust aggregation rule, iterative outlier scissor ($IOS$), is introduced in \cite{b-ZhaoxianWu-2023}.
In $IOS$, a benign agent iteratively discards $b$ messages that are the farthest from the average of the remaining ones.
In the work of \cite{b-LieHe-2022}, a benign agent employs the self-centered clipping ($SCC$) aggregation rule, which clips the received messages and calculates a weighted average.

The Byzantine-resilient decentralized consensus optimization algorithms proposed in \cite{b-Lili-Su-2021,b-Zhixiong-Yang-2019,b-ZhaoxianWu-2023,b-LieHe-2022} are not applicable to the resource allocation problem, which is not in a consensus form and has coupled constraints.
In the offline setup, the remedy is the combination of primal-dual algorithms and robust aggregation rules.
A Byzantine-resilient primal-dual algorithm is developed in \cite{b-Berkay-Turan-2021}.
Nevertheless, it must rely on a central server.
The work of \cite{b-Runhua-Wang-2022} proposes a Byzantine-resilient decentralized resource allocation algorithm (BREDA), which uses $CTM$ to defend against Byzantine attacks.
The work of \cite{b-RunhuaWang-2023} extends BREDA through incorporating a wide class of robust aggregation rules.
But unfortunately, Byzantine-resilient decentralized resource allocation algorithms in the online setup are still lacking.

In this paper, we address the less-studied Byzantine-resilient decentralized online resource allocation problem and make the following contributions.

\noindent \textbf{C1)} We propose a class of Byzantine-resilient decentralized online resource allocation algorithms that achieve Byzantine-resilience by employing a variety of well-designed appropriate robust aggregation rules on dual variables to filter out erroneous messages. 
In particular, we combine the adaptive robust clipping ($ARC$) technique with existing robust aggregation rules such as $IOS$ and $SCC$, and prove them to be applicable to the online setup in arbitrary dimensions.

\noindent \textbf{C2)} We analytically prove that the proposed algorithms have linear dynamic regret and accumulative constraint violation, which are inevitable under Byzantine attacks. The associated constants are determined by the properties of the well-designed robust aggregation rules.
We conduct numerical experiments on decentralized online economic dispatch to verify the theoretical results.

Compared to the short, preliminary version of this paper \cite{b-RunhuaWang2025}, this journal article presents some enhancements.
It incorporates thorough derivations related to algorithm development, expanded theoretical analysis, and new numerical experiments.
These additions strengthen the theoretical foundation and improve the applicability of the proposed algorithms. 

\noindent\textbf{Paper Organization:} This paper is organized as follows. In Section \ref{sec 2}, we formulate the decentralized online resource allocation problem under Byzantine attacks. Section \ref{sec 3} gives an attack-free decentralized online resource allocation algorithm, and shows its failure under Byzantine attacks. Section \ref{sec 4} further proposes a class of Byzantine-resilient decentralized online resource allocation algorithms. Section \ref{sec 5} analyzes the performance of proposed algorithms. Numerical experiments are given in Section \ref{sec 6}. Section \ref{sec 7} concludes this paper.

\noindent\textbf{Notations:} Throughout this paper, $(\cdot)^{ \top }$ stands for the transposition of a vector or a matrix, $\|\cdot\|$ stands for the $\ell_{2}$-norm of a vector, $\| \cdot \|$ is the operator norm of a matrix induced by $\ell_2$ norm, $\|\cdot\|_{F}$ denotes the Frobenius norm of a matrix, and $\left \langle \cdot , \cdot \right \rangle $ represents the inner product of vectors. We define $\widetilde{\bm{1}}\in \mathbb{R}^{M} $ and $\bm{1}\in \mathbb{R}^{H} $ as all-one column vectors while $I\in \mathbb{R}^{H \times H}$ as an identity matrix, where $M$ is the number of all agents and $H$ is the number of benign agents, respectively.

\section{PROBLEM STATEMENT}\label{sec 2}
We consider a decentralized online resource allocation problem involving $M$ agents.
The decentralized network is modeled as an undirected, connected graph $\widetilde{\mathcal{G}} (\mathcal{M}, \widetilde{\mathcal{E}})$, in which $\mathcal{M}$ represents the set of agents, and $\widetilde{\mathcal{E}}$ denotes the set of communication edges.
If two agents $i$ and $j$ can communicate with each other, then $(i,j) \in \widetilde{\mathcal{E}}$.
The set of neighbors of agent $i \in \mathcal{M}$ is denoted as $\mathcal{N}_{i}=\{j|(i,j)\in \widetilde{\mathcal{E}}\}$.
The decentralized online resource allocation problem aims to determine an optimal sequence of resource allocation strategies that minimize the sum of time-varying agent costs over a given time horizon, while satisfying both long-term global resource constraints and local resource constraints.
Denote the time horizon as $[1,T]$, where $T$ is the total number of time periods.
In each time period $t \in [1,T]$, $P_{i}^{t}\in \mathbb{R}^{d}$ represents the resource allocation strategy of agent $i \in \mathcal{M}$ and belongs to a compact, convex local resource constraint $\Omega_{i}$.
Let $\frac{1}{M}\sum_{i\in \mathcal{M}}P_{i}^{t}$ denote the average resource in each time period $t$, and let $D^{t}\in \mathbb{R}^{d}$ be the time-varying average resource constraint vector. Then, over the time horizon $[1,T]$, the long-term global resource constraint is $\sum_{t=1}^T \frac{1}{M}\sum_{i\in \mathcal{M}} P_i^t=\sum_{t=1}^T D^{t}$.
Each agent $i\in \mathcal{M}$ has a time-varying convex and continuously differentiable cost function $C_{i}^{t}(P_{i}^{t})$.

Within a given time horizon $[1,T]$, in the decentralized online resource allocation problem, uncertainties arise from a sequence of time-varying cost functions $\{C_{i}^{t}(P_{i}^{t}), \forall i \in \mathcal{M}\}_{t=1}^{T}$ and average resource constraint $\{D^{t}\}_{t=1}^{T}$, which are unknown in advance and sequentially disclosed over time.
The online resource allocation problem is formulated as
\begin{align}
\label{online-decentralized-economic-dispatch-problem}
 \min_{\{\widetilde{\bm{P}}^{t}\}_{t=1}^{T}}   & \quad \sum_{t=1}^{T}\widetilde{C}^{t}(\widetilde{\bm{P}}^{t}) ~~{\rm with}~~ \widetilde{C}^{t}(\widetilde{\bm{P}}^{t})= \sum_{i\in \mathcal{M}} C_{i}^{t}(P_{i}^{t}), \\
 s.t. & \quad \sum_{t=1}^{T}\sum_{i\in \mathcal{M}}\widetilde{G}_{i}^{t}(P_{i}^{t})=0~~{\rm with}~~
 \widetilde{G}_{i}^{t}(P_{i}^{t})=\frac{1}{M}(P_{i}^{t}-D^{t}), \notag \\
 & \quad P_{i}^{t}\in \Omega_{i}, ~\forall i \in \mathcal{M}, ~\forall t\in \{1,\cdots, T\}, \notag
\end{align}
where $\widetilde{\bm{P}}^{t}:=[\cdots,P_{i}^{t},\cdots]\in \mathbb{R}^{M d}$ concatenates all resource allocation strategies of all agents $ i \in \mathcal{M}$. Denote $\widetilde{\Omega}$ as the Cartesian product of $\Omega_{i}$ for all $i\in  \mathcal{M}$.

\textbf{Example} (Decentralized online economic dispatch): We consider a decentralized online economic dispatch problem involving $M$ generation stations, among which some are traditional thermal and the others are renewable wind generation stations.
The entire power network is modeled as an undirected, connected graph $\widetilde{\mathcal{G}} (\mathcal{M} = \mathcal{M}_{th} \cup \mathcal{M}_{wi}, \widetilde{\mathcal{E}})$, in which $\mathcal{M}_{th}$ and $\mathcal{M}_{wi}$ represent the set of thermal generation stations and the set of wind generation stations, respectively.
Each thermal generation station $i\in \mathcal{M}_{th}$ has a scheduled power generation strategy $P_{i}^{t}$ at time period $t$, which is confined by a local power capacity limit $\Omega_{i}=[P_{th}^{\min}, P_{th}^{\max}]$. Here we assume that all local power capacity limits of the thermal generation stations are the same for notational simplicity.
Considering the stability of the thermal power output, thermal generation station $i$ usually has a time-invariant quadratic cost function $C_{i}(P_{i}^{t})=\eta_{i}(P_{i}^{t})^{2}+\zeta_{i}P_{i}^{t}+\xi_{i}$, where $\eta_{i}$, $\zeta_{i}$ and $\xi_{i}$ are cost coefficients \cite{b-Zhenwei-Guo-2021,b-Jiahu-Qin-2019}.
Each wind generation station $j\in \mathcal{M}_{wi} $ has a scheduled power generation strategy $P_{j}^{t}$ at time period $t$, which is confined by a local power capacity limit $\Omega_{j}=[P_{wi}^{\min}, P_{wi}^{\max}]$.
However, since the power outputs of the wind generation stations are influenced by weather conditions, wind generation station $j$ has a time-varying cost function $C_{j}^{t}(\varsigma^{t},\phi^{t},P_{j}^{t})$, where $\varsigma^{t}$ and $\phi^{t}$ respectively denote the scale and shape factors of the Weibull distribution of the wind speed \cite{b-Fanghong-Guo-2021, b-Xian-Liu-2010, b-John-Hetzer-2008, b-Fang-Yao-2012}.
The average power demand $D^{t}$ varies over time.
In this paper, we consider a day-ahead power generation scheduling task with a time resolution of $5$ minutes, requiring $288$ dispatch decisions over a $24$-hour horizon.
The goal of the decentralized online economic dispatch problem is to determine an optimal sequence of power generation strategies $\{\{P_{i}^{t}\}_{\forall i \in \mathcal{M}_{th}},\{P_{j}^{t}\}_{\forall j \in \mathcal{M}_{wi}}\}_{t=1}^{288}$ over the time horizon $[1, 288]$ that minimizes the sum of time-varying generation costs while satisfying both long-term power demand and power capacity limits.
Hence, this online economic dispatch problem can be written as
\begin{align}
\label{online-decentralized-economic-dispatch-problem-example}
 & \hspace{-2em} \min_{\{\{P_{i}^{t}\},\{P_{j}^{t}\}\}_{t=1}^{288}}    \quad \sum_{t=1}^{288}\left[\sum_{i\in \mathcal{M}_{th}}C_{i}(P_{i}^{t})+\sum_{j\in \mathcal{M}_{wi}}C_{j}^{t}(P_{j}^{t})\right], \\
 s.t. & \quad \sum_{t=1}^{288}\frac{1}{M}\left[\sum_{i\in \mathcal{M}_{th}}P_{i}^{t}+\sum_{j\in \mathcal{M}_{wi}}P_{j}^{t}\right]-\sum_{t=1}^{288}D^{t}=0, \notag \\
 & \quad P_{i}^{t}\in [P_{th}^{\min},P_{th}^{\max}], ~\forall i \in \mathcal{M}_{th}, ~\forall t\in \{1,\cdots, 288\}, \notag\\
 & \quad P_{j}^{t}\in [P_{wi}^{\min},P_{wi}^{\max}], ~\forall j \in \mathcal{M}_{wi}, ~\forall t\in \{1,\cdots, 288\}.\notag
\end{align}

To solve \eqref{online-decentralized-economic-dispatch-problem} in a decentralized manner, agents communicate with their neighbors and exchange messages.
However, not all agents are reliable. Some of the agents are subject to external manipulations or internal damages, such that they behave maliciously and send wrong messages to neighboring agents, thereby disrupting the online resource allocation optimization process.
We refer to them as Byzantine agents, and the other ones as benign agents.

Denote the sets of the Byzantine and benign agents as $\mathcal{B}$ and $\mathcal{H}$, respectively.
Because the Byzantine agents may not adhere to the given optimization process, it is impossible to solve \eqref{online-decentralized-economic-dispatch-problem}.
Hence, when there are Byzantine agents, the oracle decentralized online resource allocation problem for the benign agents is refined to
\begin{align}
\label{online-decentralized-economic-dispatch-problem-oracle}
 \min_{\{\bm{P}^{t}\}_{t=1}^{T}}   & \quad \sum_{t=1}^{T}C^{t}(\bm{P}^{t})~~{\rm with}~~ C^{t}(\bm{P}^{t}) =\sum_{i\in \mathcal{H}} C_{i}^{t}(P_{i}^{t}), \\
 s.t. & \quad \sum_{t=1}^{T}\sum_{i\in \mathcal{H}}G_{i}^{t}(P_{i}^{t})=0 ~~{\rm with}~~G_{i}^{t}(P_{i}^{t})=\frac{1}{H}(P_{i}^{t}-D^{t}), \notag \\
 & \quad P_{i}^{t}\in \Omega_{i}, ~ \forall i \in \mathcal{H}, ~ \forall t\in \{1,\cdots, T\}, \notag
\end{align}
within which $H$ is the number of benign agents, $\bm{P}^{t}:=[\cdots,P_{i}^{t},\cdots]\in \mathbb{R}^{H d}$ concatenates all resource allocation strategies of all benign agents $ i \in \mathcal{H}$. Denote $\Omega$ as the Cartesian product of $\Omega_{i}$ for $i \in  \mathcal{H}$.

In this paper, we focus on developing Byzantine-resilient decentralized online resource allocation algorithms to tackle \eqref{online-decentralized-economic-dispatch-problem-oracle}, in the presence of Byzantine attacks.

\section{ATTACK-FREE DECENTRALIZED ONLINE RESOURCE ALLOCATION}\label{sec 3}
This section introduces a decentralized online resource allocation algorithm designed to tackle \eqref{online-decentralized-economic-dispatch-problem}, and also highlights its vulnerability to Byzantine attacks.

\subsection{Algorithm Development}
Since agents cannot access future time-varying costs and demands, \eqref{online-decentralized-economic-dispatch-problem} must be tackled in an online fashion.
The online regularized Lagrangian function corresponding to \eqref{online-decentralized-economic-dispatch-problem} at each time period $t$ is given by
\begin{align}
\label{Lagrangian-fucntion}
 \mathcal{L}^{t}_{\theta}(\widetilde{\bm{P}},\widetilde{\lambda})
= \widetilde{C}^{t}(\widetilde{\bm{P}})+ \langle \widetilde{\lambda}, \sum_{i\in \mathcal{M}}\widetilde{G}_{i}^{t}(P_{i})  \rangle -\frac{\theta}{2}\|\widetilde{\lambda}\|^{2},
\end{align}
where $\widetilde{\lambda} \in \mathbb{R}^{d}$ represents the dual variable and $\theta >0$ stands for a regularization parameter.
\begin{remark}
\label{regularization term}
Adding a regularization term to the dual variable is a classical technique in online primal-dual optimization for preventing the dual variable from becoming excessively large \cite{b-Jueyou-Li2022,b-Xinlei-Yi2020, b-Keishin-Tada2024, b-Xinlei-Yi2021a,b-Xinlei-Yi2023,b-Kunpeng-Zhang2024}. When the dual variable grows too large, it can significantly amplify the gradient of the Lagrangian with respect to the primal variable, potentially leading to unstable updates and poor regret performance.
To address this issue, we include the quadratic regularization term $-\frac{\theta}{2}\|\widetilde{\lambda}\|^2$ in the Lagrangian, which stabilizes the primal-dual updates and facilitates the convergence analysis.
\end{remark}

To find the saddle point of $\mathcal{L}^{t}_{\theta}(\widetilde{\bm{P}},\widetilde{\lambda})$ at time period $t$, the classical online primal-dual algorithm \cite{mahdavi2012trading,cao2018online} is
\begin{align}
\label{online-primal-dual-primal}
 P_{i}^{t+1}=&\arg\min_{P\in \Omega_{i}} \{\langle P-P_{i}^{t},\nabla C_{i}^{t}(P_{i}^{t}) + \frac{\widetilde{\lambda}^{t}}{M}  \rangle \notag\\
 &+\frac{1}{2 \alpha}\|P-P_{i}^{t}\|^{2}\},\\
\label{online-primal-dual-dual}
\widetilde{\lambda}^{t+1} =& \widetilde{\lambda}^{t}+\beta\cdot(\sum_{i\in \mathcal{M}}\widetilde{G}_{i}^{t}(P_{i}^{t})-\theta \widetilde{\lambda}^{t}),
\end{align}
where $\alpha>0$ and $\beta>0$ are step sizes. The above online primal-dual algorithm has been proven to attain sublinear static regret \cite{mahdavi2012trading} or dynamic regret \cite{cao2018online}, as well as sublinear accumulative constraint violation.

Nevertheless, \eqref{online-primal-dual-primal} and \eqref{online-primal-dual-dual} cannot be executed in a decentralized manner as the dual variable $\widetilde{\lambda}$ and the constraint function $\sum_{i\in \mathcal{M}}\widetilde{G}_{i}^{t}(P_{i}^{t})$ involve global information. To address this issue, we assign each agent $i$ a local dual variable $\lambda_{i}$ and approximate the global constraint function $\sum_{i\in \mathcal{M}}\widetilde{G}_{i}^{t}(P_{i}^{t})$ using $\widetilde{G}_{i}^{t}(P_{i}^{t})$. In addition, we let each agent $i$ aggregate its own local dual variable with those of its neighboring agents using a well-designed weight matrix to promote the consensus of the dual variables. Thus, we have
\begin{align}
\label{attack-free online decentralized resource allocation algorithm-primal}
P_{i}^{t+1}  =& \arg\min_{P\in \Omega_{i}} \{\langle P-P_{i}^{t},\nabla C_{i}^{t}(P_{i}^{t}) + \frac{\lambda_{i}^{t}}{M} \rangle \notag\\
&+\frac{1}{2 \alpha}\|P-P_{i}^{t}\|^{2}\},\\
\label{attack-free online decentralized resource allocation algorithm-dual}
\lambda_{i}^{t+\frac{1}{2}} =& \lambda_{i}^{t}+\beta\cdot(\widetilde{G}_{i}^{t}(P_{i}^{t})-\theta \lambda_{i}^{t}),\\
\label{attack-free online decentralized resource allocation algorithm-aggregation-dual}
\lambda_{i}^{t+1}=&\sum_{j\in \mathcal{N}_{i}\cup\{i\}}\widetilde{e}_{ij}\lambda_{j}^{t+\frac{1}{2}},
\end{align}
where $\widetilde{e}_{ij}$ is the weight assigned by agent $i$ to $j$. The weight matrix $\widetilde{E}:=[\widetilde{e}_{ij}]\in \mathbb{R}^{M\times M}$, in which $\widetilde{e}_{ij}>0$ if and only if $(i,j)\in \widetilde{\mathcal{E}}$ or $i=j$, is doubly stochastic.
The updates are summarized in Algorithm \ref{alg1}.

\begin{algorithm}[H]
\caption{Attack-free decentralized online resource allocation algorithm}\label{alg1}
\begin{algorithmic}
\STATE
\STATE $P_{i}^{0}=\lambda_{i}^{0}=D^{0}=0$ for all agents $i\in \mathcal{M}$.
\FOR{$t=0$ to $T$}
  \FOR{all agents $i\in \mathcal{M}$}
  \STATE Compute $P_{i}^{t+1}$ according to \eqref{attack-free online decentralized resource allocation algorithm-primal}.
  \STATE Compute $\lambda_{i}^{t+\frac{1}{2}}$ according to \eqref{attack-free online decentralized resource allocation algorithm-dual}.
  \STATE Broadcast $\lambda_{i}^{t+\frac{1}{2}}$ to its neighbors.
  \STATE Receive $\lambda_{j}^{t+\frac{1}{2}}$ from its neighbors.
  \STATE Compute $\lambda_{i}^{t+1} $ according to \eqref{attack-free online decentralized resource allocation algorithm-aggregation-dual}.
  \ENDFOR
\ENDFOR
\end{algorithmic}
\end{algorithm}

\subsection{Vulnerability of Algorithm \ref{alg1} under Byzantine Attacks}
Under suitable regularization conditions,
Algorithm \ref{alg1} is able to tackle \eqref{online-decentralized-economic-dispatch-problem} if all agents are benign \cite{b-Jueyou-Li2022, b-Xinlei-Yi2020}, and achieves sublinear dynamic regret and accumulative constraint violation.
However, in the presence of Byzantine attacks, such convergence guarantees no longer hold.
At each time period $t$, when agent $j$ is benign, it sends the true $\lambda_{j}^{t+\frac{1}{2}}$ to its neighbors. But, a Byzantine agent $j$ can instead send a malicious message $\dag$ to its neighbors.\footnote{In fact, it can send different wrong messages to different neighbors. We use the same $\dag$ for convenience.} Define the message sent by agent $j$ as
\begin{align}
\label{broadcasting_message}
\check{\lambda}_{j}^{t+\frac{1}{2}} =\left\{\begin{matrix}\lambda_{j}^{t+\frac{1}{2}}, \quad & j\in \mathcal{H},
 \\ \dag, \quad & j \in \mathcal{B}.
\end{matrix}\right.
\end{align}
Under \eqref{broadcasting_message}, we note that the weighted average aggregation in \eqref{attack-free online decentralized resource allocation algorithm-aggregation-dual} is wrong and controlled by malicious messages from Byzantine agents, in the sense that it incorporates arbitrary messages from Byzantine agents, which can significantly distort the result and make the aggregation deviate from the true weighted average of benign dual variables.
This yields unfavorable resource allocation strategies for the benign agents.

\section{BYZANTINE-RESILIENT DECENTRALIZED ONLINE RESOURCE ALLOCATION}
\label{sec 4}
Given that the vulnerability of Algorithm \ref{alg1} stems from its susceptible weighted average aggregation rule in the form of $\lambda_{i}^{t+1} =\sum_{j\in \mathcal{N}_{i}\cup\{i\}}\widetilde{e}_{ij}\lambda_{j}^{t+\frac{1}{2}}$, a natural idea to address this issue is to replace them with a robust aggregation rule. To this end, we consider a class of robust aggregation rules, denoted as $AGG(\cdot)$, and introduce a set of properties that such rules satisfy to support the convergence analysis of our online resource allocation algorithms.

\noindent \textbf{Properties of Robust Aggregation Rules in Online Resource Allocation.} Intuitively, for benign agent $i$, we expect the output of $AGG(\lambda_{i}^{t+\frac{1}{2}}, \{\check{\lambda}_{j}^{t+\frac{1}{2}} \}_{j\in \mathcal{N}_{i}})$ to be sufficiently close to a proper weighted average of the messages from its benign neighbors and its own local dual variable. Below, we use $\bar{\lambda}_{i}^{t+\frac{1}{2}}:=\sum_{j\in (\mathcal{N}_{i}\cap\mathcal{H})\cup {i}}e_{ij}\lambda_{j}^{t+\frac{1}{2}}$, in which the weights $\{e_{ij}\}_{j \in \mathcal{H}}$ satisfy $\sum_{j\in (\mathcal{N}_{i}\cap\mathcal{H})\cup {i}}e_{ij}=1$, to denote such a weighted average. We also use the value of $\sum_{j\in \mathcal{N}_{i}\cap\mathcal{H}\cup\{i\}}$ $e_{ij}\|\lambda_{j}^{t+\frac{1}{2}}-\bar{\lambda}_{i}^{t+\frac{1}{2}}\|^{2} $ as the standard to measure the proximity. Therefore, a set of robust aggregation rules should satisfy the following property.

\begin{property}
\label{d1}
Consider an robust aggregation rule $AGG(\cdot)$.
For any set $\{\lambda_{i},\{\check{\lambda}_{j}\}_{j\in \mathcal{N}_{i}}\}$, there exists a constant $\rho \geq 0$ and a matrix $E\in \mathbb{R}^{H \times H}$ whose elements satisfy $e_{ij} \in (0,1]$ when $j\in (\mathcal{N}_{i}\cap\mathcal{H})\cup {i}$, $e_{ij}=0$ when $j\notin (\mathcal{N}_{i}\cap\mathcal{H})\cup {i}$, and $\sum_{j\in (\mathcal{N}_{i}\cap\mathcal{H})\cup {i}}e_{ij}=1$ for any $i \in \mathcal{H}$, such that it holds
\begin{align}
    \|AGG(\lambda_{i},\{\check{\lambda}_{j}\}_{j\in \mathcal{N}_{i}})-\bar{\lambda}_{i}\|^{2} \le  \rho \sum_{j\in \mathcal{N}_{i}\cap\mathcal{H} \cup \{i\}}e_{ij}\|\lambda_{j}-\bar{\lambda}_{i}\|^{2}, \nonumber
\end{align}
for any $i\in \mathcal{H}$, with $\bar{\lambda}_{i}:=\sum_{j\in (\mathcal{N}_{i}\cap\mathcal{H})\cup {i}}e_{ij}\lambda_{j}$. Here, $\rho$ is the contraction constant and $E$ is the weight matrix associated with the robust aggregation rule $AGG(\cdot)$.
\end{property}
\begin{remark}
Property \ref{d1} in this paper is similar to the corresponding property used in \cite{b-ZhaoxianWu-2023,b-RunhuaWang-2023,b-Haoxiang-Ye-2023, b-Haoxiang-Ye-2025}.
Specifically, in Property \ref{d1}, we use the value of $\sum_{j\in \mathcal{N}_{i}\cap\mathcal{H}\cup\{i\}}e_{ij}\|\lambda_{j}-\bar{\lambda}_{i}\|^{2} $ as the standard to measure the proximity.
The works of \cite{b-ZhaoxianWu-2023,b-RunhuaWang-2023,b-Haoxiang-Ye-2023, b-Haoxiang-Ye-2025} use the value of $\max_{j\in \mathcal{N}_{i}\cap\mathcal{H}\cup\{i\}}\|\lambda_{j}-\bar{\lambda}_{i}\|^{2} $.
This adjustment is made to facilitate the convergence analysis.
\end{remark}

Simply satisfying Property \ref{d1} is insufficient to guarantee the convergence of an online resource allocation algorithm using $AGG(\cdot)$. The reason is that, bounding the benign dual variables is of paramount importance in the investigated online resource allocation problem, but with Property \ref{d1} the benign dual variables may grow to infinity at a rate of $(1+2\sqrt{\rho})^{t}$. Our analysis reveals that bounding the benign dual variables requires the output of $AGG(\lambda_{i}^{t+\frac{1}{2}}, \{\check{\lambda}_{j}^{t+\frac{1}{2}} \}_{j\in \mathcal{N}_{i}})$ to be bounded by the maximal norm of all benign neighboring dual variables, as shown in the following property.
\begin{property}
\label{d2}
Consider a robust aggregation rule $AGG(\cdot)$. For any set $\{\lambda_{i},\{\check{\lambda}_{j}\}_{j\in \mathcal{N}_{i}}\}$, it holds
\begin{align}
\|AGG(\lambda_{i},\{\check{\lambda}_{j}\}_{j\in \mathcal{N}_{i}})\| \le  \max_{j\in \mathcal{N}_{i}\cap \mathcal{H}\cup \{i\}}\|\lambda_{j}\|.\nonumber
\end{align}
\end{property}


\textbf{Robust Aggregation Rules Satisfying Both Properties \ref{d1} and \ref{d2}.} For offline resource allocation, there exist various robust aggregation rules, such as $CTM$, $IOS$ and $SCC$. These existing robust aggregation rules all satisfy Property \ref{d1} \cite{b-Haoxiang-Ye-2023}. As for Property \ref{d2}, the scalar case of $CTM$ is a special exception: when the dual variable is scalar, $CTM$ reduces to $TM$, and $TM$ satisfies both Properties \ref{d1} and \ref{d2}, with the proof of Property \ref{d2} given in Appendix A-D.
In contrast, all the remaining cases fail to satisfy Property \ref{d2}$,$ as shown by the explicit counter-examples in Appendix A-B.
Fortunately, we find that combining the existing robust aggregation rules $IOS$ and $SCC$ with the adaptive robust clipping technique (denoted as $ARC$) proposed in \cite{b-Youssef-Allouah-2024} yields two robust aggregation rules, namely $IOS(ARC(\cdot))$ and $SCC(ARC(\cdot))$, that satisfy both Properties \ref{d1} and \ref{d2}$.$ In Appendices A-C and A-D, we prove that $IOS(ARC(\cdot))$ and $SCC(ARC(\cdot))$ satisfy both Properties \ref{d1} and \ref{d2} in arbitrary dimensions.
However, the $ARC$ technique is still insufficient to guarantee that a $CTM$-based design satisfies Property \ref{d2} in higher dimensions. Due to the coordinate-wise trimming nature of $CTM$, how to design a $CTM$-based mechanism that also satisfies Property \ref{d2} for vector-valued dual variables remains an open problem.

Thus, the robust aggregation rules provably applicable to online resource allocation depend on the dual dimension. When the dual variable is scalar, the applicable choices are $TM(\cdot)$, $IOS(ARC(\cdot))$, and $SCC(ARC(\cdot))$. When the dual dimension satisfies $d\ge 2$, the applicable choices are $IOS(ARC(\cdot))$ and $SCC(ARC(\cdot))$. In the following, we describe $ARC(\cdot)$, $TM(\cdot)$, $IOS(ARC(\cdot))$, and $SCC(ARC(\cdot))$ in turn.

\noindent \textbf{Description of $ARC$}.
For any benign agent $i \in \mathcal{H}$, the $ARC$ procedure consists of the following steps:\\
\noindent \textbf{Step 1 (Sorting)}: Benign agent $i$ receives dual variables $\{\check{\lambda}_j\}_{j \in \mathcal{N}_i}$ from its neighbors and sorts them by their norms to get a permutation $\pi$ such that $\|\check{\lambda}_{\pi_{1}}\|\ge\|\check{\lambda}_{\pi_{2}}\| \cdots \ge \|\check{\lambda}_{\pi_{|\mathcal{N}_{i}|}}\|$.\\
\noindent \textbf{Step 2 (Clipping threshold selection)}: Given an upper bound $b_i$ on the number of Byzantine neighbors, benign agent $i$ selects the $(b_i + 1)$-th largest norm as the clipping threshold, as $C_{i}=\|\check{\lambda}_{\pi_{b_{i}+1}}\|$.\\
\noindent \textbf{Step 3 (Clipping)}: Agent $i$ clips each received dual variable $\check{\lambda}_j (j\in \mathcal{N}_i)$ to obtain $clip_{C_{i}}(\check{\lambda}_{j}):=\min (1, \frac{C_{i}}{\|\check{\lambda}_{j}\|}) \check{\lambda}_{j}$.

Based on the above three steps, we conclude that the norm of any clipped dual variable in the set $\{clip_{C_{i}}(\check{\lambda}_{j})\}_{j\in \mathcal{N}_{i}}$ must be smaller than the maximal norm of all benign dual variables in $\{\lambda_{j}\}_{j\in \mathcal{N}_{i}\cap \mathcal{H}}$, i.e., $\|clip_{C_{i}}(\check{\lambda}_{j})\|\le \max_{j \in \mathcal{N}_{i}\cap \mathcal{H}}\|\lambda_{j}\|, \forall j \in \mathcal{N}_{i}$.
Specifically, by Step 2, the clipping threshold $C_i$ is chosen as the $(b_i + 1)$-th largest norm among all received dual variables. Since there are at most $b_i$ Byzantine neighbors, there must exist at least one benign neighbor $j \in \mathcal{N}_i \cap \mathcal{H}$ such that $\|\lambda_j\| \ge C_i$. Then, by Step 3, each received dual variable is clipped to have norm at most $C_i$, which implies the desired bound. Finally, we denote $ARC(\lambda_{i},\{\check{\lambda}_{j}\}_{j\in \mathcal{N}_{i}})=\{\lambda_{i},\{clip_{C_{i}}(\check{\lambda}_{j})\}_{j\in \mathcal{N}_{i}}\}$.\\
\textbf{Description of $TM(\cdot)$}.
For any benign agent $i \in \mathcal{H}$, the $TM$ procedure operates on the scalar dual variables $\{\check{\lambda}_j\}_{j \in \mathcal{N}_i}$, and consists of the following steps:\\
\noindent \textbf{Step 1 (Sorting)}: Benign agent $i$ sorts the received neighboring dual variables $\{\check{\lambda}_j\}_{j \in \mathcal{N}_i}$ in ascending order.\\
\noindent \textbf{Step 2 (Outlier removal)}: Given an upper bound $b_i$ on the number of Byzantine neighbors, benign agent $i$ discards the largest $b_i$ and smallest $b_i$ values.\\
\noindent \textbf{Step 3 (Averaging)}: Benign agent $i$ computes the average of the remaining values together with its own local dual variable $\lambda_i$ to obtain the aggregation result.\\
\noindent \textbf{Description of $IOS(ARC(\cdot))$}.
For any benign agent $i \in \mathcal{H}$, the $IOS$ procedure operates on the clipped dual variables $\{clip_{C_i}(\check{\lambda}_j)\}_{j \in \mathcal{N}_i}$ produced by $ARC$ procedure, and consists of the following steps:\\
\noindent \textbf{Step 1 (Iterative outlier removal)}: Agent $i$ iteratively removes $b_i$ outliers from the set of received clipped dual variables. In each iteration, it computes the weighted average of the current set, identifies the variable farthest from this weighted average and removes it.\\
\noindent \textbf{Step 2 (Weighted averaging)}: Agent $i$ computes the weighted average of the remaining variables by re-normalizing their weights to sum to one, and returns the result as the aggregation output.\\
\noindent \textbf{Description of $SCC(ARC(\cdot))$}. For any benign agent $i \in \mathcal{H}$, the $SCC$ procedure operates on the clipped dual variables $\{clip_{C_i}(\check{\lambda}_j)\}_{j \in \mathcal{N}_i}$ produced by $ARC$, and consists of the following steps:\\
\noindent \textbf{Step 1 (Clipping)}: Agent $i$ selects a local clipping threshold $\tau_{i}$ and, for each received clipped dual variable, checks its distance to its own dual variable. If the distance exceeds $\tau_{i}$, the variable is clipped toward the agent's own value along the same direction, such that the resulting distance equals the clipping threshold. Otherwise, the variable is kept unchanged.\\
\noindent \textbf{Step 2 (Weighted averaging)}: Agent $i$ computes a weighted average over the resulting clipped variables and its own dual variable, and returns the result as the aggregation output.

Based on the applicable robust aggregation rules discussed above, we propose a class of Byzantine-resilient online resource allocation algorithms.
At time period $t$, the updates of the primal and dual variables for each benign agent $i \in \mathcal{H}$ are given by
\begin{align}
\label{Byzantine-resilient online decentralized resource allocation algorithm-primal}
P_{i}^{t+1}=&\arg\min_{P\in \Omega_{i}}\{ \langle P-P_{i}^{t},\nabla C_{i}^{t}(P_{i}^{t}) + \frac{\lambda_{i}^{t}}{M} \rangle \notag\\
&\hspace{4em}+\frac{1}{2 \alpha}\|P-P_{i}^{t}\|^{2} \},\\
\label{Byzantine-resilient online decentralized resource allocation algorithm-dual}
\lambda_{i}^{t+\frac{1}{2}} &= \lambda_{i}^{t}+\beta\cdot(\widetilde{G}_{i}^{t}(P_{i}^{t})-\theta \lambda_{i}^{t}),\\
\label{Byzantine-resilient online decentralized resource allocation algorithm-aggregation-dual}
\lambda_{i}^{t+1}& =AGG(\lambda_{i}^{t+\frac{1}{2}}, \{\check{\lambda}_{j}^{t+\frac{1}{2}}\}_{j\in \mathcal{N}_{i}}),
\end{align}
where $AGG(\cdot)$ denotes an applicable robust aggregation rule for the present online problem: when the dual variable is scalar, $AGG(\cdot)$ can be chosen as $TM(\cdot)$, $IOS(ARC(\cdot))$, or $SCC(ARC(\cdot))$; when the dual dimension satisfies $d\ge 2$, $AGG(\cdot)$ can be chosen as $IOS(ARC(\cdot))$ or $SCC(ARC(\cdot))$.
The updates are summarized in Algorithm \ref{alg2}.
\begin{remark}
    \label{r-OFF(ARC)} The idea of combining an existing robust aggregation rule with $ARC$ technique has also been explored in~\cite{b-Youssef-Allouah-2024}. However, \cite{b-Youssef-Allouah-2024} focuses on Byzantine-resilient consensus optimization coordinated by a central server. In contrast, our work considers Byzantine-resilient online resource allocation that is non-consensual and decentralized.
\end{remark}
\begin{algorithm}[H]
\caption{Byzantine-resilient decentralized online resource allocation algorithm}\label{alg2}
\begin{algorithmic}
\STATE
\STATE $P_{i}^{0}=\lambda_{i}^{0}=D^{0}=0$ for all benign agents $i\in \mathcal{H}$.
\FOR{$t=0$ to $T$}
  \FOR{all benign agents $i\in \mathcal{H}$}
  \STATE Compute $P_{i}^{t+1}$ according to \eqref{Byzantine-resilient online decentralized resource allocation algorithm-primal}.
  \STATE Compute $\lambda_{i}^{t+\frac{1}{2}}$ according to \eqref{Byzantine-resilient online decentralized resource allocation algorithm-dual}.
  \STATE Broadcast $\lambda_{i}^{t+\frac{1}{2}}$ to its neighbors.
  \STATE Receive $\check{\lambda}_{j}^{t+\frac{1}{2}}$ from its neighbors.
  \STATE Compute $\lambda_{i}^{t+1}$ according to \eqref{Byzantine-resilient online decentralized resource allocation algorithm-aggregation-dual}.
  \ENDFOR
  \FOR{all Byzantine agents $i\in \mathcal{B}$}
  \STATE Broadcast $\check{\lambda}_{i}^{t+\frac{1}{2}}=\dag$ to its neighbors
  \ENDFOR
\ENDFOR
\end{algorithmic}
\end{algorithm}

\section{THEORETICAL ANALYSIS}\label{sec 5}
This section analyzes the performance of attack-free and Byzantine-resilient decentralized online resource allocation algorithms, respectively. We begin with several assumptions.
\begin{assumption}
\label{a1}
For any agent $i\in \mathcal{M}$ at each time period $t$, the local cost function $C_{i}^{t}(\cdot)$ is convex and bounded,
and the local constraint set $\Omega_{i}$ is compact and convex. Specifically, there exist positive constants $F$ and $R$ such that $|C_{i}^{t}(\cdot)|\le F$ and $|x-y|\le R$ for all $x, y \in \Omega_{i}$. The gradient of $C_{i}^{t}(\cdot)$ is bounded. Namely, there exists a positive constant $\varphi$ such that $|\nabla C_{i}^{t}(\cdot) |\le \varphi$.
Furthermore, for any agent $i\in \mathcal{M}$ at each time period $t$, the local constraints $\widetilde{G}_i^t(\cdot)$ and $G_i^t(\cdot)$ are both bounded, i.e., $|\widetilde{G}_i^t(\cdot)|\le \widetilde{\psi}$ and $|G_i^t(\cdot)|\le \psi$, where $\widetilde{\psi}$ and $\psi$ are positive constants.
\end{assumption}

Assumption \ref{a1} is common in analyzing the convergence of online primal-dual algorithms \cite{b-Jueyou-Li2022,b-Xinlei-Yi2020, b-Keishin-Tada2024, b-Xinlei-Yi2021a,b-Xinlei-Yi2023,b-Kunpeng-Zhang2024}.

\begin{assumption}
\label{a3}
Consider a subgraph $\mathcal{G} (\mathcal{H},\mathcal{E})$ of $\widetilde{\mathcal{G}} (\mathcal{M},\widetilde{\mathcal{E}})$, where $\mathcal{E}$ is the set of edges between the benign agents. Both graphs $\widetilde{\mathcal{G}} (\mathcal{M},\widetilde{\mathcal{E}})$ and $\mathcal{G} (\mathcal{H},\mathcal{E})$ are undirected and connected. The weight matrices $\widetilde{E}$ and $E$ are doubly stochastic and row stochastic, respectively, and also satisfy
\begin{align}
\widetilde{\kappa} & := \|\widetilde{E}-\frac{1}{M}\widetilde{\bm{1}}\widetilde{\bm{1}}^{\top}\|^{2} < 1, \label{eq_a6_tttt} \\
\kappa & := \|E-\frac{1}{H}\bm{1}\bm{1}^{\top} E\|^{2} < 1, \label{eq_a6}
\end{align}
in which $\widetilde{\bm{1}}\in \mathbb{R}^{M}$ and $\bm{1}\in \mathbb{R}^{H}$ are both all-one column vectors.
\end{assumption}

Assumption \ref{a3} describes the connectivity of the communication topology. Similar assumptions have been widely adopted in prior works on Byzantine-resilient decentralized optimization \cite{b-ZhaoxianWu-2023,b-RunhuaWang-2023, b-Haoxiang-Ye-2023}.

\subsection{Attack-free Decentralized Online Resource Allocation Algorithm}
We shall use two commonly used performance metrics for online constrained optimization:
(i) dynamic regret: $$\widetilde{R}^{T}_{\mathcal{M}}:=\sum_{t=1}^{T}\sum_{i\in \mathcal{M}}C_{i}^{t}(P_{i}^{t})-\sum_{t=1}^{T}\sum_{i\in \mathcal{M}}C_{i}^{t}(\widetilde{P}_{i}^{t*}),$$ where $\widetilde{P}_{i}^{t*}$ is the $i$th element of $\widetilde{\bm{P}}^{t*}:=\arg\min_{\widetilde{\bm{P}}\in \widetilde{\Omega}} \sum_{i\in \mathcal{M}}$ $C_{i}^{t}(P_{i})$, $s. t. ~ \sum_{i\in \mathcal{M}}\widetilde{G}_{i}^{t}(P_{i})=0$, the instantaneous optimal solution to \eqref{online-decentralized-economic-dispatch-problem} at time period $t$; (ii) accumulative constraint violation: $$\widetilde{V}^{T}_{\mathcal{M}}:=\|\sum_{t=1}^{T}\sum_{i\in \mathcal{M}}\widetilde{G}_{i}^{t}(P_{i}^{t})\|.$$

\begin{theorem}
\label{t1}
Suppose that Assumptions \ref{a1}--\ref{a3} hold and that the instantaneous optimal solutions to \eqref{online-decentralized-economic-dispatch-problem} satisfy $\sum_{t=1}^{T} $ $\sum_{i\in \mathcal{M}}$ $\|\widetilde{P}_{i}^{t*}-\widetilde{P}_{i}^{t-1*}\|=O(T^{\gamma})$, where $\gamma\in [0,1)$.
Set the step sizes $\alpha$, $\beta$ and the regularization parameter $\theta$ as
$\alpha=$ $T^{\frac{\gamma-1}{2}}$, $\beta=T^{-\frac{1}{2}}$ and $\theta=T^{-c}$, where $c\in (0,\frac{1-\gamma}{4})$. For the sequences $\{P_{i}^{t+1}\}_{i\in \mathcal{M}}$ generated by Algorithm \ref{a1}, we have
\begin{align}
\label{t1-1}
\widetilde{R}^{T}_{\mathcal{M}}& \le O(T^{\frac{1+\gamma}{2}+2c}),\\
\label{t1-2}
\widetilde{V}^{T}_{\mathcal{M}} & \le O(T^{\max\{1-\frac{c}{2},\frac{3+\gamma}{4}+\frac{c}{2}\}}).
\end{align}
\end{theorem}

\begin{remark}
\label{r1-t1}
Theorem \ref{t1} demonstrates that the attack-free decentralized online resource allocation algorithm achieves sublinear dynamic regret and accumulative constraint violation, aligning with the existing results for online convex optimization \cite{b-Jueyou-Li2022, b-Xinlei-Yi2020}. 
To achieve sublinear dynamic regret and accumulative constraint violation, it is required that the optimal solutions do not change too rapidly over time (i.e., the accumulative variation grows no faster than $O(T^\gamma)$ with $\gamma<1$). Additionally, the algorithm's step sizes and regularization parameter must be carefully chosen in accordance with this variation rate so that the regret exponent $\frac{1+\gamma}{2} + 2c$ and $\max\{1-\frac{c}{2},\frac{3+\gamma}{4}+\frac{c}{2}\}$ remains below 1. The
proof of Theorem \ref{t1} is in Appendix \ref{B}.
\end{remark}

\subsection{Byzantine-resilient Decentralized Online Resource Allocation Algorithm}
To evaluate the Byzantine-resilient decentralized online resource allocation algorithms, the performance metrics are modified to:
(i) dynamic regret: $${R}^{T}_{\mathcal{H}}:=\sum_{t=1}^{T}\sum_{i\in \mathcal{H}}C_{i}^{t}(P_{i}^{t})-\sum_{t=1}^{T}\sum_{i\in \mathcal{H}}C_{i}^{t}({P}_{i}^{t*}),$$ where ${P}_{i}^{t*}$ is the $i$th element of ${\bm{P}}^{t*}:=\arg\min_{{\bm{P}}\in {\Omega}} \sum_{i\in \mathcal{H}}$ $C_{i}^{t}(P_{i})$, $s. t. ~ \sum_{i\in \mathcal{H}}{G}_{i}^{t}(P_{i})=0$, namely the instantaneous optimal solution to \eqref{online-decentralized-economic-dispatch-problem-oracle} at time period $t$; (ii) accumulative constraint violation: $${V}^{T}_{\mathcal{H}}:=\|\sum_{t=1}^{T}\sum_{i\in \mathcal{H}} {G}_{i}^{t}(P_{i}^{t})\|.$$

\begin{theorem}
\label{t2}
Suppose that Assumptions \ref{a1}--\ref{a3} hold and that the instantaneous optimal solutions to \eqref{online-decentralized-economic-dispatch-problem-oracle} satisfy $\sum_{t=1}^{T} \sum_{i\in \mathcal{H}}$ $\|P_{i}^{t*}-P_{i}^{t-1*}\|=O(T^{\gamma})$, where $\gamma\in [0,1)$. Set the step sizes $\alpha$ and $\beta$ as
$\alpha=T^{\frac{\gamma-1}{2}}$ and $\beta=T^{-\frac{1}{2}}$. If the robust aggregation rule $AGG(\cdot)$ satisfies Properties \ref{d1} and \ref{d2} and the contraction constant satisfies $\rho \le \frac{(1-\kappa)^{2}}{64H}$, then for the sequences $\{P_{i}^{t+1}\}_{i\in \mathcal{H}}$ generated by Algorithm \ref{alg2}, we have
\begin{align}
\label{t2-1}
R^{T}_{\mathcal{H}}&\le O((\rho+\chi)\cdot \frac{T}{\theta}+\frac{T^{\frac{1+\gamma}{2}}}{\theta^{2}}),\\
\label{t2-2}
V^{T}_{\mathcal{H}}
&  \le O((\rho+\chi)T+\sqrt{\theta}\cdot T+\frac{T^{\frac{3+\gamma}{4}}}{\sqrt{\theta}}),
\end{align}
\end{theorem}
where $\chi=\frac{1}{H}\|E^{\top}\bm{1}-\bm{1}\|^{2}$ quantifies the skewness of the weight matrix $E$ associated with the online robust aggregation rules.

\begin{remark}
\label{r2-t2}
Choosing a regularization parameter $\theta=\frac{1}{2HF}$ leads to linear dynamic regret and accumulative constraint violation, as implied by Theorem \ref{t2}.
These linear bounds, although not promising, are inevitable and tight results in the presence of Byzantine attacks.
The underlying reason is that, the heterogeneity of agents' cost functions and the presence of Byzantine agents cannot guarantee perfectly accurate aggregation. As a result, non-vanishing aggregation errors accumulate over time, leading to linear dynamic regret and accumulative constraint violation.
Several recent studies support this viewpoint.
In particular, \cite{shi2025optimal} considers offline consensus optimization under Byzantine attacks and proves that, in the presence of data heterogeneity, any Byzantine-resilient first-order algorithm suffers from an unavoidable convergence error, for which a tight lower bound is established.
The work of \cite{dong2023byzantine} investigates online consensus optimization and shows that a range of Byzantine-resilient algorithms necessarily incur tight linear regret.
Although we consider a non-consensus resource allocation problem, the Byzantine-resilient operations in our algorithm are applied to the consensus dual variables. Consequently, the unavoidable aggregation errors in the dual space are propagated to the primal variables, ultimately leading to the observed linear regret and constraint violation.
In fact, for the offline setup, it has been proved in \cite{b-RunhuaWang-2023} that a class of Byzantine-resilient decentralized resource allocation algorithms converge to neighborhoods of the optimal resource allocation strategy, and the errors are in the order of $O(\rho+\chi)$.
Intuitively, for the online setup, such errors accumulate over the time horizon $[1,T]$, resulting in linear dynamic regret and accumulative constraint violation.

While this result may seem pessimistic, it is important to emphasize that it stems from the theoretical analysis that necessarily considers the worst-case scenario---specifically, when the heterogeneity of agents' cost functions and the presence of Byzantine agents make perfectly accurate aggregation impossible to guarantee. Nevertheless, in practical scenarios where the cost heterogeneity is low and the wrong messages are only outliers, perfectly accurate aggregation can often be achieved. In such cases, the theoretical parameter reduces to $\rho=0$ and $E$ becomes doubly stochastic (namely, $\chi =0$), so that Theorem \ref{t2} reduces to Theorem \ref{t1}, leading to sublinear dynamic regret and accumulative constraint violation.
\end{remark}

\begin{figure}[htbp]
\centerline{\includegraphics[width=9cm]{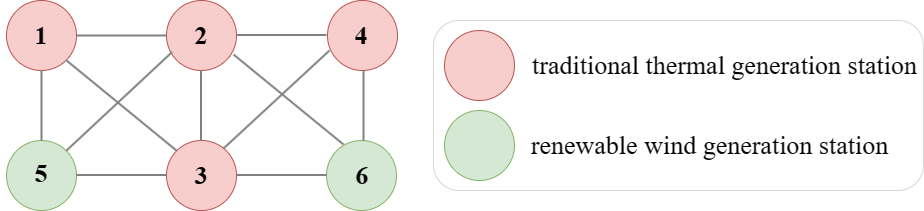}}
\caption{The communication graph of synthetic problem.}
\label{fig_1}
\end{figure}

\begin{figure*}[htbp]
\centerline{\includegraphics[width=17cm]{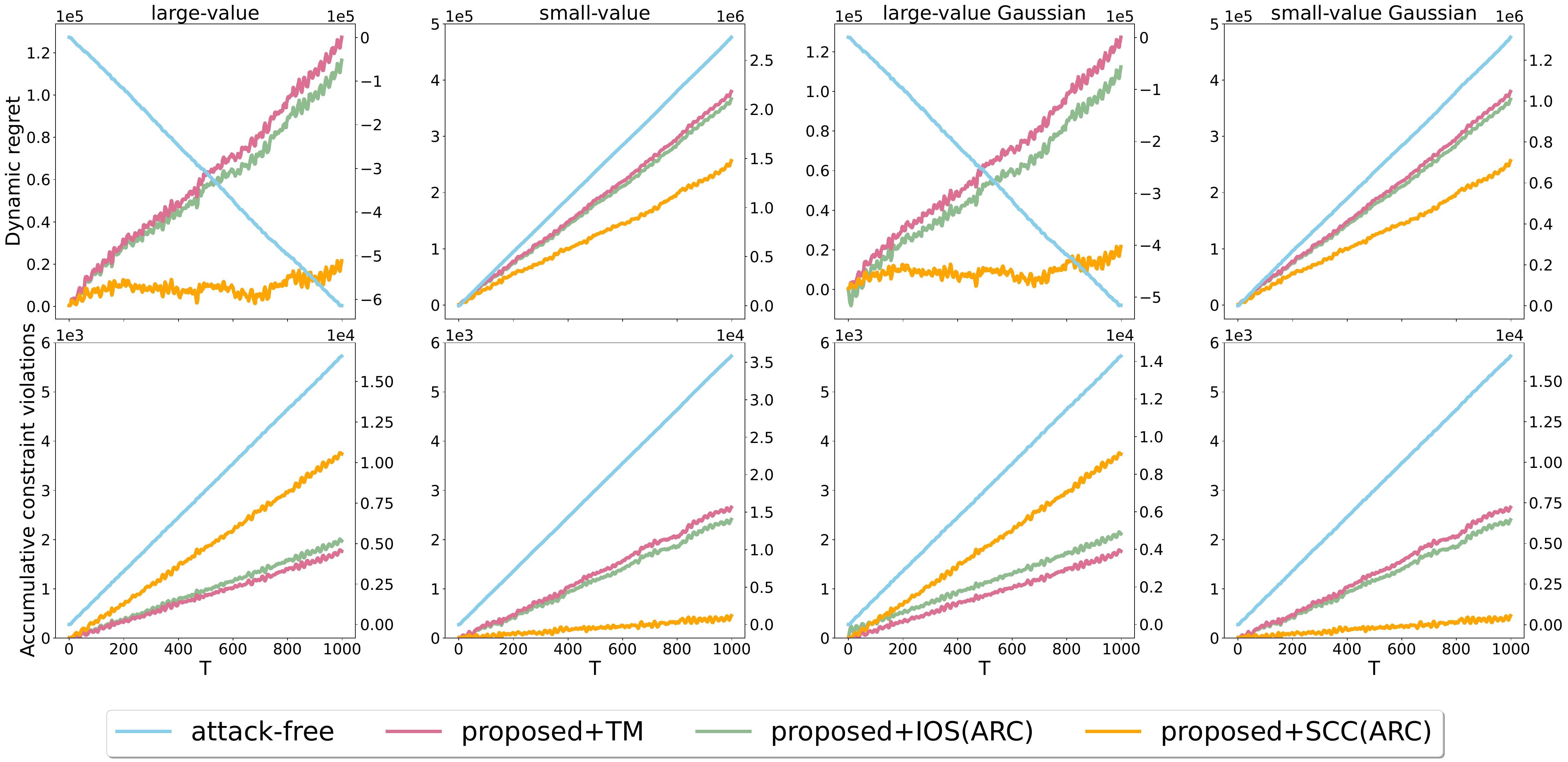}}
\caption{Dynamic regret and accumulative constraint violations of the compared algorithms under various Byzantine attacks.}
\label{fig_2}
\end{figure*}

\begin{table*}[!ht]
\centering
\renewcommand{\arraystretch}{1.18}
\caption{The parameters of traditional thermal generation stations for case 1 \cite{b-Fanghong-Guo-2021}}
\begin{tabular}{|c|c|c|c|c|c|} \hline Thermal generation station No. & $\eta_{i}$ & $\zeta_{i}$ & $\xi_{i}$ & $P_{th, i}^{min}$ & $P_{th, i}^{max}$\\ \hline
1  & 0.0675 &	2 & 0 & 50 & 200 \\ \hline
2  & 0.0675 & 1.75 & 0 & 20 & 120\\ \hline
3  & 0.0925 & 1 & 0 & 15 & 80 \\ \hline
4  & 0.0625 & 3 & 0 & 10 & 100\\ \hline
\end{tabular}\\
\label{table-1}
\end{table*}

\begin{table*}[!ht]
\centering
\renewcommand{\arraystretch}{1.18}
\caption{The parameters of renewable wind generation stations for case 1 \cite{b-Fanghong-Guo-2021, b-Fang-Yao-2012} }
\begin{tabular}{|c|c|c|c|c|c|c|c|c|c|} \hline
Wind generation station No. & $\varrho_{i}$ & $v_{in,i}$ & $v_{out,i}$ & $v_{r,i}$ & $\sigma_{ue,i}$ & $\sigma_{oe,i}$ & $P_{r,i}$ & $P_{wi,i}^{min}$ & $P_{wi,i}^{max}$\\ \hline
 5  &	 1  &	3 & 25 & 13 & 5 & 30& 160& 0 & 160\\ \hline
6  & 6 & 5 & 45 & 15 & 5 & 20 & 160 & 0 &160 \\ \hline
\end{tabular}\\
\label{table-2}
\end{table*}

\section{NUMERICAL EXPERIMENTS}\label{sec 6}
In this section,  we conduct numerical experiments on decentralized online economic dispatch  to validate the effectiveness of our proposed algorithms. The code is available online. \renewcommand{\thefootnote}{2}\footnote{\url{https://github.com/RunhuaWang}}
\subsection{Case 1: Synthetic Problem}
We first consider a power system comprising $4$ traditional thermal and $2$ renewable wind power stations.
The communication graph of the power system is shown in Fig. \ref{fig_1}. Based on the communication graph, we use the Metropolis constant weight rule \cite{b-Wei-Shi-2015} to generate a doubly stochastic weight matrix $\widetilde{E}$.
Each traditional thermal power station $i\in \{1,2,3,4\}$ possesses a cost function $C_{i}(P_{i})=\eta_{i}(P_{i})^{2}+\zeta_{i}P_{i}+\xi_{i}$ \cite{b-Zhenwei-Guo-2021, b-Jiahu-Qin-2019}, and is subject to the local constraint $[P_{th,i}^{\min},P_{th,i}^{\max}]$, where $\eta_{i}$, $\zeta_{i}$ and $\xi_{i}$ are cost coefficients, while $P_{th,i}^{\min}$ and $P_{th,i}^{\max}$ represent the low and upper bounds of power output, respectively.
The settings of these parameters are outlined in TABLE \ref{table-1}.
Each renewable wind power station $i\in \{5,6\}$ is governed by a time-varying cost function in the form of $C_{i}^{t}(P_{i})$ $=\varrho_{i} P_{i}+C_{ue,i}^{t}(\varsigma^{t},\phi^{t},\sigma_{ue,i},v_{in,i},v_{out,i},v_{r,i},P_{r,i}, P_{i})+ C_{oe,i}^{t}(\varsigma^{t}, \phi^{t},\sigma_{oe,i},v_{in,i},v_{out,i},v_{r,i},P_{r,i}, P_{i})$, where $\varrho_{i}$ denotes the cost coefficient, and $C_{ue,i}^{t}(\cdot)$ and $C_{oe,i}^{t}(\cdot)$ represent the underestimation and overestimation costs, $\sigma_{ue,i}$ and $\sigma_{oe,i}$ are underestimation and overestimation penalty cost coefficients, $v_{in,i},v_{out,i}$ and $v_{r,i}$ are cut-in, cut-out and rate wind speeds, while $P_{r,i}$ is the rate power output. The specific forms of cost function $C_{ue,i}^{t}(\cdot)$ and $C_{oe,i}^{t}(\cdot)$ can be found in \cite{b-Fanghong-Guo-2021,b-Xian-Liu-2010}. The settings of these parameters are shown in TABLE \ref{table-2}.
The uncertainties in the cost function of a wind generation station arise from the scale factor $\varsigma^{t}$ and the shape factor $\phi^{t}$ of the Weibull distribution of wind speed.
Here, $\varsigma^{t}$ is drawn from a uniform distribution in the range $[3,25]$, and $\phi^{t}$ is drawn from a uniform distribution in the range $[2,3]$.
The time-varying power demand $D^{t}$ is drawn from a Gaussian distribution with mean $70$ and variance $5^{2}$.

We randomly select $|\mathcal{B}|=1$ Byzantine wind generation station and investigate four Byzantine attacks: large-value, small-value, large-value Gaussian, and small-value Gaussian attacks. Specifically, in large-value Byzantine attacks, the Byzantine wind generation station sets its message as $-0.01$. In small-value Byzantine attacks, the Byzantine wind generation station sets its message as $-300$. In large-value Gaussian Byzantine attacks, the Byzantine wind generation station sets its message following a Gaussian distribution with mean $-10$ and variance $5$. In small-value Gaussian Byzantine attacks, the Byzantine wind generation station sets its message following a Gaussian distribution with mean $-150$ and variance $5$.
We consider three robust aggregation rules: $TM(\cdot)$, $IOS(ARC(\cdot))$ and $SCC(ARC(\cdot))$. In $TM(\cdot)$ and $IOS(ARC(\cdot))$, benign generation stations set the parameters $b$ as the number of Byzantine neighbors. In $SCC(ARC(\cdot))$, the clipping threshold $\tau$ is set according to Theorem 3 in \cite{b-LieHe-2022}.
The primal step size is $\alpha=1$, and the dual step size is $\beta=3$. The regularization parameter is $\theta=0.001$.

The numerical results are shown in Fig.  \ref{fig_2}. We use the performance of the attack-free decentralized online resource allocation algorithm under various Byzantine attacks as the baseline. Given the significant differences in terms of dynamic regret and accumulative constraint violation between the attack-free and Byzantine-resilient decentralized online resource allocation algorithms, we depict the numerical results of Byzantine-resilient and attack-free algorithms on the left and right ordinates, respectively.

Under small-value and small-value Gaussian attacks, the attack-free algorithm exhibits both linear dynamic regret and accumulative constraint violation. Similarly, the proposed Byzantine-resilient algorithms equipped with robust aggregation rules $TM(\cdot)$, $IOS(ARC(\cdot))$ and $SCC(ARC(\cdot))$ demonstrate linear dynamic regret and accumulative constraint violation, but the values are significantly smaller than that of the attack-free algorithm. This observation highlights the advantages of the proposed Byzantine-resilient algorithms.

Under large-value and large-value Gaussian attacks, the dynamic regret of the attack-free algorithm decreases linearly with respect to $T$. The reason is that the Byzantine generation station consistently sends wrong but large dual variables to the benign generation stations, resulting in smaller power generation strategies for the benign generation stations. Consequently, the instantaneous costs are always smaller than the optimal cost. Nevertheless, the accumulative constraint violation remain linear. The conclusions drawn from the small-value and small-value Gaussian attacks still hold true under the large-value and large-value Gaussian attacks.

\begin{table*}[!ht]
\centering
\renewcommand{\arraystretch}{1.18}
\caption{The parameters of renewable wind generation stations for case 2 \cite{b-Fanghong-Guo-2021, b-Fang-Yao-2012} }
\begin{tabular}{|c|c|c|c|c|c|c|c|c|c|} \hline
Wind generation station No. & $\varrho_{i}$ & $v_{in,i}$ & $v_{out,i}$ & $v_{r,i}$ & $\sigma_{ue,i}$ & $\sigma_{oe,i}$ & $P_{r,i}$ & $P_{wi,i}^{min}$ & $P_{wi,i}^{max}$\\ \hline
 55  &	 1  &	3 & 25 & 13 & 3 & 20& 150& 0 & 500\\ \hline
56  & 6 & 4 & 45 & 15 & 5 & 30 & 160 & 0 &  300 \\ \hline
57  & 1 & 5 & 25 & 16 & 3 & 20 & 150 & 0 & 400 \\ \hline
58  & 6 & 3 & 45 & 13 & 5 & 30 & 160 & 0 & 200 \\ \hline
59  & 1 & 4 & 25 & 15 & 3 & 20 & 150 & 0 & 300 \\ \hline
60  & 6 & 5 & 45 & 16 & 5 & 30 & 160 & 0 & 200 \\ \hline
\end{tabular}\\
\label{table-3}
\end{table*}

\begin{figure}[htbp]
\centerline{\includegraphics[width=8.5cm]{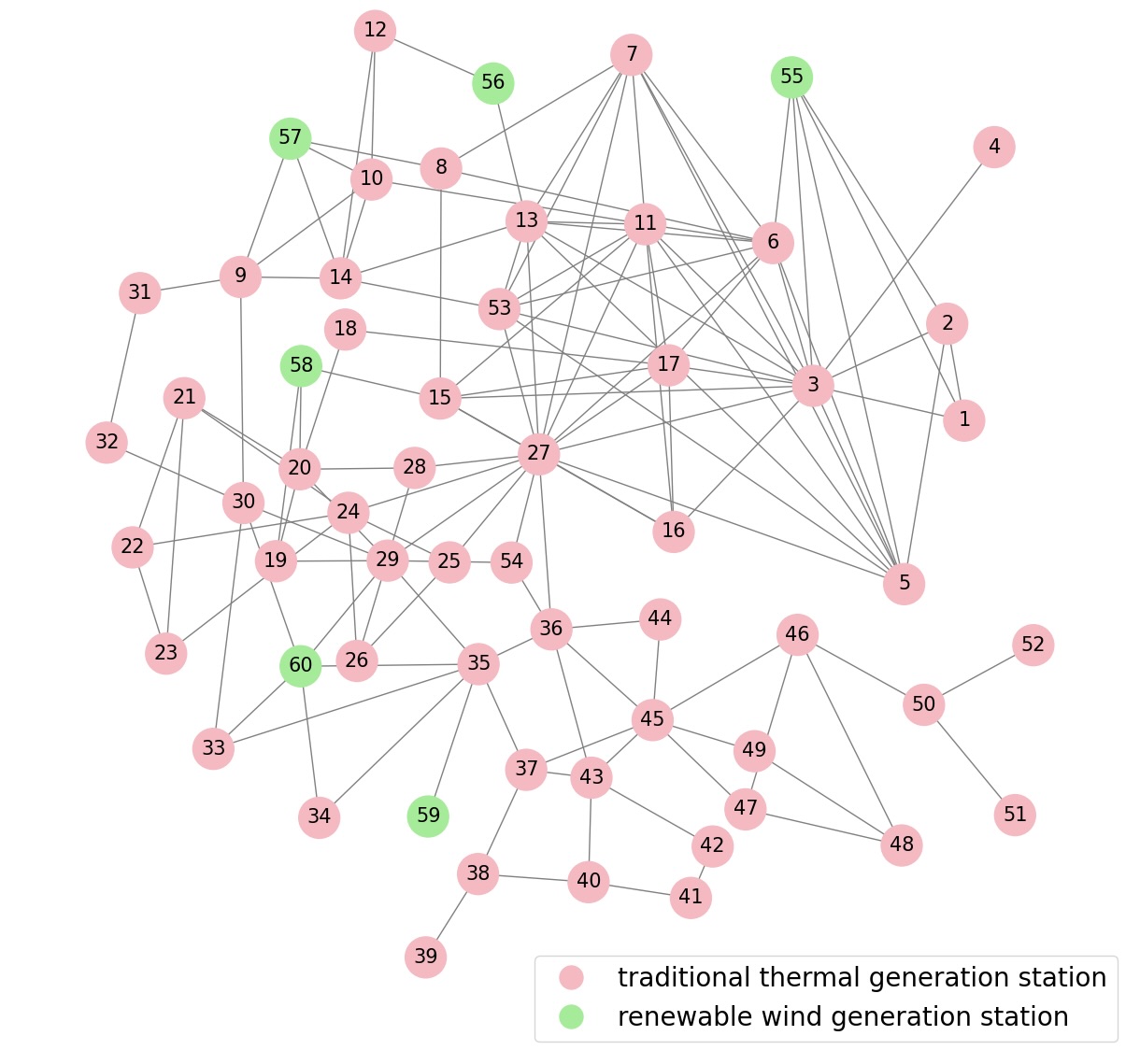}}
\caption{The communication graph of IEEE 118-Bus test system with 6 wind generation stations.}
\label{fig_3}
\end{figure}

\begin{figure*}[htbp]
\centerline{\includegraphics[width=17cm]{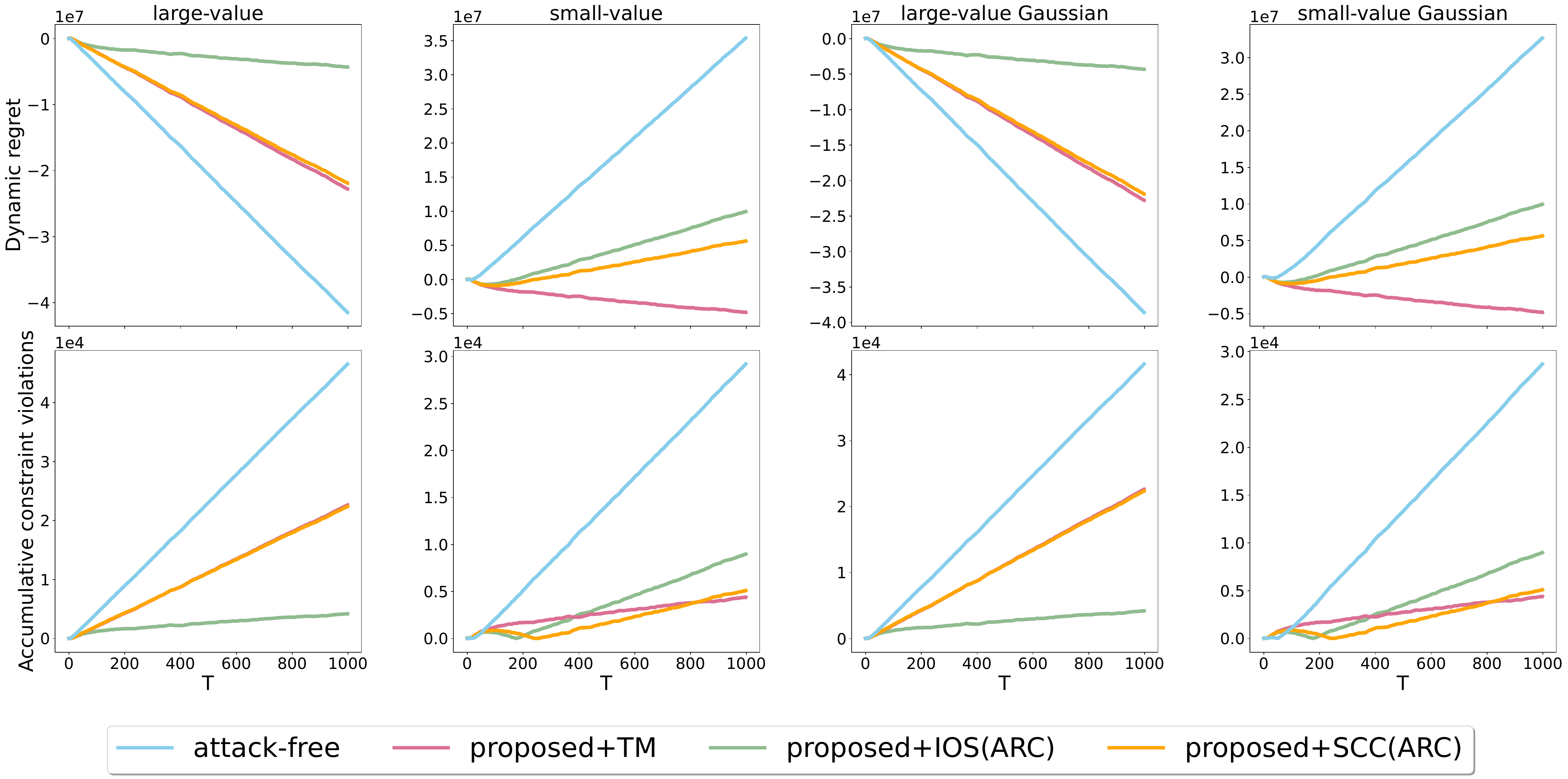}}
\caption{Dynamic regret and accumulative constraint violations of the compared algorithms under various Byzantine attacks.}
\label{fig_4}
\end{figure*}

\subsection{Case 2: IEEE 118-Bus Test System with 6 Wind Generation Stations}
Next, we consider the IEEE 118-bus test system which contains 54 traditional thermal generation stations \cite{b-IEEE-118}. We randomly select $6$ buses to deploy $6$ renewable wind generation stations. The resultant communication graph is shown in Fig. \ref{fig_3}. According to the communication graph, we use the Metropolis constant weight rule \cite{b-Wei-Shi-2015} to generate a doubly stochastic weight matrix $\widetilde{E}$.
Each traditional thermal generation station $i\in \{1,2,\cdots,54\}$ has a cost function $C_{i}(P_{i})=\eta_{i}(P_{i})^{2}+\zeta_{i}P_{i}+\xi_{i}$, where $\eta_{i}\in [0.0024, 0.0697]$, $\zeta_{i}\in [8.3391,37.6968]$, and $\zeta_{i}\in [6.78, 74.33]$. The local power constraint of each traditional thermal generation station $i$ is $P_{i}\in [P_{th,i}^{\min}, P_{th,i}^{\max}]$, where $P_{th,i}^{\min}\in [5,150]$ and $P_{th,i}^{\max}\in [30,420]$. The time-varying cost function of each renewable wind generation station $i\in \{55,56,\cdots,60\}$ is $C_{i}^{t}(P_{i})=\varrho_{i} P_{i}+C_{ue,i}^{t}(\varsigma^{t},\kappa^{t},\sigma_{ue,i},v_{in,i},v_{out,i},v_{r,i},P_{r,i}, P_{i})+ C_{oe,i}^{t}(\varsigma^{t}, \kappa^{t},\sigma_{oe,i},v_{in,i},v_{out,i},v_{r,i},P_{r,i}, P_{i})$. The settings of the cost parameters $\varrho_{i},v_{in,i},v_{out,i},v_{r,i},\sigma_{ue,i},\sigma_{oe,i},P_{r,i}$ and the local constraint parameters $P_{wi,i}^{\min}, P_{wi,i}^{\max}$ are shown in TABLE \ref{table-3}. The time-varying cost parameters $\varsigma^{t}$ and $\kappa^{t}$ are from the actual hourly wind speed data of the continental United States \cite{b-wind-speed-data}. The time-varying power demand $D^{t}$ is drawn from Gaussian distribution with mean $100$ and variance $10^{2}$.

We randomly select $|\mathcal{B}|=2$ Byzantine generation stations out of $60$ generation stations, and test the performance of dynamic regret and accumulative constraint violations of proposed algorithms under four types of Byzantine attacks, including large-value, small-value, large-value Gaussian, and small-value Gaussian. For large-value Byzantine attacks, Byzantine generation stations set their messages as $-0.01$. For small Byzantine attacks, Byzantine generation stations set their messages as $-2000$. For large-value Gaussian Byzantine attacks, Byzantine generation stations set their messages following a Gaussian distribution with mean $-500$ and variance $30^{2}$. For small-value Gaussian Byzantine attacks, Byzantine generation stations set their messages following a Gaussian distribution with mean $-1500$ and variance $30^{2}$. The primal and dual step sizes are $\alpha=\beta=5$. The regularization parameter is $\theta=0.00001$.

Fig. \ref{fig_4}. shows the performance of dynamic regret and accumulative constraint violations of the attack-free algorithm and our proposed algorithms with different robust aggregations, i.e., $TM(\cdot)$, $IOS(ARC(\cdot))$, and $SCC(ARC(\cdot))$. Under small-value and small-value Gaussian Byzantine attacks, the attack-free decentralized online resource allocation algorithm has a linear and large dynamic regret and accumulative constraint violations. However, the proposed Byzantine-resilient decentralized online resource allocation algorithms with robust aggregation rules $TM(\cdot)$, $IOS(ARC(\cdot))$ and $SCC(ARC(\cdot))$ have much smaller linear dynamic regret and accumulative constraint violations. Hence, the proposed algorithms are resilient.

Considering the characteristics of large-value and large-value Gaussian Byzantine attacks, we only focus on the performance comparison of the attack-free and the proposed Byzantine-resilient decentralized online resource allocation algorithms in terms of accumulative constraint violation. It is observed that the accumulative constraint violations of the attack-free decentralized online resource allocation algorithm are much larger than those of the proposed Byzantine-resilient decentralized online resource allocation algorithms.
\section{CONCLUSIONS AND FUTURE WORK}
\label{sec 7}
In this paper, we investigate decentralized online resource allocation under Byzantine attacks. We propose a class of Byzantine-resilient decentralized online resource allocation algorithms equipped with robust aggregation rules. Theoretically, when the robust aggregation rules are properly designed, the proposed algorithms will achieve linear dynamic regret and accumulative constraint violations. Experimental results corroborate our theoretically findings.

\begin{appendices}
\section{Proof of Theorem \ref{t2}}\label{appendix1}
\subsection{Proof of Theorem \ref{t2}}\label{aa1}
\begin{IEEEproof}
\label{proof-theorem2}
For notational convenience, we define a function given by $L_{i}^{t}(P):=\left \langle P-P_{i}^{t},\nabla C_{i}^{t}(P_{i}^{t}) +\frac{\lambda_{i}^{t}}{M}  \right \rangle +\frac{1}{2 \alpha}\|P-P_{i}^{t}\|^{2}$. Therefore, the update of primal variables $P_{i}^{t+1}$ in Algorithm \ref{alg2} can be rewritten as $P_{i}^{t+1}=\arg\min_{P\in \Omega_{i}} L_{i}^{t}(P)$. Given the definition of $L_{i}^{t}(P)$, we have $\nabla^{2}L_{i}^{t}(P)=\frac{1}{\alpha}>0$. Therefore, the function $L_{i}^{t}(P)$ is $\frac{1}{\alpha}$-strongly convex. According to the definition of a strongly convex function, we have
\begin{align}
\label{proof-t2-1}
L_{i}^{t}(P_{i}^{t*})\ge & L_{i}^{t}(P_{i}^{t+1})+\left \langle \nabla L_{i}^{t}(P_{i}^{t+1})  , P_{i}^{t*}-P_{i}^{t+1}  \right \rangle \\
&+ \frac{1}{2\alpha}\|P_{i}^{t*}-P_{i}^{t+1}\|^{2}. \notag
\end{align}
Since $P_{i}^{t+1}=\arg\min_{P\in \Omega_{i}} L_{i}^{t}(P)$, we obtain the optimality condition $\left \langle \nabla L_{i}^{t}(P_{i}^{t+1}), P_{i}^{t*}-P_{i}^{t+1}  \right \rangle \ge 0$. Hence, we have
\begin{align}
\label{proof-t2-2}
L_{i}^{t}(P_{i}^{t*})\ge L_{i}^{t}(P_{i}^{t+1})+ \frac{1}{2\alpha}\|P_{i}^{t*}-P_{i}^{t+1}\|^{2}.
\end{align}
By the definition $L_{i}^{t}(P):=\left \langle P-P_{i}^{t},\nabla C_{i}^{t}(P_{i}^{t}) + \frac{\lambda_{i}^{t} }{M} \right \rangle +\frac{1}{2 \alpha}\|P-P_{i}^{t}\|^{2}$, we can rewritten \eqref{proof-t2-2} as
\begin{align}
\label{proof-t2-3}
&\left \langle P_{i}^{t*}-P_{i}^{t},\nabla C_{i}^{t}(P_{i}^{t}) + \frac{\lambda_{i}^{t}}{M} \right \rangle +\frac{1}{2 \alpha}\|P_{i}^{t*}-P_{i}^{t}\|^{2}\ge  \\
&\left \langle P_{i}^{t+1}-P_{i}^{t},\nabla C_{i}^{t}(P_{i}^{t}) +\frac{\lambda_{i}^{t}}{M} \right \rangle +\frac{1}{2 \alpha}\|P_{i}^{t+1}-P_{i}^{t}\|^{2}\notag\\
&+ \frac{1}{2\alpha}\|P_{i}^{t*}-P_{i}^{t+1}\|^{2}. \notag
\end{align}
Adding $C_{i}^{t}(P_{i}^{t})$ to both sides of \eqref{proof-t2-3} and rearranging the terms, we obtain
\begin{align}
\label{proof-t2-4}
&C_{i}^{t}(P_{i}^{t})+\left \langle P_{i}^{t+1}-P_{i}^{t},\nabla C_{i}^{t}(P_{i}^{t}) + \frac{\lambda_{i}^{t}}{M} \right \rangle \\
&+\frac{1}{2 \alpha}\|P_{i}^{t+1}-P_{i}^{t}\|^{2}\notag\\
\le &C_{i}^{t}(P_{i}^{t})+\left \langle P_{i}^{t*}-P_{i}^{t},\nabla C_{i}^{t}(P_{i}^{t}) + \frac{\lambda_{i}^{t}}{M} \right \rangle \notag\\
&+\frac{1}{2 \alpha}(\|P_{i}^{t*}-P_{i}^{t}\|^{2}- \|P_{i}^{t*}-P_{i}^{t+1}\|^{2})\notag \\
\le & C_{i}^{t}(P_{i}^{t*})+\left \langle P_{i}^{t*}-P_{i}^{t},\frac{\lambda_{i}^{t}}{M} \right \rangle \notag\\
&+\frac{1}{2 \alpha}(\|P_{i}^{t*}-P_{i}^{t}\|^{2}- \|P_{i}^{t*}-P_{i}^{t+1}\|^{2}), \notag
\end{align}
where the last inequality holds because $C_{i}^{t}(\cdot)$ is convex, i.e., $C_{i}^{t}(P_{i}^{t})+\left \langle P_{i}^{t*}-P_{i}^{t},\nabla C_{i}^{t}(P_{i}^{t}) \right \rangle\le C_{i}^{t}(P_{i}^{t*})$. Rearranging \eqref{proof-t2-4}, we have
\begin{align}
\label{proof-t2-5}
&C_{i}^{t}(P_{i}^{t})-C_{i}^{t}(P_{i}^{t*})\\
\le& \underbrace{\left \langle P_{i}^{t*}-P_{i}^{t+1}, \frac{\lambda_{i}^{t}}{M} \right \rangle}_{A_{1}} \underbrace{-\left \langle P_{i}^{t+1}-P_{i}^{t},\nabla C_{i}^{t}(P_{i}^{t}) \right \rangle}_{A_{2}} \notag\\
&+\underbrace{\frac{1}{2 \alpha}(\|P_{i}^{t*}-P_{i}^{t}\|^{2}- \|P_{i}^{t*}-P_{i}^{t+1}\|^{2})}_{A_{3}}-\frac{1}{2 \alpha}\|P_{i}^{t+1}-P_{i}^{t}\|^{2}.\notag
\end{align}

Next, we analyze $A_{1}$, $A_{2}$ and $A_{3}$ in turn.\\
\noindent\textbf{Bounding $A_{1}$:}
According to the definition of $G_{i}^{t}(P_{i})=\frac{1}{H} P_{i}-\frac{1}{H}D^{t}$, we obtain
\begin{align}
\label{proof-t2-6}
A_{1}=&\left \langle P_{i}^{t*}-P_{i}^{t+1},\frac{\lambda_{i}^{t}}{M}\right \rangle \\
=&\frac{H}{M}\left \langle G_{i}^{t}(P_{i}^{t*})-G_{i}^{t}(P_{i}^{t+1}), \lambda_{i}^{t} \right \rangle\notag\\
=& \frac{H}{M} \left \langle \lambda_{i}^{t} ,G_{i}^{t}(P_{i}^{t*}) \right \rangle-\frac{H}{M}\left \langle \lambda_{i}^{t} ,G_{i}^{t}(P_{i}^{t+1}) \right \rangle\notag\\
&+\frac{H}{M}\left \langle \bar{\lambda}^{t} ,G_{i}^{t}(P_{i}^{t+1}) \right \rangle-\frac{H}{M}\left \langle \bar{\lambda}^{t} ,G_{i}^{t}(P_{i}^{t+1}) \right \rangle\notag\\
=& \frac{H}{M}\left \langle \lambda_{i}^{t} ,G_{i}^{t}(P_{i}^{t*}) \right \rangle+\frac{H}{M}\left \langle \bar{\lambda}^{t}-\lambda_{i}^{t} ,G_{i}^{t}(P_{i}^{t+1}) \right \rangle\notag\\
&-\frac{H}{M}\left \langle \bar{\lambda}^{t} ,G_{i}^{t}(P_{i}^{t+1}) \right \rangle+ \frac{H}{M}\left \langle \bar{\lambda}^{t} ,G_{i}^{t}(P_{i}^{t*}) \right \rangle\notag\\
&-\frac{H}{M}\left \langle \bar{\lambda}^{t} ,G_{i}^{t}(P_{i}^{t*}) \right \rangle\notag \\
=& \frac{H}{M}\left \langle \lambda_{i}^{t}-\bar{\lambda}^{t} ,G_{i}^{t}(P_{i}^{t*}) \right \rangle+\frac{H}{M}\left \langle \bar{\lambda}^{t}-\lambda_{i}^{t} ,G_{i}^{t}(P_{i}^{t+1}) \right \rangle\notag\\
&+ \frac{H}{M}\left \langle \bar{\lambda}^{t} ,G_{i}^{t}(P_{i}^{t*}) \right \rangle-\frac{H}{M}\left \langle \bar{\lambda}^{t} ,G_{i}^{t}(P_{i}^{t+1}) \right \rangle. \notag
\end{align}
\noindent\textbf{Bounding $A_{2}$:}
Under Assumption \ref{a1}, we obtain
\begin{align}
\label{proof-t2-7}
A_{2}=&-\left \langle P_{i}^{t+1}-P_{i}^{t},\nabla C_{i}^{t}(P_{i}^{t}) \right \rangle \\
\le & \|P_{i}^{t+1}-P_{i}^{t}\|\nabla C_{i}^{t}(P_{i}^{t})\|\notag\\
\le & \frac{u_{1}}{2}\cdot \|P_{i}^{t+1}-P_{i}^{t}\|^{2}+\frac{1}{2u_{1}}\cdot \|\nabla C_{i}^{t}(P_{i}^{t})\|^{2}\notag\\
\le & \frac{u_{1}}{2}\cdot \|P_{i}^{t+1}-P_{i}^{t}\|^{2}+\frac{\varphi^{2}}{2u_{1}}, \notag
\end{align}
where $u_{1}>0$ is any positive constant.

\noindent\textbf{Bounding $A_{3}$:} Under Assumption \ref{a1}, we obtain
\begin{align}
\label{proof-t2-8}
A_{3}=& \frac{1}{2 \alpha}(\|P_{i}^{t*}-P_{i}^{t}\|^{2}- \|P_{i}^{t*}-P_{i}^{t+1}\|^{2})\\
=& \frac{1}{2 \alpha}(\|P_{i}^{t*}-P_{i}^{t}\|^{2}- \|P_{i}^{t}-P_{i}^{t-1*}\|^{2}+\|P_{i}^{t}-P_{i}^{t-1*}\|^{2} \notag\\
&-\|P_{i}^{t*}-P_{i}^{t+1}\|^{2})\notag\\
=& \frac{1}{2 \alpha}(\|P_{i}^{t*}-P_{i}^{t-1*}\|\cdot\|P_{i}^{t*}-2P_{i}^{t}+P_{i}^{t-1*}\|\notag\\
&+\|P_{i}^{t}-P_{i}^{t-1*}\|^{2} -\|P_{i}^{t*}-P_{i}^{t+1}\|^{2})\notag\\
=& \frac{1}{2 \alpha}[\|P_{i}^{t*}-P_{i}^{t-1*}\|\cdot(\|P_{i}^{t*}-P_{i}^{t}\|+\|P_{i}^{t}-P_{i}^{t-1*}\|)\notag\\
&+\|P_{i}^{t}-P_{i}^{t-1*}\|^{2} -\|P_{i}^{t*}-P_{i}^{t+1}\|^{2}]\notag\\
\le & \frac{R}{\alpha} \|P_{i}^{t*}-P_{i}^{t-1*}\|+\frac{1}{2\alpha}(\|P_{i}^{t}-P_{i}^{t-1*}\|^{2} -\|P_{i}^{t*}-P_{i}^{t+1}\|^{2})\notag.
\end{align}

Substituting \eqref{proof-t2-6}, \eqref{proof-t2-7} and \eqref{proof-t2-8} into \eqref{proof-t2-5} and rearranging
the terms, we have
\begin{align}
\label{proof-t2-9}
&C_{i}^{t}(P_{i}^{t})-C_{i}^{t}(P_{i}^{t*})\\
\le &(\frac{u_{1}}{2}-\frac{1}{2\alpha})\|P_{i}^{t+1}-P_{i}^{t}\|^{2}+\frac{R}{\alpha} \|P_{i}^{t*}-P_{i}^{t-1*}\|\notag\\
&+\frac{1}{2\alpha}(\|P_{i}^{t}-P_{i}^{t-1*}\|^{2} -\|P_{i}^{t*}-P_{i}^{t+1}\|^{2})\notag\\
&+\frac{H}{M}\left \langle \bar{\lambda}^{t} ,G_{i}^{t}(P_{i}^{t*}) \right \rangle-\frac{H}{M}\left \langle \bar{\lambda}^{t} ,G_{i}^{t}(P_{i}^{t+1}) \right \rangle\notag\\
&+\frac{H}{M}\left \langle \lambda_{i}^{t}-\bar{\lambda}^{t} ,G_{i}^{t}(P_{i}^{t*}) \right \rangle+\frac{H}{M}\left \langle \bar{\lambda}^{t}-\lambda_{i}^{t} ,G_{i}^{t}(P_{i}^{t+1}) \right \rangle+\frac{\varphi^{2}}{2u_{1}}.\notag
\end{align}
Since $\bm{P}^{t*}:=[P_{1}^{t*},\cdots,P_{H}^{t*}]$  is the optimal solution of problem \eqref{online-decentralized-economic-dispatch-problem-oracle} at each time period $t$, we have $\sum_{i\in \mathcal{H}}G_{i}^{t}(P_{i}^{t*})=0$. Summing over $i \in  \mathcal{H}$ on both sides of \eqref{proof-t2-9}, we have
\begin{align}
\label{proof-t2-10}
&\sum_{i\in \mathcal{H}}C_{i}^{t}(P_{i}^{t})-\sum_{i\in \mathcal{H}}C_{i}^{t}(P_{i}^{t*})\\
\le &(\frac{u_{1}}{2}-\frac{1}{2\alpha})\sum_{i\in \mathcal{H}}\|P_{i}^{t+1}-P_{i}^{t}\|^{2}+\frac{R}{\alpha} \sum_{i\in \mathcal{H}}\|P_{i}^{t*}-P_{i}^{t-1*}\|\notag\\
&+\frac{1}{2\alpha}\sum_{i\in \mathcal{H}}(\|P_{i}^{t}-P_{i}^{t-1*}\|^{2} -\|P_{i}^{t*}-P_{i}^{t+1}\|^{2})\notag\\
&-\frac{H}{M}\sum_{i\in \mathcal{H}}\left \langle \bar{\lambda}^{t} ,G_{i}^{t}(P_{i}^{t+1}) \right \rangle+\underbrace{\frac{H}{M}\sum_{i\in \mathcal{H}}\left \langle \lambda_{i}^{t}-\bar{\lambda}^{t} ,G_{i}^{t}(P_{i}^{t*}) \right \rangle}_{A_{4}}\notag\\
&+\underbrace{\frac{H}{M}\sum_{i\in \mathcal{H}}\left \langle \bar{\lambda}^{t}-\lambda_{i}^{t} ,G_{i}^{t}(P_{i}^{t+1}) \right \rangle}_{A_{5}}+\frac{\varphi^{2}H}{2u_{1}}.\notag
\end{align}

Next, we analyze $A_{4}$ and $A_{5}$ in turn.

\noindent\textbf{Bounding $A_{4}$:}
Based on Assumption \ref{a1}, Lemma \ref{lemma4} and the fact $\sum_{i\in \mathcal{H}} \|\lambda_{i}^{t}-\bar{\lambda}^{t}\|\le \sqrt{H}\cdot\|\Lambda^{t}-\frac{1}{H}\bm{1}\bm{1}^{\top}\Lambda^{t}
\|_{F}$, we obtain
\begin{align}
\label{proof-t2-11}
A_{4}=&\frac{H}{M}\sum_{i\in \mathcal{H}}\left \langle \lambda_{i}^{t}-\bar{\lambda}^{t} ,G_{i}^{t}(P_{i}^{t*}) \right \rangle \\
\le & \frac{H}{M}\sum_{i\in \mathcal{H}} \|\lambda_{i}^{t}-\bar{\lambda}^{t}\|\|G_{i}^{t}(P_{i}^{t*})\|\notag \\
\le& \frac{2H^{3}\psi^{2}\beta}{\epsilon\sqrt{\epsilon}M^{2}}. \notag
\end{align}

\noindent\textbf{Bounding $A_{5}$:}
Similar to the derivation of \eqref{proof-t2-11}, we obtain
\begin{align}
\label{proof-t2-12}
A_{5}=&\frac{H}{M}\sum_{i\in \mathcal{H}}\left \langle \bar{\lambda}^{t}-\lambda_{i}^{t} ,G_{i}^{t}(P_{i}^{t+1}) \right \rangle \\
\le& \frac{2H^{3}\psi^{2}\beta}{\epsilon\sqrt{\epsilon}M^{2}}. \notag
\end{align}
Substituting \eqref{proof-t2-11} and \eqref{proof-t2-12} into \eqref{proof-t2-10}, we have
\begin{align}
\label{proof-t2-13}
&\sum_{i\in \mathcal{H}}C_{i}^{t}(P_{i}^{t})-\sum_{i\in \mathcal{H}}C_{i}^{t}(P_{i}^{t*})\\
\le &(\frac{u_{1}}{2}-\frac{1}{2\alpha})\sum_{i\in \mathcal{H}}\|P_{i}^{t+1}-P_{i}^{t}\|^{2}+\frac{R}{\alpha} \sum_{i\in \mathcal{H}}\|P_{i}^{t*}-P_{i}^{t-1*}\|\notag\\
&+\frac{1}{2\alpha}\sum_{i\in \mathcal{H}}(\|P_{i}^{t}-P_{i}^{t-1*}\|^{2} -\|P_{i}^{t*}-P_{i}^{t+1}\|^{2})\notag\\
&-\frac{H}{M}\sum_{i\in \mathcal{H}}\left \langle \bar{\lambda}^{t} ,G_{i}^{t}(P_{i}^{t+1}) \right \rangle+\frac{4H^{3}\psi^{2}\beta}{\epsilon\sqrt{\epsilon}M^{2}}+\frac{\varphi^{2}H}{2u_{1}}.\notag
\end{align}
Combining \eqref{proof-t2-13} and Lemma \ref{lemma5}, we have
\begin{align}
\label{proof-t2-14}
&\frac{1-\beta\theta}{2\beta}\cdot\Delta^{t}+\sum_{i\in \mathcal{H}}C_{i}^{t}(P_{i}^{t})-\sum_{i\in \mathcal{H}}C_{i}^{t}(P_{i}^{t*})\\
\le &(\frac{u_{1}}{2}-\frac{1}{2\alpha})\sum_{i\in \mathcal{H}}\|P_{i}^{t+1}-P_{i}^{t}\|^{2}+\frac{R}{\alpha} \sum_{i\in \mathcal{H}}\|P_{i}^{t*}-P_{i}^{t-1*}\|\notag\\
&+\frac{1}{2\alpha}\sum_{i\in \mathcal{H}}(\|P_{i}^{t}-P_{i}^{t-1*}\|^{2} -\|P_{i}^{t*}-P_{i}^{t+1}\|^{2})\notag \\
&+\underbrace{\frac{H}{M} \sum_{i\in \mathcal{H}} \left \langle  \bar{\lambda}^{t}  , G_{i}^{t}(P_{i}^{t})  \right \rangle-\frac{H}{M}\sum_{i\in \mathcal{H}}\left \langle \bar{\lambda}^{t} ,G_{i}^{t}(P_{i}^{t+1}) \right \rangle}_{A_{6}}\notag\\
&-\frac{ H}{M} \sum_{i\in \mathcal{H}} \left \langle  \lambda  , G_{i}^{t}(P_{i}^{t})  \right \rangle+(\frac{4H^{3}\psi^{2}}{\epsilon\sqrt{\epsilon}M^{2}}+\frac{H^{3}\psi^{2}}{M^{2}}+\psi^{2}H)\cdot\beta\notag\\
&+(1+\frac{1}{\epsilon^{3}})\cdot (4\rho H+\chi) \cdot \frac{8 H^{3}\psi^{2}}{M^{2}\theta}+\frac{\varphi^{2}H}{2u_{1}}+\frac{\theta H}{2}\|\lambda\|^{2}.\notag
\end{align}

Next we analyze the term $A_{6}$.\\
\noindent\textbf{Bounding $A_{6}$:}
According to the definition $G_{i}^{t}(P_{i})=\frac{1}{H} P_{i}-\frac{1}{H}D^{t}$, Assumption \ref{a1}, Lemma \ref{lemma2} and Lemma \ref{lemma4}, we have
\begin{align}
\label{proof-t2-15}
A_{6}=&\frac{H}{M} \sum_{i\in \mathcal{H}} \left \langle  \bar{\lambda}^{t}  , G_{i}^{t}(P_{i}^{t})  \right \rangle-\frac{H}{M}\sum_{i\in \mathcal{H}}\left \langle \bar{\lambda}^{t} ,G_{i}^{t}(P_{i}^{t+1}) \right \rangle\\
=& \frac{1}{M}\sum_{i\in \mathcal{H}}\left \langle \bar{\lambda}^{t}  ,P_{i}^{t}- P_{i}^{t+1} \right \rangle\notag\\
\le & \frac{u_{2}}{2M}\sum_{i\in \mathcal{H}} \|P_{i}^{t+1}- P_{i}^{t}\|^{2}+\frac{H}{2u_{2}M} \|\bar{\lambda}^{t}\|^{2}\notag\\
\le & \frac{u_{2}}{2M}\sum_{i\in \mathcal{H}} \|P_{i}^{t+1}- P_{i}^{t}\|^{2}+\frac{H}{2u_{2}M}\cdot\frac{\psi^{2} }{\theta^{2}},\notag
\end{align}
where $u_{2}>0$ is any positive constant. Letting $u_{2}=\frac{M}{2\alpha}$, we can rewrite \eqref{proof-t2-15} as
\begin{align}
\label{proof-t2-16}
A_{6}\le \frac{1}{4\alpha}\sum_{i\in \mathcal{H}} \|P_{i}^{t+1}- P_{i}^{t}\|^{2}+\frac{H}{M^{2}}\cdot\frac{\psi^{2} \alpha}{\theta^{2}}.
\end{align}
Substituting \eqref{proof-t2-16} into \eqref{proof-t2-14} and rearranging the terms, we have
\begin{align}
\label{proof-t2-17}
&\frac{1-\beta\theta}{2\beta}\cdot\Delta^{t}+\sum_{i\in \mathcal{H}}C_{i}^{t}(P_{i}^{t})-\sum_{i\in \mathcal{H}}C_{i}^{t}(P_{i}^{t*})\\
\le &(\frac{u_{1}}{2}+\frac{1}{4\alpha}-\frac{1}{2\alpha})\sum_{i\in \mathcal{H}}\|P_{i}^{t+1}-P_{i}^{t}\|^{2}+\frac{R}{\alpha} \sum_{i\in \mathcal{H}}\|P_{i}^{t*}-P_{i}^{t-1*}\|\notag\\
&+\frac{1}{2\alpha}\sum_{i\in \mathcal{H}}(\|P_{i}^{t}-P_{i}^{t-1*}\|^{2} -\|P_{i}^{t*}-P_{i}^{t+1}\|^{2})\notag\\
&+(\frac{4H^{3}\psi^{2}}{\epsilon\sqrt{\epsilon}M^{2}}+\frac{H^{3}\psi^{2}}{M^{2}}+\psi^{2}H)\cdot\beta+\frac{\varphi^{2}H}{2u_{1}}+\frac{H}{M^{2}}\cdot\frac{\psi^{2} \alpha}{\theta^{2}}\notag\\
&+(1+\frac{1}{\epsilon^{3}})\cdot (4\rho H+\chi) \cdot \frac{8 H^{3}\psi^{2}}{M^{2}\theta}\notag\\
&-\frac{ H}{M} \sum_{i\in \mathcal{H}} \left \langle  \lambda  , G_{i}^{t}(P_{i}^{t})  \right \rangle+\frac{\theta H}{2}\|\lambda\|^{2}\notag\\
= & \frac{R}{\alpha} \sum_{i\in \mathcal{H}}\|P_{i}^{t*}-P_{i}^{t-1*}\|\notag\\
&+\frac{1}{2\alpha}\sum_{i\in \mathcal{H}}(\|P_{i}^{t}-P_{i}^{t-1*}\|^{2} -\|P_{i}^{t*}-P_{i}^{t+1}\|^{2})\notag\\
&+\varphi^{2}H\cdot \alpha
+(\frac{4H^{3}\psi^{2}}{\epsilon\sqrt{\epsilon}M^{2}}+\frac{H^{3}\psi^{2}}{M^{2}}+\psi^{2}H)\cdot\beta + \frac{H}{M^{2}}\cdot\frac{\psi^{2} \alpha}{\theta^{2}} \notag\\
&+(1+\frac{1}{\epsilon^{3}})\cdot (4\rho H+\chi) \cdot \frac{8 H^{3}\psi^{2}}{M^{2}\theta}\notag\\
&-\frac{ H}{M} \sum_{i\in \mathcal{H}} \left \langle  \lambda  , G_{i}^{t}(P_{i}^{t})  \right \rangle+\frac{\theta H}{2}\|\lambda\|^{2},\notag
\end{align}
where the last equality holds by setting $u_{1}=\frac{1}{2\alpha}$.  Summing over $t \in  [1,T]$ on both sides of \eqref{proof-t2-17}, we have
\begin{align}
\label{proof-t2-18}
\hspace{-1em} &Reg^{T}_{\mathcal{H}}\\
\hspace{-1em}\le &\frac{R}{\alpha} \sum_{t=1}^{T} \sum_{i\in \mathcal{H}}\|P_{i}^{t*}-P_{i}^{t-1*}\|\notag\\
\hspace{-1em}&+\underbrace{\frac{1}{2\alpha}\sum_{t=1}^{T}\sum_{i\in \mathcal{H}}(\|P_{i}^{t}-P_{i}^{t-1*}\|^{2} -\|P_{i}^{t*}-P_{i}^{t+1}\|^{2})}_{A_{7}}\notag\\
\hspace{-1em}&+\varphi^{2}H\cdot \alpha T
+(\frac{4H^{3}\psi^{2}}{\epsilon\sqrt{\epsilon}M^{2}}+\frac{H^{3}\psi^{2}}{M^{2}}+\psi^{2}H)\cdot\beta T + \frac{H\psi^{2}}{M^{2}}\cdot\frac{\alpha T}{\theta^{2}} \notag\\
\hspace{-1em}&+(1+\frac{1}{\epsilon^{3}})\cdot (4\rho H+\chi) \cdot \frac{8 H^{3}\psi^{2}}{M^{2}}\cdot \frac{T}{\theta}\notag\\
\hspace{-1em}&-\frac{ H}{M}\sum_{t=1}^{T} \sum_{i\in \mathcal{H}} \left \langle  \lambda  , G_{i}^{t}(P_{i}^{t})  \right \rangle+\frac{\theta H T}{2}\|\lambda\|^{2}\underbrace{-\frac{1-\beta\theta}{2\beta}\cdot \sum_{t=1}^{T}\Delta^{t}}_{A_{8}}.\notag
\end{align}

Next, we analyze the terms $A_{7}$ and $A_{8}$ in turn. \\
\noindent\textbf{Bounding $A_{7}$:} It holds that
\begin{align}
\label{proof-t2-19}
A_{7}=&\frac{1}{2\alpha}\sum_{t=1}^{T}\sum_{i\in \mathcal{H}}(\|P_{i}^{t}-P_{i}^{t-1*}\|^{2} -\|P_{i}^{t*}-P_{i}^{t+1}\|^{2}) \\
=& \frac{1}{2\alpha} \sum_{i\in \mathcal{H}}[\|P_{i}^{1}-P_{i}^{0*}\|^{2} -\|P_{i}^{1*}-P_{i}^{2}\|^{2}+\cdots \notag\\
&+\|P_{i}^{T}-P_{i}^{T-1*}\|^{2} -\|P_{i}^{T*}-P_{i}^{T+1}\|^{2}]\notag\\
=& \frac{1}{2\alpha} \sum_{i\in \mathcal{H}}[\|P_{i}^{1}-P_{i}^{0*}\|^{2}-\|P_{i}^{T*}-P_{i}^{T+1}\|^{2}]\notag\\
\le & \frac{1}{2\alpha} \sum_{i\in \mathcal{H}}\|P_{i}^{1}-P_{i}^{0*}\|^{2}. \notag
\end{align}
\noindent\textbf{Bounding $A_{8}$:}
According to the definition $\Delta^{t}:=H \|\bar{\lambda}^{t+1}-\lambda\|^{2}-H \|\bar{\lambda}^{t}-\lambda\|^{2}$, we have
\begin{align}
\label{proof-t2-20}
A_{8}=&-\frac{(1-\beta\theta)\cdot H}{2\beta}\cdot \sum_{t=1}^{T}\Delta^{t} \\
=&-\frac{(1-\beta\theta)\cdot H}{2\beta}[\|\bar{\lambda}^{2}-\lambda\|^{2}-\|\bar{\lambda}^{1}-\lambda\|^{2}\notag\\
&+\|\bar{\lambda}^{3}-\lambda\|^{2}-\|\bar{\lambda}^{2}-\lambda\|^{2}+\cdots+\|\bar{\lambda}^{T}-\lambda\|^{2}\notag\\
&-\|\bar{\lambda}^{T-1}-\lambda\|^{2}+\|\bar{\lambda}^{T+1}-\lambda\|^{2}-\|\bar{\lambda}^{T}-\lambda\|^{2}]\notag\\
=&-\frac{(1-\beta\theta)\cdot H}{2\beta}[\|\bar{\lambda}^{T+1}-\lambda\|^{2}-\|\bar{\lambda}^{1}-\lambda\|^{2}]\notag\\
\le &\frac{(1-\beta\theta)\cdot H}{2\beta}\|\bar{\lambda}^{1}-\lambda\|^{2}\notag\\
\le & \frac{H}{2\beta}\|\lambda\|^{2}, \notag
\end{align}
where the last inequality holds, since $\bar{\lambda}^{1}=0$ which is true based on initialization $P_{i}^{0}=\lambda_{i}^{0}=D^{0}=0$.

Substituting \eqref{proof-t2-19} and \eqref{proof-t2-20} into \eqref{proof-t2-18} and rearranging the terms, we have
\begin{align}
\label{proof-t2-21}
&\hspace{-1em} Reg^{T}_{\mathcal{H}}+\frac{ H}{M}\sum_{t=1}^{T} \sum_{i\in \mathcal{H}} \left \langle  \lambda  , G_{i}^{t}(P_{i}^{t})  \right \rangle- (\frac{\theta H T}{2}+\frac{H}{2\beta})\|\lambda\|^{2} \\
\le &\frac{R}{\alpha} \sum_{t=1}^{T} \sum_{i\in \mathcal{H}}\|P_{i}^{t*}-P_{i}^{t-1*}\|+\frac{1}{2\alpha} \sum_{i\in \mathcal{H}}\|P_{i}^{1}-P_{i}^{0*}\|^{2}\notag\\
+&\varphi^{2}H\cdot \alpha T
+(\frac{4H^{3}\psi^{2}}{\epsilon\sqrt{\epsilon}M^{2}}+\frac{H^{3}\psi^{2}}{M^{2}}+\psi^{2}H)\cdot\beta T + \frac{H\psi^{2}}{M^{2}}\cdot\frac{\alpha T}{\theta^{2}} \notag\\
+&(1+\frac{1}{\epsilon^{3}})\cdot (4\rho H+\chi) \cdot \frac{8 H^{3}\psi^{2}}{M^{2}}\cdot \frac{T}{\theta}. \notag
\end{align}
i) Substituting $\lambda=0$ into \eqref{proof-t2-21} and rearranging the terms, we have
\begin{align}
\label{proof-t2-22}
&Reg^{T}_{\mathcal{H}} \\
\le &\frac{R}{\alpha} \sum_{t=1}^{T} \sum_{i\in \mathcal{H}}\|P_{i}^{t*}-P_{i}^{t-1*}\|+\frac{1}{2\alpha} \sum_{i\in \mathcal{H}}\|P_{i}^{1}-P_{i}^{0*}\|^{2}\notag\\
&+\varphi^{2}H\cdot \alpha T
+(\frac{4H^{3}\psi^{2}}{\epsilon\sqrt{\epsilon}M^{2}}+\frac{H^{3}\psi^{2}}{M^{2}}+\psi^{2}H)\cdot\beta T + \frac{H\psi^{2}}{M^{2}}\cdot\frac{\alpha T}{\theta^{2}} \notag\\
&+(1+\frac{1}{\epsilon^{3}})\cdot (4\rho H+\chi) \cdot \frac{8 H^{3}\psi^{2}}{M^{2}}\cdot \frac{T}{\theta}. \notag
\end{align}
ii) Substituting $\lambda=\frac{\sum_{t=1}^{T}\sum_{i\in \mathcal{H}}G_{i}^{t}(P_{i}^{t})}{2(\frac{\theta H T}{2}+\frac{H}{2 \beta})}$ into \eqref{proof-t2-21} and rearranging the terms, we have
\begin{align}
\label{proof-t2-23}
&\|\sum_{t=1}^{T}\sum_{i\in \mathcal{H}}G_{i}^{t}(P_{i}^{t})\|^{2} \\
\le &  \frac{M}{2H-M}\cdot [\underline{2H \theta T}+\frac{2H}{\beta}] \cdot [\frac{R}{\alpha} \sum_{t=1}^{T} \sum_{i\in \mathcal{H}}\|P_{i}^{t*}-P_{i}^{t-1*}\|\notag\\
&+\frac{1}{2\alpha} \sum_{i\in \mathcal{M}}\|P_{i}^{1}-P_{i}^{0*}\|^{2} +\varphi^{2}H\cdot\alpha T+\frac{H\psi^{2}}{M^{2}}\cdot \frac{\alpha T}{\theta^{2}} \notag\\
&+(\frac{4H^{3}\psi^{2}}{\epsilon\sqrt{\epsilon}M^{2}}+\frac{H^{3}\psi^{2}}{M^{2}}+\psi^{2}H)\cdot \beta T+2HF\cdot T\notag\\
&+(1+\frac{1}{\epsilon^{3}})\cdot (4\rho H+\chi) \cdot \frac{8 H^{3}\psi^{2}}{M^{2}}\cdot \frac{T}{\theta}].\notag
\end{align}
\end{IEEEproof}

\noindent\textbf{Supporting Lemmas for Proof of Theorem \ref{t2}}
\begin{lemma}
\label{lemma2}
Suppose that the robust aggregation rule $AGG(\cdot)$ satisfies Property \ref{d2}. Under Assumptions \ref{a1} and \ref{a3}, for any agent $i\in \mathcal{H}$ and $t\in [0,\cdots,T]$, $\lambda_{i}^{t+\frac{1}{2}}$ and $\lambda_{i}^{t+1}$ generated by Algorithm \ref{alg2} satisfy
\begin{align}
\label{l2}
\|\lambda_{i}^{t+\frac{1}{2}}\|\le \frac{\psi}{\theta}, \quad \|\lambda_{i}^{t+1}\|\le \frac{\psi}{\theta}.
\end{align}
\begin{IEEEproof}
\label{proof-l2}
Combining the initializations $P_{i}^{0}=\lambda_{i}^{0}=D^{0}=0$ and the updates of $\lambda_{i}^{t+\frac{1}{2}}$ and $\lambda_{i}^{t+1}$ in Algorithm \ref{alg2}, we are able to derive $\|\lambda_{i}^{0+\frac{1}{2}}\|=0\le \frac{\psi}{\theta}$ and $\|\lambda_{i}^{0+1}\|=AGG(\lambda_{i}^{0+\frac{1}{2}}, \{\check{\lambda}_{j}^{0+\frac{1}{2}}\}_{j\in \mathcal{N}_{i}})\le \max_{j\in (\mathcal{N}_{i}\cap \mathcal{H})\cup \{i\}}\|\lambda_{j}^{0+\frac{1}{2}}\|=0\le \frac{\psi}{\theta}$. Therefore, when $t=0$, the propositions $\|\lambda_{i}^{t+\frac{1}{2}}\|\le \frac{\psi}{\theta}$ and $\|\lambda_{i}^{t+1}\|\le \frac{\psi}{\theta}$ hold.

Next, we prove the conclusion by mathematical induction. Suppose that when $t=t^{'}$, the propositions  $\|\lambda_{i}^{t^{'}+\frac{1}{2}}\|\le \frac{\psi}{\theta}$ and $\|\lambda_{i}^{t^{'}+1}\|\le \frac{\psi}{\theta}$ hold. We analyze when $t=t^{'}+1$, whether $\|\lambda_{i}^{t^{'}+1+\frac{1}{2}}\|\le \frac{\psi}{\theta}$ and $\|\lambda_{i}^{t^{'}+1+1}\|\le \frac{\psi}{\theta}$ hold. According to the update of $\lambda_{i}^{t+\frac{1}{2}}$ in Algorithm \ref{alg2} and the relationship $\widetilde{G}_{i}^{t}(P_{i})=\frac{H}{M}G_{i}^{t}(P_{i})$, we have
\begin{align}
\label{proof-l2-1}
\|\lambda_{i}^{t^{'}+1+\frac{1}{2}}\|&=\|\lambda_{i}^{t^{'}+1}+\beta\cdot(\widetilde{G}_{i}^{t^{'}+1}(P_{i}^{t^{'}+1})-\theta \lambda_{i}^{t^{'}+1})\|\\
&\le (1-\beta \theta)\|\lambda_{i}^{t^{'}+1}\|+\beta \frac{H}{M}\|G_{i}^{t^{'}+1}(P_{i}^{t^{'}+1})\|\notag\\
&\le (1-\beta \theta)\cdot \frac{\psi}{\theta}+ \frac{\beta H}{M} \psi \notag\\
&=\frac{\psi}{\theta}+(\frac{ H}{M}-1)\beta \psi \notag \\
& \le \frac{\psi}{\theta},\notag
\end{align}
where the second inequality holds based on $\|\lambda_{i}^{t^{'}+1}\|\le \frac{\psi}{\theta}$ and Assumption \ref{a1}. To derive the last inequality, we use the fact that $\frac{ H}{M}-1 \le 0$. According to the update of $\lambda_{i}^{t+1}$ in Algorithm \ref{alg2} and Property \ref{d2}, we have
\begin{align}
\label{proof-l2-2}
\|\lambda_{i}^{t^{'}+1+1}\|&=\|AGG(\lambda_{i}^{t^{'}+1+\frac{1}{2}}, \{\check{\lambda}_{j}^{t^{'}+1+\frac{1}{2}}\}_{j\in \mathcal{N}_{i}})\|\\
&\le \max_{j\in (\mathcal{N}_{i}\cap \mathcal{H})\cup \{i\}} \|\lambda_{j}^{t^{'}+1+\frac{1}{2}}\|\le\frac{\psi}{\theta}, \notag
\end{align}
where the second inequality holds based on \eqref{proof-l2-1}.  Hence, when $t=t^{'}+1$, $\|\lambda_{i}^{t^{'}+1+\frac{1}{2}}\|\le \frac{\psi}{\theta}$ and $\|\lambda_{i}^{t^{'}+1+1}\|\le \frac{\psi}{\theta}$ hold.
\end{IEEEproof}
\end{lemma}

\begin{lemma}
\label{lemma3}
Define a matrix $\Lambda^{t+\frac{1}{2}}=[\cdots,\bm{\lambda}_{i}^{t+\frac{1}{2}},\cdots]\in \mathbb{R}^{H\times d} $ that collects the dual variables $\bm{\lambda}_{i}^{t+\frac{1}{2}}$ of all benign agents $i\in \mathcal{H}$ generated by Algorithm \ref{alg2}. Under Assumption \ref{a1}, we have
\begin{align}
\label{l3-1}
&\|\Lambda^{t+\frac{1}{2}}-\frac{1}{H}\bm{1}\bm{1}^{\top}\Lambda^{t+\frac{1}{2}}
\|^{2}_{F}\\
\le&\frac{1}{1-u}\|\Lambda^{t}-\frac{1}{H}\bm{1}\bm{1}^{\top}\Lambda^{t}
\|^{2}_{F}+\frac{4H^{3}\beta^{2}\psi^{2}}{uM^{2}},\notag
\end{align}
where $u$ is any positive constant in $(0,1)$. If $u=\frac{1}{2}$, this further yields
\begin{align}
\label{l3-2}
&\|\Lambda^{t+\frac{1}{2}}-\frac{1}{H}\bm{1}\bm{1}^{\top}\Lambda^{t+\frac{1}{2}}
\|^{2}_{F}\\
\le&2\|\Lambda^{t}-\frac{1}{H}\bm{1}\bm{1}^{\top}\Lambda^{t}
\|^{2}_{F}+\frac{8H^{3}\beta^{2}\psi^{2}}{M^{2}}. \notag
\end{align}
\begin{IEEEproof}
\label{proof-l3}
According to the update of $\bm{\lambda}_{i}^{t+\frac{1}{2}}$ in Algorithm \ref{alg2}, the fact $\|\Lambda^{t+\frac{1}{2}}-\frac{1}{H}\bm{1}\bm{1}^{\top}\Lambda^{t+\frac{1}{2}}
\|^{2}_{F}=\sum_{i\in \mathcal{H}}\|\bm{\lambda}_{i}^{t+\frac{1}{2}}-\bar{\bm{\lambda}}^{t+\frac{1}{2}}\|^{2}$ and the relationship $\widetilde{G}_{i}^{t}(P_{i})=\frac{H}{M}G_{i}^{t}(P_{i})$, we have
\begin{align}
\label{proof-l3-1}
&\|\Lambda^{t+\frac{1}{2}}-\frac{1}{H}\bm{1}\bm{1}^{\top}\Lambda^{t+\frac{1}{2}}
\|^{2}_{F}\\
=&\sum_{i\in \mathcal{H}}\|\bm{\lambda}_{i}^{t+\frac{1}{2}}-\bar{\bm{\lambda}}^{t+\frac{1}{2}}\|^{2}\notag \\
=&\sum_{i\in \mathcal{H}}\|\lambda_{i}^{t}+\beta\cdot(\widetilde{G}_{i}^{t}(P_{i}^{t})-\theta \lambda_{i}^{t})\notag\\
&-\bar{\lambda}^{t}-\beta\cdot\frac{1}{H}\sum_{i\in \mathcal{H}}(\widetilde{G}_{i}^{t}(P_{i}^{t})-\theta \lambda_{i}^{t})\|^{2}\notag\\
=&\sum_{i\in \mathcal{H}}\|(1-\beta \theta )\cdot(\lambda_{i}^{t}-\bar{\lambda}^{t})+\beta\cdot[\widetilde{G}_{i}^{t}(P_{i}^{t})-\frac{1}{H}\sum_{i\in \mathcal{H}}\widetilde{G}_{i}^{t}(P_{i}^{t})]\|^{2}\notag\\
\le& \frac{(1-\beta \theta)^{2}}{1-u}\sum_{i\in \mathcal{H}}\|\lambda_{i}^{t}-\bar{\lambda}^{t}\|^{2}\notag\\
&+\frac{\beta^{2}}{u}\sum_{i\in \mathcal{H}}\|\widetilde{G}_{i}^{t}(P_{i}^{t})-\frac{1}{H}\sum_{i\in \mathcal{H}}\widetilde{G}_{i}^{t}(P_{i}^{t})\|^{2}\notag\\
\le& \frac{(1-\beta \theta)^{2}}{1-u}\sum_{i\in \mathcal{H}}\|\lambda_{i}^{t}-\bar{\lambda}^{t}\|^{2}+\frac{2\beta^{2}}{u}\sum_{i\in \mathcal{H}}\|\widetilde{G}_{i}^{t}(P_{i}^{t})\|^{2}\notag\\
&+\frac{2\beta^{2} H}{u}\|\frac{1}{H}\sum_{i\in \mathcal{H}}\widetilde{G}_{i}^{t}(P_{i}^{t})\|^{2}\notag\\
\le& \frac{(1-\beta \theta)^{2}}{1-u}\sum_{i\in \mathcal{H}}\|\lambda_{i}^{t}-\bar{\lambda}^{t}\|^{2}+\frac{2\beta^{2}H^{2}}{M^{2}u}\sum_{i\in \mathcal{H}}\|G_{i}^{t}(P_{i}^{t})\|^{2}\notag\\
&+\frac{2\beta^{2}H^{2}}{M^{2}u}\sum_{i\in \mathcal{H}}\|G_{i}^{t}(P_{i}^{t})\|^{2}\notag\\
\le& \frac{(1-\beta \theta)^{2}}{1-u}\sum_{i\in \mathcal{H}}\|\lambda_{i}^{t}-\bar{\lambda}^{t}\|^{2}+\frac{4H^{3}\beta^{2}\psi^{2}}{uM^{2}}\notag\\
\le &\frac{1}{1-u}\sum_{i\in \mathcal{H}}\|\lambda_{i}^{t}-\bar{\lambda}^{t}\|^{2}+\frac{4H^{3}\beta^{2}\psi^{2}}{uM^{2}}\notag\\
=&\frac{1}{1-u}\|\Lambda^{t}-\frac{1}{H}\bm{1}\bm{1}^{\top}\Lambda^{t}
\|^{2}_{F}+\frac{4H^{3}\beta^{2}\psi^{2}}{uM^{2}},\notag
\end{align}
where $u$ is any positive constant in $(0,1)$. To derive the first and second inequalities, we use $\|a+b\|^{2}\le \frac{1}{1-u}\|a\|^{2}+\frac{1}{u}\|b\|^{2} (u\in (0,1))$. The third inequality holds, since the mean inequality $\|a_{1}+a_{2}+\cdots+a_{H}\|^{2}\le H (\|a_{1}\|^{2}+\|a_{2}\|^{2}+\cdots+\|a_{H}\|^{2})$. Based on Assumption \ref{a1}, the forth inequality holds.
\end{IEEEproof}
\end{lemma}

\begin{lemma}
\label{lemma4}
Define a matrix $\Lambda^{t+1}=[\cdots,\bm{\lambda}_{i}^{t+1},\cdots]\in \mathbb{R}^{H\times d} $ that collects the dual variables $\bm{\lambda}_{i}^{t+1}$ of all benign agents $i\in \mathcal{H}$ generated by Algorithm \ref{alg2}. Suppose that the robust aggregation rule $AGG(\cdot)$ satisfies Property \ref{d1}. Under Assumptions \ref{a1} and \ref{a3}, if the contraction constant $\rho$ satisfies
$\rho<\frac{(1-\kappa)^{2}}{64H}$, we have
\begin{align}
\label{l4}
\|\Lambda^{t+1}-\frac{1}{H}\bm{1}\bm{1}^{\top}\Lambda^{t+1}
\|^{2}_{F}\le \frac{4H^{3}\beta^{2}\psi^{2}}{\epsilon^{3}M^{2}},
\end{align}
where $\epsilon:=1-\kappa-8\sqrt{\rho H}$.
\begin{IEEEproof}
For any positive constant $w\in (0,1)$, we have
\begin{align}
\label{proof-l4-1}
 &\|\Lambda^{t+1}-\frac{1}{H}\bm{1}\bm{1}^{\top}\Lambda^{t+1}
\|^{2}_{F} \\
=& \|\Lambda^{t+1}-\frac{1}{H}\bm{1}\bm{1}^{\top}\Lambda^{t+1}+E\Lambda^{t+\frac{1}{2}}-E\Lambda^{t+\frac{1}{2}}\notag\\
&+\frac{1}{H}\bm{1}\bm{1}^{\top}E\Lambda^{t+\frac{1}{2}}-\frac{1}{H}\bm{1}\bm{1}^{\top}E\Lambda^{t+\frac{1}{2}}\|^{2}_{F} \notag\\
\le& \underbrace{\frac{1}{1-w}\|E\Lambda^{t+\frac{1}{2}}-\frac{1}{H}\bm{1}\bm{1}^{\top}E\Lambda^{t+\frac{1}{2}}\|^{2}_{F}}_{A_{1}}+\underbrace{\frac{2}{w}\|\Lambda^{t+1}-E\Lambda^{t+\frac{1}{2}}\|^{2}_{F}}_{A_{2}}\notag\\
&+\underbrace{\frac{2}{w}\|\frac{1}{H}\bm{1}\bm{1}^{\top}\Lambda^{t+1}-\frac{1}{H}\bm{1}\bm{1}^{\top}E\Lambda^{t+\frac{1}{2}}\|^{2}_{F}}_{A_{3}}. \notag
\end{align}
Next, we analyze $A_{1}$, $A_{2}$ and $A_{3}$ in turn.

\noindent\textbf{Bounding $A_{1}$:}
According to Assumption \ref{a3}, we have
\begin{align}
\label{proof-l4-2}
A_{1}=&\frac{1}{1-w}\|E\Lambda^{t+\frac{1}{2}}-\frac{1}{H}\bm{1}\bm{1}^{\top}E\Lambda^{t+\frac{1}{2}}\|^{2}_{F}\\
=&\frac{1}{1-w}\|(I-\frac{1}{H}\bm{1}\bm{1}^{\top})E\Lambda^{t+\frac{1}{2}}\|^{2}_{F} \notag\\
=&\frac{1}{1-w}\|(I-\frac{1}{H}\bm{1}\bm{1}^{\top})E(I-\frac{1}{H}\bm{1}\bm{1}^{\top})\Lambda^{t+\frac{1}{2}}\|^{2}_{F}\notag\\
\le&\frac{1}{1-w}\|(I-\frac{1}{H}\bm{1}\bm{1}^{\top})E\|^{2}\|(I-\frac{1}{H}\bm{1}\bm{1}^{\top})\Lambda^{t+\frac{1}{2}}\|^{2}_{F} \notag\\
=&\frac{\kappa}{1-w}\|\Lambda^{t+\frac{1}{2}}-\frac{1}{H}\bm{1}\bm{1}^{\top}\Lambda^{t+\frac{1}{2}}
\|^{2}_{F},  \notag
\end{align}
where the last inequality holds because of Assumption \ref{a3} and the fact that $\|AB\|_{F}^{2}\le \|A\|^{2}\|B\|_{F}^{2}$.

\noindent\textbf{Bounding $A_{2}$:} According to the update of $\bm{\lambda}_{i}^{t+1}$ in Algorithm \ref{alg2}, Property \ref{d1} and Lemma \ref{lemma-AGG-ARC-definition1}, we have
\begin{align}
\label{proof-l4-3}
A_{2}=&\frac{2}{w}\|\Lambda^{t+1}-E\Lambda^{t+\frac{1}{2}}\|^{2}_{F}\\
=&\frac{2}{u}\sum_{i\in\mathcal{H} }\|\bm{\lambda}_{i}^{t+1}-\bar{\bm{\lambda}}_{i}^{t+\frac{1}{2}}\|^{2} \notag\\
=&\frac{2}{w}\sum_{i\in\mathcal{H} }\|AGG (\bm{\lambda}_{i}^{t+\frac{1}{2}},\{\check{\bm{\lambda}}_{j}^{t+\frac{1}{2}}\}_{j\in \mathcal{N}_{i}})-\bar{\bm{\lambda}}_{i}^{t+\frac{1}{2}}\|^{2} \notag
\end{align}
\begin{align}
\le&\frac{2}{w}\sum_{i\in \mathcal{H}} \rho \cdot \sum_{j\in \mathcal{N}_{i}\cap\mathcal{H}\cup i}e_{ij}\|\bm{\lambda}_{j}^{t+\frac{1}{2}}-\bar{\bm{\lambda}}_{i}^{t+\frac{1}{2}}\|^{2} \notag\\
=& \frac{2}{w}\sum_{i\in \mathcal{H}}\rho \cdot \max_{j\in \mathcal{N}_{i}\cap\mathcal{H}\cup i}\|\bm{\lambda}_{j}^{t+\frac{1}{2}}-\bar{\bm{\lambda}}^{t+\frac{1}{2}}+\bar{\bm{\lambda}}^{t+\frac{1}{2}}-\bar{\bm{\lambda}}_{i}^{t+\frac{1}{2}}\|^{2} \notag\\
\le& \frac{4}{w}\sum_{i\in \mathcal{H}} \rho \cdot [\max_{j\in \mathcal{N}_{i}\cap\mathcal{H}\cup i}\|\bm{\lambda}_{j}^{t+\frac{1}{2}}-\bar{\bm{\lambda}}^{t+\frac{1}{2}}\|^{2}+\|\bar{\bm{\lambda}}^{t+\frac{1}{2}}-\bar{\bm{\lambda}}_{i}^{t+\frac{1}{2}}\|^{2}]\notag\\
\le&\frac{4}{w}\sum_{i\in \mathcal{H}} \rho \cdot [\max_{i\in \mathcal{H}}\|\bm{\lambda}_{i}^{t+\frac{1}{2}}-\bar{\bm{\lambda}}^{t+\frac{1}{2}}\|^{2}+\max_{i\in \mathcal{H}}\|\bar{\bm{\lambda}}^{t+\frac{1}{2}}-\bar{\bm{\lambda}}_{i}^{t+\frac{1}{2}}\|^{2}]\notag\\
=&\frac{8 \rho H}{w}\cdot\max_{i\in \mathcal{H}}\|\bm{\lambda}_{i}^{t+\frac{1}{2}}-\bar{\bm{\lambda}}^{t+\frac{1}{2}}\|^{2}\notag\\
\le&\frac{8 \rho H}{w} \|\Lambda^{t+\frac{1}{2}}-\frac{1}{H}\bm{1}\bm{1}^{\top}\Lambda^{t+\frac{1}{2}}
\|^{2}_{F}\notag,
\end{align}
where the last inequality holds as $\max_{i\in \mathcal{H}}\|\bm{\lambda}_{i}^{t+\frac{1}{2}}-\bar{\bm{\lambda}}^{t+\frac{1}{2}}\|^{2}\le \|\Lambda^{t+\frac{1}{2}}-\frac{1}{H}\bm{1}\bm{1}^{\top}\Lambda^{t+\frac{1}{2}}\|^{2}_{F}$.

\noindent\textbf{Bounding $A_{3}$:}
Likewise, according to the update of $\bm{\lambda}_{i}^{t+1}$ in Algorithm \ref{alg2} and Property \ref{d1}, we have
\begin{align}
\label{proof-l4-4}
A_{3}=&\frac{2}{w}\|\frac{1}{H}\bm{1}\bm{1}^{\top}\Lambda^{t+1}-\frac{1}{H}\bm{1}\bm{1}^{\top}E\Lambda^{t+\frac{1}{2}}\|^{2}_{F}\\
=&\frac{2}{w}\|\frac{1}{H}\bm{1}\bm{1}^{\top}(\Lambda^{t+1}-E\Lambda^{t+\frac{1}{2}})\|^{2}_{F} \notag\\
\le & \frac{2}{w} \|\frac{1}{H}\bm{1}\bm{1}^{\top}\|^{2}_{F}\|\Lambda^{t+1}-E\Lambda^{t+\frac{1}{2}}\|^{2}_{F}\notag\\
=&\frac{2}{w}\|\Lambda^{t+1}-E\Lambda^{t+\frac{1}{2}}\|^{2}_{F}\notag\\
\le &\frac{8\rho H}{w}\max_{i\in \mathcal{H}}\|\bm{\lambda}_{i}^{t+\frac{1}{2}}-\bar{\bm{\lambda}}^{t+\frac{1}{2}}\|^{2}\notag\\
\le&\frac{8\rho H}{w} \|\Lambda^{t+\frac{1}{2}}-\frac{1}{H}\bm{1}\bm{1}^{\top}\Lambda^{t+\frac{1}{2}}
\|^{2}_{F}.\notag
\end{align}
To derive the last equality, we use the fact $\|\frac{1}{H}\bm{1}\bm{1}^{\top}\|^{2}_{F}=1$. From the last equality to the last inequality, we use the same technique in deriving \eqref{proof-l4-3}.

Therefore, substituting \eqref{proof-l4-2}, \eqref{proof-l4-3} and \eqref{proof-l4-4} into \eqref{proof-l4-1} and rearranging the terms, we obtain
\begin{align}
\label{proof-l4-5}
&\|\Lambda^{t+1}-\frac{1}{H}\bm{1}\bm{1}^{\top}\Lambda^{t+1}
\|^{2}_{F}\\
\le&(\frac{\kappa}{1-w}+\frac{16\rho H}{w})\|\Lambda^{t+\frac{1}{2}}-\frac{1}{H}\bm{1}\bm{1}^{\top}\Lambda^{t+\frac{1}{2}}
\|^{2}_{F}. \notag
\end{align}
Substituting \eqref{l3-1} in Lemma \ref{lemma3} into \eqref{proof-l4-5} and rearranging the terms, we obtain
\begin{align}
\label{proof-l4-6}
&\|\Lambda^{t+1}-\frac{1}{H}\bm{1}\bm{1}^{\top}\Lambda^{t+1}
\|^{2}_{F}\\
\le&(\frac{\kappa}{1-w}+\frac{16\rho H}{w})\cdot \frac{1}{1-u}\|\Lambda^{t}-\frac{1}{H}\bm{1}\bm{1}^{\top}\Lambda^{t}
\|^{2}_{F}\notag\\
&+(\frac{\kappa}{1-w}+\frac{16\rho H}{w})\cdot\frac{4H^{3}\beta^{2}\psi^{2}}{uM^{2}}. \notag
\end{align}

By setting the constant $w=4\sqrt{\rho H}\le 1-\kappa$, $\frac{\kappa}{1-w}\le \kappa+w$ holds. Therefore, we can rewrite \eqref{proof-l4-6} as
\begin{align}
\label{proof-l4-7}
&\|\Lambda^{t+1}-\frac{1}{H}\bm{1}\bm{1}^{\top}\Lambda^{t+1}
\|^{2}_{F} \\
\le&(\kappa+8\sqrt{\rho H})\cdot \frac{1}{1-u}\|\Lambda^{t}-\frac{1}{H}\bm{1}\bm{1}^{\top}\Lambda^{t}
\|^{2}_{F}\notag\\
&+(\kappa+8\sqrt{\rho H})\cdot\frac{4H^{3}\beta^{2}\psi^{2}}{uM^{2}} \notag\\
=&(1-\epsilon) \cdot \frac{1}{1-u}\|\Lambda^{t}-\frac{1}{H}\bm{1}\bm{1}^{\top}\Lambda^{t}
\|^{2}_{F}+(1-\epsilon)\cdot\frac{4H^{3}\beta^{2}\psi^{2}}{uM^{2}},\notag
\end{align}
where $\epsilon:=1-\kappa-8\sqrt{\rho H}$. The parameter $\rho$ should satisfy $\rho<\frac{(1-\kappa)^{2}}{64H}$ to guarantee $\epsilon>0$.

Set $u=\frac{\epsilon}{1+\epsilon}$. Therefore, we have $\frac{1}{1-u}= 1+\epsilon$. In consequence, \eqref{proof-l4-7} can be rewritten as
\begin{align}
\label{proof-l4-8}
&\|\Lambda^{t+1}-\frac{1}{H}\bm{1}\bm{1}^{\top}\Lambda^{t+1}
\|^{2}_{F}\notag \\
\le&(1-\epsilon^{2})\cdot \|\Lambda^{t}-\frac{1}{H}\bm{1}\bm{1}^{\top}\Lambda^{t}
\|^{2}_{F}+\frac{1-\epsilon^{2}}{\epsilon}\cdot \frac{4H^{3}\beta^{2}\psi^{2}}{M^{2}} \notag\\
\le&(1-\epsilon^{2})\cdot \|\Lambda^{t}-\frac{1}{H}\bm{1}\bm{1}^{\top}\Lambda^{t}
\|^{2}_{F}+\frac{4H^{3}\beta^{2}\psi^{2}}{\epsilon M^{2}}.
\end{align}

Under the conditions $\rho<\frac{(1-\kappa)^{2}}{64H}$ and $\epsilon\in(0,1)$, we write \eqref{proof-l4-8} recursively to yield
\begin{align}
\label{proof-l4-9}
&\|\Lambda^{t+1}-\frac{1}{H}\bm{1}\bm{1}^{\top}\Lambda^{t+1}
\|^{2}_{F}\\
\le& (1-\epsilon^{2})^{t+1}\|\Lambda^{0}-\frac{1}{H}\bm{1}\bm{1}^{\top}\Lambda^{0}
\|^{2}_{F}+\sum_{l=0}^{t}(1-\epsilon^{2})^{t-l}\cdot\frac{4H^{3}\beta^{2}\psi^{2}}{\epsilon M^{2}}\notag.
\end{align}

With the same initialization $\bm{\lambda}_{i}^{0}$ for all benign agents $i\in \mathcal{H}$, we can rewrite \eqref{proof-l4-9} as
\begin{align}
\label{proof-l4-10}
&\|\Lambda^{t+1}-\frac{1}{H}\bm{1}\bm{1}^{\top}\Lambda^{t+1}
\|^{2}_{F}\\
\le& \sum_{l=0}^{t}(1-\epsilon^{2})^{t-l}\cdot\frac{4H^{3}\beta^{2}\psi^{2}}{\epsilon M^{2}}\le  \frac{4H^{3}\beta^{2}\psi^{2}}{\epsilon^{3}M^{2}}.\notag
\end{align}
\end{IEEEproof}
\end{lemma}

\begin{lemma}
\label{lemma5}
Suppose that the robust aggregation rule $AGG(\cdot)$ satisfies Property \ref{d1}. For any agent $i\in \mathcal{H}$ and $t\in [1,\cdots,T]$, consider $\lambda_{i}^{t+1}$ generated by Algorithm \ref{alg2}. Under Assumptions \ref{a1} and \ref{a3}, we have
\begin{align}
\label{l5}
&\frac{1-\beta\theta}{2\beta}\cdot\Delta^{t}\le  \frac{H}{M} \sum_{i\in \mathcal{H}} \left \langle  \bar{\lambda}^{t}  , G_{i}^{t}(P_{i}^{t})  \right \rangle - \frac{H}{M} \sum_{i\in \mathcal{H}} \left \langle  \lambda  , G_{i}^{t}(P_{i}^{t})  \right \rangle \notag\\
&+\frac{\theta H}{2} \|\lambda\|^{2}+\frac{H^{3}\psi^{2}\beta}{M^{2}}+\psi^{2}H \beta \\
&+(1+\frac{1}{\epsilon^{3}})\cdot (4\rho^{2}H+\chi^{2}) \cdot \frac{8 H^{3}\psi^{2}}{M^{2}\theta},\notag
\end{align}
where $\Delta^{t}:=H \|\bar{\lambda}^{t+1}-\lambda\|^{2}-H\|\bar{\lambda}^{t}-\lambda\|^{2}$.

\begin{IEEEproof}
According to the update of $\lambda_{i}^{t+1} $ in Algorithm \ref{alg2}, we have
\begin{align}
\label{proof-l5-1}
&H \|\bar{\lambda}^{t+1}-\lambda\|^{2} \\
=&\|\frac{1}{H}\sum_{i\in \mathcal{H}} AGG (\lambda_{i}^{t+\frac{1}{2}}, \{\check{\lambda}_{j}^{t+\frac{1}{2}}\}_{j\in \mathcal{N}_{i}})-\lambda\|^{2}\notag
\end{align}
\begin{align}
=&H\|\frac{1}{H}\sum_{i\in \mathcal{H}} AGG (\lambda_{i}^{t+\frac{1}{2}}, \{\check{\lambda}_{j}^{t+\frac{1}{2}}\}_{j\in \mathcal{N}_{i}})-\frac{1}{H}\sum_{i\in \mathcal{H}} \bar{\lambda}_{i}^{t+\frac{1}{2}}\notag\\
&+\frac{1}{H}\sum_{i\in \mathcal{H}} \bar{\lambda}_{i}^{t+\frac{1}{2}}-\bar{\lambda}^{t+\frac{1}{2}}+\bar{\lambda}^{t+\frac{1}{2}}-\lambda\|^{2}\notag\\
\le & \underbrace{\frac{2H}{u} \|\frac{1}{H}\sum_{i\in \mathcal{H}} AGG(\lambda_{i}^{t+\frac{1}{2}}, \{\check{\lambda}_{j}^{t+\frac{1}{2}}\}_{j\in \mathcal{N}_{i}})-\frac{1}{H}\sum_{i\in \mathcal{H}} \bar{\lambda}_{i}^{t+\frac{1}{2}}\|^{2}}_{A_{1}}\notag\\
&+\underbrace{\frac{2H}{u} \|\frac{1}{H}\sum_{i\in \mathcal{H}}  \bar{\lambda}_{i}^{t+\frac{1}{2}}-\bar{\lambda}^{t+\frac{1}{2}}\|^{2}}_{A_{2}}+\underbrace{\frac{H}{1-u}\|\bar{\lambda}^{t+\frac{1}{2}}-\lambda \|^{2}}_{A_{3}},\notag
\end{align}
where $u$ is any positive constant in $(0,1)$.

Next, we analyze $A_{1}$, $A_{2}$ and $A_{3}$ in turn.

\noindent\textbf{Bounding $A_{1}$:}
Similar to the derivation of \eqref{proof-l4-3} in Lemma \ref{lemma4}, we obtain
\begin{align}
\label{proof-l5-2}
A_{1}=& \frac{2H}{u} \|\frac{1}{H}\sum_{i\in \mathcal{H}}  AGG(\lambda_{i}^{t+\frac{1}{2}}, \{\check{\lambda}_{j}^{t+\frac{1}{2}}\}_{j\in \mathcal{N}_{i}})-\frac{1}{H}\sum_{i\in \mathcal{H}} \bar{\lambda}_{i}^{t+\frac{1}{2}}\|^{2}\notag \\
\le & \frac{2}{u} \sum_{i\in \mathcal{H}} \|AGG (\lambda_{i}^{t+\frac{1}{2}}, \{\check{\lambda}_{j}^{t+\frac{1}{2}}\}_{j\in \mathcal{N}_{i}})-\bar{\lambda}_{i}^{t+\frac{1}{2}}\|^{2}\notag\\
\le & \frac{8\rho H}{u} \cdot \|\Lambda^{t+\frac{1}{2}}-\frac{1}{H}\bm{1}\bm{1}^{\top}\Lambda^{t+\frac{1}{2}}
\|^{2}_{F}.
\end{align}

\noindent\textbf{Bounding $A_{2}$:} It holds that
\begin{align}
\label{proof-l5-4}
A_{2}=&\frac{2H}{u} \|\frac{1}{H}\sum_{i\in \mathcal{H}}  \bar{\lambda}_{i}^{t+\frac{1}{2}}-\bar{\lambda}^{t+\frac{1}{2}}\|^{2}\\
=& \frac{2H}{u}\|\frac{1}{H}\bm{1}^{\top}E\Lambda^{t+\frac{1}{2}}-\frac{1}{H}\bm{1}^{\top}\Lambda^{t+\frac{1}{2}}\|^{2}\notag\\
=& \frac{2}{u H} \cdot \|(\bm{1}^{\top}E-\bm{1}^{\top})(\Lambda^{t+\frac{1}{2}}-\frac{1}{H}\bm{1}\bm{1}^{\top}\Lambda^{t+\frac{1}{2}}
)\|^{2}\notag\\
\le & \frac{2}{uH} \cdot \|E^{\top}\bm{1}-\bm{1}\|^{2}\|\Lambda^{t+\frac{1}{2}}-\frac{1}{H}\bm{1}\bm{1}^{\top}\Lambda^{t+\frac{1}{2}}
\|^{2}.\notag
\end{align}
To derive the last equality, we use Property \ref{d1} that the virtual weight matrix $E$ is row stochastic.

Define $\chi=\frac{1}{H}\|E^{\top}\bm{1}-\bm{1}\|^{2}$ to quantify how non-column stochastic the virtual weight matrix $E$ is. Applying the fact $\|\cdot\|^{2}\le \|\cdot\|_{F}^{2}$ to the right-hand side of \eqref{proof-l5-4}, we have
\begin{align}
\label{proof-l5-5}
A_{2}\le \frac{2\chi}{u} \cdot \|\Lambda^{t+\frac{1}{2}}-\frac{1}{H}\bm{1}\bm{1}^{\top}\Lambda^{t+\frac{1}{2}}
\|_{F}^{2}.
\end{align}

\noindent\textbf{Bounding $A_{3}$:}
According to the update of $\lambda_{i}^{t+\frac{1}{2}}$ in Algorithm \ref{alg2} and the relationship $\widetilde{G}_{i}^{t}(P_{i})=\frac{H}{M}G_{i}^{t}(P_{i})$, we have
\begin{align}
\label{proof-l5-7}
A_{3}=&\frac{H}{1-u}\|\bar{\lambda}^{t+\frac{1}{2}}-\lambda \|^{2}\\
=& \frac{H}{1-u}\|\frac{1}{H}\sum_{i\in \mathcal{H}}[\lambda_{i}^{t}-\lambda+\beta\cdot(\widetilde{G}_{i}^{t}(P_{i}^{t})-\theta \lambda_{i}^{t})] \|^{2}\notag
\end{align}
\begin{align}
=&\frac{H}{1-u}\|\bar{\lambda}^{t}-\lambda+\frac{\beta}{H}\sum_{i\in \mathcal{H}}(\widetilde{G}_{i}^{t}(P_{i}^{t})-\theta \lambda_{i}^{t}) \|^{2}\notag\\
=& \frac{H}{1-u}\|\bar{\lambda}^{t}-\lambda\|^{2}+
\frac{\beta^{2}}{1-u}\sum_{i\in \mathcal{H}}\|\widetilde{G}_{i}^{t}(P_{i}^{t})-\theta \lambda_{i}^{t} \|^{2}\notag\\
&+\frac{2H\beta}{1-u} \left \langle  \bar{\lambda}^{t}-\lambda  , \frac{1}{H}\sum_{i\in \mathcal{H}}(\widetilde{G}_{i}^{t}(P_{i}^{t})-\theta \lambda_{i}^{t})   \right \rangle \notag\\
=& \frac{H}{1-u}\|\bar{\lambda}^{t}-\lambda\|^{2}+\frac{2H\beta}{1-u} \left \langle  \bar{\lambda}^{t}  , \frac{1}{H}\sum_{i\in \mathcal{H}}\widetilde{G}_{i}^{t}(P_{i}^{t})  \right \rangle\notag\\
&- \frac{2\beta}{1-u} \sum_{i\in \mathcal{H}} \left \langle  \lambda  , \widetilde{G}_{i}^{t}(P_{i}^{t})  \right \rangle-\frac{2H\beta\theta }{1-u}\left \langle \bar{\lambda}^{t}-\lambda  , \bar{\lambda}^{t}  \right \rangle \notag\\
&+\frac{\beta^{2}}{1-u}\sum_{i\in \mathcal{H}}\|\widetilde{G}_{i}^{t}(P_{i}^{t})-\theta \lambda_{i}^{t} \|^{2}\notag\\
=& \frac{H}{1-u}\|\bar{\lambda}^{t}-\lambda\|^{2}+ \frac{2\beta H}{(1-u)M} \sum_{i\in \mathcal{H}} \left \langle  \bar{\lambda}^{t}  , G_{i}^{t}(P_{i}^{t})  \right \rangle\notag\\
&+\underbrace{\frac{\beta^{2}}{1-u}\sum_{i\in \mathcal{H}}\|\frac{H}{M}G_{i}^{t}(P_{i}^{t})-\theta \lambda_{i}^{t} \|^{2}}_{A_{3-1}}\underbrace{-\frac{2H\beta\theta}{1-u} \left \langle \bar{\lambda}^{t}-\lambda  , \bar{\lambda}^{t} \right \rangle  }_{A_{3-2}} \notag\\
&- \frac{2\beta H}{(1-u)M} \sum_{i\in \mathcal{H}} \left \langle  \lambda  , G_{i}^{t}(P_{i}^{t})  \right \rangle.\notag
\end{align}

Next, we analyze the terms $A_{3-1}$ and $A_{3-2}$ in turn. Based on the mean inequality $\|a+b\|^{2}\le 2\|a\|^{2}+2\|b\|^{2}$, we have
\begin{align}
\label{proof-l5-8}
A_{3-1}=&\frac{\beta^{2}}{1-u}\sum_{i\in \mathcal{H}}\|\frac{H}{M}G_{i}^{t}(P_{i}^{t})-\theta \lambda_{i}^{t} \|^{2}\\
\le & \frac{2\beta^{2}}{1-u}\sum_{i\in \mathcal{H}} \|\frac{H}{M}G_{i}^{t}(P_{i}^{t})\|^{2}+\frac{2\beta^{2}\theta^{2}}{1-u}\sum_{i\in \mathcal{H}}\|\lambda_{i}^{t} \|^{2}\notag\\
\le & \frac{2H^{3}\beta^{2}\psi^{2}}{M^{2}(1-u)}+\frac{2\beta^{2}\psi^{2}H}{1-u},\notag
\end{align}
where the last inequality holds, since Assumption \ref{a1} and the conclusions in Lemma \ref{lemma2}.

Based on the inequality $-2\left \langle a-b  , a  \right \rangle \le \|b\|^{2}-\|a-b\|^{2}$, we obtain
\begin{align}
\label{proof-l5-9}
A_{3-2}=&-\frac{2H\beta\theta}{1-u} \left \langle \bar{\lambda}^{t}-\lambda  , \bar{\lambda}^{t}  \right \rangle\\
\le& \frac{H \beta \theta}{1-u} [\|\lambda\|^{2}-\|\bar{\lambda}^{t}-\lambda\|^{2}]\notag\\
=& \frac{\beta\theta H}{1-u} \|\lambda\|^{2}-\frac{\beta\theta H}{1-u}\|\bar{\lambda}^{t}-\lambda\|^{2}.\notag
\end{align}

Substituting \eqref{proof-l5-8} and \eqref{proof-l5-9} into \eqref{proof-l5-7} and rearranging
the terms, we have
\begin{align}
\label{proof-l5-10}
A_{3}\le&\frac{(1-\beta \theta)H}{1-u}\|\bar{\lambda}^{t}-\lambda\|^{2}+ \frac{2\beta H}{(1-u)M} \sum_{i\in \mathcal{H}} \left \langle  \bar{\lambda}^{t}  , G_{i}^{t}(P_{i}^{t})  \right \rangle\notag\\
&- \frac{2\beta H}{(1-u)M} \sum_{i\in \mathcal{H}} \left \langle  \lambda  , G_{i}^{t}(P_{i}^{t})  \right \rangle+\frac{\beta\theta H}{1-u} \|\lambda\|^{2}\notag\\
&+\frac{2H^{3}\beta^{2}\psi^{2}}{M^{2}(1-u)}+\frac{2\beta^{2}\psi^{2}H}{1-u}.
\end{align}

Substituting \eqref{proof-l5-2}, \eqref{proof-l5-5} and \eqref{proof-l5-10} into \eqref{proof-l5-1} and rearranging
the terms, we have
\begin{align}
\label{proof-l5-11}
&H \|\bar{\lambda}^{t+1}-\lambda\|^{2}
\\
\le&\frac{(1-\beta \theta)H}{1-u}\|\bar{\lambda}^{t}-\lambda\|^{2}+ \frac{2\beta H}{(1-u)M} \sum_{i\in \mathcal{H}} \left \langle  \bar{\lambda}^{t}  , G_{i}^{t}(P_{i}^{t})  \right \rangle\notag\\
&- \frac{2\beta H}{(1-u)M} \sum_{i\in \mathcal{H}} \left \langle  \lambda  , G_{i}^{t}(P_{i}^{t})  \right \rangle+\frac{\beta\theta H}{1-u} \|\lambda\|^{2}\notag\\
&+\frac{2H^{3}\beta^{2}\psi^{2}}{M^{2}(1-u)}+\frac{2\beta^{2}\psi^{2}H}{1-u}\notag\\
&+\frac{2}{u}\cdot (4\rho H+\chi) \cdot \|\Lambda^{t+\frac{1}{2}}-\frac{1}{H}\bm{1}\bm{1}^{\top}\Lambda^{t+\frac{1}{2}}
\|_{F}^{2}.\notag
\end{align}

Based on Lemma \ref{lemma3} and Lemma \ref{lemma4}, we obtain
\begin{align}
\label{proof-l5-12}
& H \|\bar{\lambda}^{t+1}-\lambda\|^{2}
-\frac{(1-\beta \theta)}{1-u}H\|\bar{\lambda}^{t}-\lambda\|^{2}\\
\le& \frac{2\beta H}{(1-u)M} \sum_{i\in \mathcal{H}} \left \langle  \bar{\lambda}^{t}  , G_{i}^{t}(P_{i}^{t})  \right \rangle - \frac{2\beta H}{(1-u)M} \sum_{i\in \mathcal{H}} \left \langle  \lambda  , G_{i}^{t}(P_{i}^{t})  \right \rangle \notag\\
&+\frac{\beta\theta H}{1-u} \|\lambda\|^{2}+\frac{2H^{3}\beta^{2}\psi^{2}}{M^{2}(1-u)}+\frac{2\beta^{2}\psi^{2}H}{1-u}\notag\\
&+\frac{16}{u}\cdot(1+\frac{1}{\epsilon^{3}})\cdot (4\rho H+\chi) \cdot \frac{H^{3}\beta^{2}\psi^{2}}{M^{2}}.\notag
\end{align}

Setting $u=\beta\theta$, we rewritten \eqref{proof-l5-12} as
\begin{align}
\label{proof-l5-13}
&H \|\bar{\lambda}^{t+1}-\lambda\|^{2}
-H\|\bar{\lambda}^{t}-\lambda\|^{2}\\
\le& \frac{2\beta H}{(1-\beta\theta)M} \sum_{i\in \mathcal{H}} \left \langle  \bar{\lambda}^{t}  , G_{i}^{t}(P_{i}^{t})  \right \rangle - \frac{2\beta H}{(1-\beta\theta)M} \sum_{i\in \mathcal{H}} \left \langle  \lambda  , G_{i}^{t}(P_{i}^{t})  \right \rangle \notag\\
&+\frac{\beta\theta H}{1-\beta\theta} \|\lambda\|^{2}+\frac{2H^{3}\beta^{2}\psi^{2}}{M^{2}(1-\beta\theta)}+\frac{2\beta^{2}\psi^{2}H}{1-\beta\theta}\notag\\
&+\frac{16}{\beta\theta}\cdot(1+\frac{1}{\epsilon^{3}})\cdot (4\rho H+\chi) \cdot \frac{H^{3}\beta^{2}\psi^{2}}{M^{2}}.\notag
\end{align}

Defining $\Delta^{t}:=H \|\bar{\lambda}^{t+1}-\lambda\|^{2}
-H\|\bar{\lambda}^{t}-\lambda\|^{2}$ and multiplying both sides of \eqref{proof-l5-13} by $\frac{1-\beta\theta}{2\beta}$, we have
\begin{align}
\label{proof-l5-14}
\frac{1-\beta\theta}{2\beta}\cdot\Delta^{t}\le&  \frac{H}{M} \sum_{i\in \mathcal{H}} \left \langle  \bar{\lambda}^{t}  , G_{i}^{t}(P_{i}^{t})  \right \rangle - \frac{H}{M} \sum_{i\in \mathcal{H}} \left \langle  \lambda  , G_{i}^{t}(P_{i}^{t})  \right \rangle \notag\\
&+\frac{\theta H}{2} \|\lambda\|^{2}+\frac{H^{3}\psi^{2}\beta}{M^{2}}+\psi^{2}H \beta \notag\\
&+(1+\frac{1}{\epsilon^{3}})\cdot (4\rho H+\chi) \cdot \frac{8 H^{3}\psi^{2}}{M^{2}\theta}.
\end{align}
\end{IEEEproof}
\end{lemma}

\subsection{Violation of Property \ref{d2} for Existing Robust Aggregation Rules}\label{aa2}
\noindent\textbf{Robust Aggregation Rule $CTM(\cdot)$.} We provide a counter-example to demonstrate that the robust aggregation rule $CTM(\cdot)$ does not satisfy Property \ref{d2}.
Specifically, consider a benign agent $i$ whose local dual variable is given by \[\lambda_{i} = [1,0]^\top.\]
Agent $i$ receives dual variables from three neighbors $\mathcal{N}_{i}=\{j_{1},j_{2},j_{3}\} $.
The received dual variables are
\[\check{\lambda}_{j_{1}} = [1,0]^\top,
\check{\lambda}_{j_{2}} = [0,1]^\top,
\check{\lambda}_{j_{3}} = [100,100]^\top.\]
Among them, $\check{\lambda}_{j_{3}}$ is from a Byzantine neighbor, while $\check{\lambda}_{j_{1}}$ and $\check{\lambda}_{j_{2}}$ are from benign neighbors. All benign dual variables, including $\lambda_i$, have a norm of $1$, such that
\[\max_{j\in\mathcal{N}_i\cap\mathcal{H}\cup\{i\}}\|\lambda_j\|=\|\lambda_{i}\| = \|\check{\lambda}_{j_{1}}\| = \|\check{\lambda}_{j_{2}}\| = 1.\]

The $CTM(\cdot)$ operator performs trimmed mean aggregation in a coordinate-wise manner.
Given four total inputs (including $\lambda_i$) and a Byzantine upper bound $b_i = 1$,
$CTM(\cdot)$ removes the $b_i$ largest and $b_i$ smallest values in each coordinate
among the received messages $\{\check{\lambda}_j\}_{j\in \mathcal{N}_i}$,
and then averages the remaining values together with the agent's own value $\lambda_i$, which is never subject to trimming.

\textit{First coordinate:}
Received dual variables are $\{1,\,0,\,100\}$ (from $\check{\lambda}_{j_1},\,\check{\lambda}_{j_2},\,\check{\lambda}_{j_3}$).
After discarding the smallest ($0$) and largest ($100$), the remaining is $\{1\}$.
Combining with $\lambda_i=1$, the averaged value is $(1+1)/2 = 1$.

\textit{Second coordinate:}
Received dual variables are $\{0,\,1,\,100\}$ (from $\check{\lambda}_{j_1},\,\check{\lambda}_{j_2},\,\check{\lambda}_{j_3}$).
After discarding the smallest ($0$) and largest ($100$), the remaining is $\{1\}$.
Combining with $\lambda_i = 0$, the averaged value is $(0 + 1)/2 = 0.5$.

Then, the final output of $CTM(\cdot)$ is given in the form of
$CTM(\lambda_i,\{\check{\lambda}_j\}_{j \in \mathcal{N}_i}) = [1,0.5]^\top$.
Hence, we obtain \[\|CTM(\lambda_i,\{\check{\lambda}_j\}_{j \in \mathcal{N}_i})\|=\sqrt{1.25} > 1=\max_{j\in\mathcal{N}_i\cap\mathcal{H}\cup\{i\}}\|\lambda_j\| .\]

Therefore, we can prove that the robust aggregation rule $CTM(\cdot)$ does not satisfy Property \ref{d2}.

\noindent\textbf{Robust Aggregation Rule $IOS(\cdot)$.} We provide a counter-example to demonstrate that the  robust aggregation rule $IOS(\cdot)$ does not satisfy Property \ref{d2} even in the one-dimensional case.
Specifically, consider a benign agent $i$ whose local dual variable is given by \[\lambda_{i} = -9 . \]
Benign agent $i$ receives dual variables from three neighbors $\mathcal{N}_{i} = \{j_1, j_2, j_3\}$.
The received dual variables are
\[
\check{\lambda}_{j_1} = -5,\quad
\check{\lambda}_{j_2} = 10,\quad
\check{\lambda}_{j_3} = -20.
\]
Among them, $\check{\lambda}_{j_3}$ is from a Byzantine neighbor, while $\check{\lambda}_{j_1}$ and $\check{\lambda}_{j_2}$ are from benign neighbors.
All benign dual variables, including $\lambda_i$, have an absolute value of at most $10$, such that
\[
\max_{j\in\mathcal{N}_i\cap\mathcal{H}\cup\{i\}}\|\lambda_j\|=\|\check{\lambda}_{j_2}\| =10.
\]
Given four total inputs (including $\lambda_i$) and an upper bound of Byzantine neighbors $b_i = 1$, the operator $IOS(\cdot)$ computes a weighted average of all values, identifies the point with the largest distance from this average, discards it, and then aggregates the remaining values using their normalized weights.

Assuming the $i$-th row of weight matrix $\widetilde{e}_{i,:} = [\frac{1}{4}, \frac{1}{4}, \frac{1}{4}, \frac{1}{4}]$, agent $i$ first computes the weighted average of the received dual variables:
\[
\bar{\lambda} = \frac{(-9) + (-5) + 10 +(-20)}{4} = -6.
\]

Next, the distance from each received dual variable to $\bar{\lambda}$ is calculated as
\[
\begin{aligned}
\|\check{\lambda}_{j_1} - \bar{\lambda}\| &= \|-5 - (-6)\| = 1, \\
\|\check{\lambda}_{j_2} - \bar{\lambda}\| &= \|10 - (-6)\| = 16, \\
\|\check{\lambda}_{j_3} - \bar{\lambda}\| &= \|-20 - (-6)\| = 14.
\end{aligned}
\]

The dual variable $\check{\lambda}_{j_2} = 10$ has the largest deviation and is hence discarded.
The remaining dual variables are $\{-9,\ -5,\ -20\}$ with weights $[\frac{1}{4},\frac{1}{4},\frac{1}{4}]$, which are further normalized to $\left[\frac{1}{3},\ \frac{1}{3},\ \frac{1}{3}\right]$.

Then, the final output of $IOS(\cdot)$ is 
$IOS(\lambda_i, \{\check{\lambda}_j\}_{j \in \mathcal{N}_i}) = \frac{(-9) + (-5) +(-20)}{3}  = -11.\dot{3}.$ Hence, we obtain
\[
\|IOS(\lambda_i, \{\check{\lambda}_j\}_{j \in \mathcal{N}_i})\| = 11.\dot{3} > 10 = \max_{j\in\mathcal{N}_i\cap\mathcal{H} \cup \{i\}} \|\lambda_j\|.
\]

Therefore, we prove that the robust aggregation rule $IOS(\cdot)$ does not satisfy Property \ref{d2}.

\noindent\textbf{Robust Aggregation Rule $SCC(\cdot)$} We provide a counter-example to demonstrate that the robust aggregation rule $SCC(\cdot)$ does not satisfy Property \ref{d2} even in the one-dimensional case.
Specifically, consider a benign agent $i$ whose local dual variable is given by
\[
\lambda_i = -1.
\]
Benign Agent $i$ receives dual variables from three neighbors $\mathcal{N}_i = \{j_1, j_2, j_3\}$. The received dual variables are
\[
\check{\lambda}_{j_1} = 1,\quad
\check{\lambda}_{j_2} = 1,\quad
\check{\lambda}_{j_3} = 100.
\]
Among them, $\check{\lambda}_{j_3}$ is from a Byzantine neighbor, while $\check{\lambda}_{j_1}$ and $\check{\lambda}_{j_2}$ are from benign neighbors. All benign dual variables, including $\lambda_i$, have an absolute value of $1$ such that
\[
\max_{j \in \mathcal{N}_i \cap \mathcal{H} \cup \{i\}} \|\lambda_j\|
=
\|\lambda_i\|
=
\|\check{\lambda}_{j_1}\|
=
\|\check{\lambda}_{j_2}\|
=1.
\]

The $SCC(\cdot)$ operator aggregates dual variables using a clipping-based strategy. It first computes a clipping threshold $\tau_i$, and then clips each received dual variable before performing weighted averaging.

Assuming that the $i$-th row of the weight matrix is
$\widetilde{e}_{i,:} = [0.05,\ 0.425,\ 0.425,\ 0.1]$,
agent $i$ computes the clipping threshold $\tau_i$ in $SCC(\cdot)$ by following the ideal threshold selection rule proposed in \cite{b-LieHe-2022}, as
\begin{align*}
\tau_i
&=
\sqrt{
\frac{1}{\widetilde{e}_{ij_3}}
\sum_{j \in \mathcal{N}_i \cap \mathcal{H} \cup \{i\}}
\widetilde{e}_{ij}\|\lambda_i-\check{\lambda}_j\|^2
} \\
&=
\sqrt{
\frac{1}{0.1}
\left(
0.05\cdot 0 + 0.425\cdot 4 + 0.425\cdot 4
\right)
} =
\sqrt{34}
\approx 5.83.
\end{align*}

Next, benign Agent $i$ clips its received dual variables according to the clipping threshold $\tau_i \approx 5.83$.
\begin{align*}
\check{\lambda}_{j_1}^{\text{clipped}}
&=
\lambda_i + clip(\check{\lambda}_{j_1}-\lambda_i,\tau_i) \\
&=
\lambda_i+\min\left(1,\frac{\tau_i}{\|\check{\lambda}_{j_1}-\lambda_i\|}\right)\cdot(\check{\lambda}_{j_1}-\lambda_i) \\
&=
-1+\min\left(1,\frac{5.83}{2}\right)\cdot(1-(-1)) \\
&=1, 
\end{align*}
\begin{align*}
\check{\lambda}_{j_2}^{\text{clipped}}
&=
\lambda_i + clip(\check{\lambda}_{j_2}-\lambda_i,\tau_i) \\
&=
\lambda_i+\min\left(1,\frac{\tau_i}{\|\check{\lambda}_{j_2}-\lambda_i\|}\right)\cdot(\check{\lambda}_{j_2}-\lambda_i) \\
&=
-1+\min\left(1,\frac{5.83}{2}\right)\cdot(1-(-1)) \\
&=1, \\
\check{\lambda}_{j_3}^{\text{clipped}}
&=
\lambda_i + clip(\check{\lambda}_{j_3}-\lambda_i,\tau_i) \\
&=
\lambda_i+\min\left(1,\frac{\tau_i}{\|\check{\lambda}_{j_3}-\lambda_i\|}\right)\cdot(\check{\lambda}_{j_3}-\lambda_i) \\
&=
-1+\min\left(1,\frac{5.83}{101}\right)\cdot(100-(-1)) \\
&=
-1+\sqrt{34}
\approx 4.83.
\end{align*}

Then, the final output of $SCC(\cdot)$ is
$SCC(\lambda_i, \{\check{\lambda}_j\}_{j \in \mathcal{N}_i})
=
0.05 \cdot (-1) + 0.425 \cdot 1 + 0.425 \cdot 1 + 0.1 \cdot (-1+\sqrt{34})
=
0.7 + 0.1\sqrt{34}$.
Hence, we obtain
\[
\|SCC(\lambda_i, \{\check{\lambda}_j\}_{j \in \mathcal{N}_i})\|
\approx 1.283
>
1
=
\max_{j \in \mathcal{N}_i \cap \mathcal{H} \cup \{i\}} \|\lambda_j\|.
\]
Therefore, we prove that the robust aggregation rule $SCC(\cdot)$ does not satisfy Property \ref{d2}.

\subsection{Satisfaction of Property \ref{d1} by Robust Aggregation Rules $AGG(\cdot)$}\label{aa3}

\noindent\textbf{Robust Aggregation Rule $CTM(\cdot)$}
\begin{lemma}
\label{Definition1-CTM}
For any benign agent $i$, suppose in $CTM(\cdot)$ the number of discarding messages   $2b_{i}$ as $2|\mathcal{N}_{i}\cap \mathcal{B}|$. The associated  weight matrix $E$ is row stochastic and its each elements $e_{ij}$ is given by
$$e_{ij}=\frac{1}{|\mathcal{N}_{i}\cap \mathcal{H}\cup \{i\}|}=\frac{1}{|\mathcal{N}_{i}|-b_{i}+1}.$$
Then, the robust aggregation rule $CTM(\cdot)$ satisfies Property \ref{d1} with the contraction constant $$\rho \le \max_{i\in \mathcal{H} }[\frac{6b_{i}(|\mathcal{N}_{i}|-b_{i}+1)}{(|\mathcal{N}_{i}|-2b_{i}+1)^{2}}].$$
\begin{IEEEproof}
\label{proof-supp-definition1-CTM}
Given that $CTM(\cdot)$ is coordinate-wise, we consider dimension $d^{'}$ firstly. Next, we analyze $CTM(\cdot)$ from the following two cases. Denote the remaining agents after $CTM(\cdot)$ in dimension $d^{'}$ as $[U_{i}]_{d^{'}} \subset \mathcal{N}_{i} \cup \{i\} $.

\textit{Case 1: Benign agent $i$ removes all Byzantine messages.} Benign agent $i$ removes all Byzantine messages in dimension $d^{'}$ means $[U_{i}]_{d^{'}}\cap \mathcal{B}= \varnothing$ and $[U_{i}]_{d^{'}}\subset \mathcal{N}_{i}\cap \mathcal{H} \cup \{i\}$. Thus, we have
\begin{align}
\label{proof-supp-definition1-CTM-1}
&\|[CTM(\lambda_{i}, \{\check{\lambda}_{j}\}_{j\in \mathcal{N}_{i} })]_{d^{'}}-[\bar{\lambda}_{i}]_{d^{'}}\|^{2}\\
=& \|\frac{1}{|\mathcal{N}_{i}|-2b_{i}+1}\sum_{j\in [U_{i}]_{d^{'}}}[\lambda_{j}]_{d^{'}}-[\bar{\lambda}_{i}]_{d^{'}}\|^{2}\notag\\
=&\|\frac{1}{|\mathcal{N}_{i}|-2b_{i}+1}\sum_{j\in [U_{i}]_{d^{'}}}([\lambda_{j}]_{d^{'}}-[\bar{\lambda}_{i}]_{d^{'}})\|^{2} \notag\\
=&\|\frac{1}{|\mathcal{N}_{i}|-2b_{i}+1}\sum_{j\in [U_{i}]_{d^{'}}}([\lambda_{j}]_{d^{'}}-[\bar{\lambda}_{i}]_{d^{'}})\notag
\end{align}
\begin{align}
&-\frac{1}{|\mathcal{N}_{i}|-2b_{i}+1}\sum_{j\in \mathcal{N}_{i}\cap \mathcal{H}\cup \{i\}}([\lambda_{j}]_{d^{'}}-[\bar{\lambda}_{i}]_{d^{'}})\|^{2}\notag\\
=& \|-\frac{1}{|\mathcal{N}_{i}|-2b_{i}+1}\sum_{j\in \mathcal{N}_{i}\cap \mathcal{H}\cup \{i\}\setminus [U_{i}]_{d^{'}}}([\lambda_{j}]_{d^{'}}-[\bar{\lambda}_{i}]_{d^{'}})\|^{2}\notag\\
=&\frac{1}{(|\mathcal{N}_{i}|-2b_{i}+1)^{2}}\|\sum_{j\in \mathcal{N}_{i}\cap \mathcal{H}\cup \{i\}\setminus [U_{i}]_{d^{'}}}([\lambda_{j}]_{d^{'}}-[\bar{\lambda}_{i}]_{d^{'}})\|^{2} \notag\\
\le& \frac{|\mathcal{N}_{i}\cap \mathcal{H}\cup \{i\}\setminus [U_{i}]_{d^{'}}|}{(|\mathcal{N}_{i}|-2b_{i}+1)^{2}}\sum_{j\in \mathcal{N}_{i}\cap \mathcal{H}\cup \{i\}\setminus [U_{i}]_{d^{'}}}\|[\lambda_{j}]_{d^{'}}-[\bar{\lambda}_{i}]_{d^{'}}\|^{2} \notag,
\end{align}
in which the third equality holds true because of the fact $\sum_{j\in \mathcal{N}_{i}\cap \mathcal{H}\cup \{i\}}([\lambda_{j}]_{d^{'}}-[\bar{\lambda}_{i}]_{d^{'}})=0$. To derive the last inequality, we use the mean inequality. Combining the facts $|\mathcal{N}_{i}\cap \mathcal{H}\cup \{i\}\setminus [U_{i}]_{d^{'}}|=b_{i}$ and $e_{ij}=\frac{1}{|\mathcal{N}_{i}|-b_{i}+1}$, we get
\begin{align}
\label{proof-supp-definition1-CTM-2}
&\|[CTM(\lambda_{i}, \{\lambda_{j}\}_{j\in \mathcal{N}_{i} })]_{d^{'}}-[\bar{\lambda}_{i}]_{d^{'}}\|^{2}\\
\le &  \frac{b_{i}}{(|\mathcal{N}_{i}|-2b_{i}+1)^{2}}\sum_{j\in \mathcal{N}_{i}\cap \mathcal{H}\cup \{i\}\setminus [U_{i}]_{d^{'}}}\|[\lambda_{j}]_{d^{'}}-[\bar{\lambda}_{i}]_{d^{'}}\|^{2} \notag\\
\le & \frac{b_{i}(|\mathcal{N}_{i}|-b_{i}+1)}{(|\mathcal{N}_{i}|-2b_{i}+1)^{2}}\cdot \frac{1}{|\mathcal{N}_{i}|-b_{i}+1}\sum_{j\in \mathcal{N}_{i}\cap \mathcal{H}\cup \{i\}}\|[\lambda_{j}]_{d^{'}}-[\bar{\lambda}_{i}]_{d^{'}}\|^{2} \notag\\
=& \frac{b_{i}(|\mathcal{N}_{i}|-b_{i}+1)}{(|\mathcal{N}_{i}|-2b_{i}+1)^{2}}\sum_{j\in \mathcal{N}_{i}\cap \mathcal{H}\cup \{i\}}e_{ij}\|[\lambda_{j}]_{d^{'}}-[\bar{\lambda}_{i}]_{d^{'}}\|^{2}  \notag.
\end{align}

\textit{Case 2: Benign agent $i$ cannot remove all Byzantine messages.}
Benign agent $i$ cannot remove all Byzantine messages means $[U_{i}]_{d^{'}}\cap \mathcal{B}\ne \varnothing$. Thus, we have
\begin{align}
\label{proof-supp-definition1-CTM-3}
&\|[CTM(\lambda_{i}, \{\check{\lambda}_{j}\}_{j\in \mathcal{N}_{i} })]_{d^{'}}-[\bar{\lambda}_{i}]_{d^{'}}\|^{2}\\
=& \|\frac{1}{|\mathcal{N}_{i}|-2b_{i}+1}\sum_{j\in [U_{i}]_{d^{'}}}[\check{\lambda}_{j}]_{d^{'}}-[\bar{\lambda}_{i}]_{d^{'}}\|^{2}\notag\\
=&\|\frac{1}{|\mathcal{N}_{i}|-2b_{i}+1}\sum_{j\in [U_{i}]_{d^{'}}}([\check{\lambda}_{j}]_{d^{'}}-[\bar{\lambda}_{i}]_{d^{'}})\|^{2} \notag\\
=&\|\frac{1}{|\mathcal{N}_{i}|-2b_{i}+1}\sum_{j\in [U_{i}]_{d^{'}}}([\check{\lambda}_{j}]_{d^{'}}-[\bar{\lambda}_{i}]_{d^{'}})\notag\\
&-\frac{1}{|\mathcal{N}_{i}|-2b_{i}+1}\sum_{j\in \mathcal{N}_{i}\cap \mathcal{H}\cup \{i\}}([\lambda_{j}]_{d^{'}}-[\bar{\lambda}_{i}]_{d^{'}})\|^{2}\notag\\
=& \|\frac{1}{|\mathcal{N}_{i}|-2b_{i}+1}\sum_{j\in [U_{i}]_{d^{'}}\setminus (\mathcal{N}_{i}\cap \mathcal{H}\cup \{i\})}([\check{\lambda}_{j}]_{d^{'}}-[\bar{\lambda}_{i}]_{d^{'}})\notag\\
&-\frac{1}{|\mathcal{N}_{i}|-2b_{i}+1}\sum_{j\in \mathcal{N}_{i}\cap \mathcal{H}\cup \{i\}\setminus [U_{i}]_{d^{'}}}([\lambda_{j}]_{d^{'}}-[\bar{\lambda}_{i}]_{d^{'}})\|^{2}\notag\\
\le &  \frac{ 2}{(|\mathcal{N}_{i}|-2b_{i}+1)^{2}}[\|\sum_{j\in [U_{i}]_{d^{'}}\setminus (\mathcal{N}_{i}\cap \mathcal{H}\cup \{i\})}([\check{\lambda}_{j}]_{d^{'}}-[\bar{\lambda}_{i}]_{d^{'}})\|^{2}\notag\\
&+\|\sum_{j\in \mathcal{N}_{i}\cap \mathcal{H}\cup \{i\}\setminus [U_{i}]_{d^{'}}}([\lambda_{j}]_{d^{'}}-[\bar{\lambda}_{i}]_{d^{'}})\|^{2}] \notag\\
\le & \frac{2|[U_{i}]_{d^{'}}\setminus (\mathcal{N}_{i}\cap \mathcal{H}\cup \{i\})|}{(|\mathcal{N}_{i}|-2b_{i}+1)^{2}}\sum_{j\in [U_{i}]_{d^{'}}\setminus (\mathcal{N}_{i}\cap \mathcal{H}\cup \{i\})}\|[\check{\lambda}_{j}]_{d^{'}}-[\bar{\lambda}_{i}]_{d^{'}}\|^{2}\notag\\
&+ \frac{2| \mathcal{N}_{i}\cap \mathcal{H}\cup \{i\}\setminus [U_{i}]_{d^{'}}|}{(|\mathcal{N}_{i}|-2b_{i}+1)^{2}} \sum_{j\in \mathcal{N}_{i}\cap \mathcal{H}\cup \{i\}\setminus [U_{i}]_{d^{'}}}\|[\lambda_{j}]_{d^{'}}-[\bar{\lambda}_{i}]_{d^{'}}\|^{2} \notag,
\end{align}
in which the third equality holds true because of the fact $\sum_{j\in \mathcal{N}_{i}\cap \mathcal{H}\cup \{i\}}([\lambda_{j}]_{d^{'}}-[\bar{\lambda}_{i}]_{d^{'}})=0$. To derive the last inequality, we use the mean inequality. From the scheme of $CTM(\cdot)$, we can obtain $\sum_{j\in [U_{i}]_{d^{'}}\setminus (\mathcal{N}_{i}\cap \mathcal{H}\cup \{i\})}\|[\check{\lambda}_{j}]_{d^{'}}-[\bar{\lambda}_{i}]_{d^{'}}\|^{2}\le \sum_{j\in \mathcal{N}_{i}\cap \mathcal{H}\cup \{i\}\setminus [U_{i}]_{d^{'}}}\|[\lambda_{j}]_{d^{'}}-[\bar{\lambda}_{i}]_{d^{'}}\|^{2} $. Therefore, we have
\begin{align}
\label{proof-supp-definition1-CTM-4}
&\|[CTM(\lambda_{i}, \{\lambda_{j}\}_{j\in \mathcal{N}_{i} })]_{d^{'}}-[\bar{\lambda}_{i}]_{d^{'}}\|^{2}\\
\le &\frac{2| \mathcal{N}_{i}\cap \mathcal{H}\cup \{i\}\setminus [U_{i}]_{d^{'}}|+2|[U_{i}]_{d^{'}}\setminus (\mathcal{N}_{i}\cap \mathcal{H}\cup \{i\})|}{(|\mathcal{N}_{i}|-2b_{i}+1)^{2}} \cdot \notag\\
&\sum_{j\in \mathcal{N}_{i}\cap \mathcal{H}\cup \{i\}\setminus [U_{i}]_{d^{'}}}\|[\lambda_{j}]_{d^{'}}-[\bar{\lambda}_{i}]_{d^{'}}\|^{2} \notag\\
\le & \frac{2b_{i}+4b_{i}}{(|\mathcal{N}_{i}|-2b_{i}+1)^{2}} \sum_{j\in \mathcal{N}_{i}\cap \mathcal{H}\cup \{i\}\setminus [U_{i}]_{d^{'}}}\|[\lambda_{j}]_{d^{'}}-[\bar{\lambda}_{i}]_{d^{'}}\|^{2}\notag\\
\le & \frac{6b_{i}(|\mathcal{N}_{i}|-b_{i}+1)}{(|\mathcal{N}_{i}|-2b_{i}+1)^{2}}\cdot \frac{1}{|\mathcal{N}_{i}|-b_{i}+1} \sum_{j\in \mathcal{N}_{i}\cap \mathcal{H}\cup \{i\}}\|[\lambda_{j}]_{d^{'}}-[\bar{\lambda}_{i}]_{d^{'}}\|^{2} \notag\\
=& \frac{6b_{i}(|\mathcal{N}_{i}|-b_{i}+1)}{(|\mathcal{N}_{i}|-2b_{i}+1)^{2}}\sum_{j\in \mathcal{N}_{i}\cap \mathcal{H}\cup \{i\}}e_{ij}\|[\lambda_{j}]_{d^{'}}-[\bar{\lambda}_{i}]_{d^{'}}\|^{2}\notag,
\end{align}
where the second inequality holds since $| \mathcal{N}_{i}\cap \mathcal{H}\cup \{i\}\setminus [U_{i}]_{d^{'}}|=|\mathcal{N}_{i}\cap \mathcal{H}\cup \{i\}\cup [U_{i}]_{d^{'}}|-|[U_{i}]_{d^{'}}|\le |\mathcal{N}_{i}|+1-(|\mathcal{N}_{i}|-2b_{i}+1)=2b_{i} $ and $| [U_{i}]_{d^{'}} \setminus \mathcal{N}_{i}\cap \mathcal{H}\cup \{i\}|=|\mathcal{N}_{i}\cap \mathcal{H}\cup \{i\}\cup [U_{i}]_{d^{'}}|-|\mathcal{N}_{i}\cap \mathcal{H}\cup \{i\}|\le |\mathcal{N}_{i}|+1-(|\mathcal{N}_{i}|-b_{i}+1)=b_{i} $.

Combining \eqref{proof-supp-definition1-CTM-2}
and \eqref{proof-supp-definition1-CTM-4}, we have
\begin{align}
\label{proof-supp-definition1-CTM-5}
&\|[CTM(\lambda_{i}, \{\check{\lambda}_{j}\}_{j\in \mathcal{N}_{i} })]_{d^{'}}-[\bar{\lambda}_{i}]_{d^{'}}\|^{2}\\
\le& \frac{6b_{i}(|\mathcal{N}_{i}|-b_{i}+1)}{(|\mathcal{N}_{i}|-2b_{i}+1)^{2}}\sum_{j\in \mathcal{N}_{i}\cap \mathcal{H}\cup \{i\}}e_{ij}\|[\lambda_{j}]_{d^{'}}-[\bar{\lambda}_{i}]_{d^{'}}\|^{2}\notag,
\end{align}
Extending \eqref{proof-supp-definition1-CTM-5} into high dimension, we have
\begin{align}
\label{proof-supp-definition1-CTM-6}
&\|CTM(\lambda_{i}, \{\check{\lambda}_{j}\}_{j\in \mathcal{N}_{i} })-\bar{\lambda}_{i}\|^{2}\\
=& \sum_{d^{'}=1}^{d}\|[CTM(\lambda_{i}, \{\check{\lambda}_{j}\}_{j\in \mathcal{N}_{i} })]_{d^{'}}-[\bar{\lambda}_{i}]_{d^{'}}\|^{2}\notag\\
\le& \frac{6b_{i}(|\mathcal{N}_{i}|-b_{i}+1)}{(|\mathcal{N}_{i}|-2b_{i}+1)^{2}}\sum_{d^{'}=1}^{d} \sum_{j\in \mathcal{N}_{i}\cap \mathcal{H}\cup \{i\}}e_{ij}\|[\lambda_{j}]_{d^{'}}-[\bar{\lambda}_{i}]_{d^{'}}\|^{2}\notag\\
\le &\frac{6b_{i}(|\mathcal{N}_{i}|-b_{i}+1)}{(|\mathcal{N}_{i}|-2b_{i}+1)^{2}} \sum_{j\in \mathcal{N}_{i}\cap \mathcal{H}\cup \{i\}}e_{ij}\|\lambda_{j}-\bar{\lambda}_{i}]\|^{2} \notag,
\end{align}
which shows that for any benign agent $i$ the robust aggregation rule $CTM(\cdot)$ satisfies Property \ref{d1} with the contraction constant $\rho=\frac{6b_{i}(|\mathcal{N}_{i}|-b_{i}+1)}{(|\mathcal{N}_{i}|-2b_{i}+1)^{2}} $. 
\end{IEEEproof}
\end{lemma}

\noindent\textbf{Robust Aggregation Rule $SCC(ARC(\cdot))$}
\begin{lemma}
\label{Definition1-SCC}
For any benign agent $i$, suppose the clipping parameter  $\tau_{i}=\sqrt{\frac{1}{\sum_{j\in \mathcal{N}_{i}\cap \mathcal{B} }\tilde{e}_{ij}}\sum_{j\in \mathcal{N}_{i}\cap \mathcal{H}\cup\{i\}}\tilde{e}_{ij}\|\lambda_{i}-\lambda_{j}\|^{2}}$  in $SCC(\cdot)$  and the number of clipping messages  $b_{i}=|\mathcal{N}_{i}\cap \mathcal{B}|$  in $ARC(\cdot)$. The associated  weight matrix $E$ is doubly stochastic and its each elements $e_{ij}$ is given by
$$e_{ij}=\left\{\begin{matrix}
\tilde{e}_{ij}+\sum_{j^{'}\in \mathcal{N}_{i} \cap \mathcal{B}}\tilde{e}_{ij^{'}}, &\text{ if } i=j,
 \\
\tilde{e}_{ij},  &\text{ if } i\ne j.
\end{matrix}\right.$$
Then the robust aggregation rule $SCC(ARC(\cdot))$ satisfies Property \ref{d1} with the contraction constant
\begin{align*}
\rho\le  \max_{i\in \mathcal{H}}[\frac{8 \sum_{j\in \mathcal{N}_{i}\cap \mathcal{B}}\tilde{e}_{ij}(1+\min_{j\in \mathcal{N}_{i}\cap \mathcal{H}\cup \{i\}}\tilde{e}_{ij})}{\min_{j\in \mathcal{N}_{i}\cap \mathcal{H}\cup \{i\}}\tilde{e}_{ij}}\\
+\frac{|\mathcal{N}_{i}\cap \mathcal{B}|\cdot \max_{j\in \mathcal{N}_{i}\cap \mathcal{H}\cup \{i\}} \tilde{e}_{ij}}{1-|\mathcal{N}_{i}\cap \mathcal{B}|\cdot \max_{j\in \mathcal{N}_{i}\cap \mathcal{H}\cup \{i\}} \tilde{e}_{ij}}].
\end{align*}
\end{lemma}
\begin{IEEEproof}
\label{proof-supp-SCC}
First, we prove that the robust aggregation rule $SCC(\cdot)$ satisfies Property \ref{d1} and analyze the corresponding contraction constant $\tilde{\rho}$. Then, combining the conclusion of the robust aggregation rule $SCC(\cdot)$ and Lemma \ref{lemma-AGG-ARC-definition1}, we prove that the robust aggregation rule $AGG(\cdot)=SCC(ARC(\cdot))$ satisfies Property \ref{d1} and analyze the corresponding contraction constant $\rho$.

According to $(68)$ in \cite{b-Haoxiang-Ye-2023}, we have
\begin{align}
\label{proof-supp-definition1-SCC-1}
&\|SCC(\lambda_{i}, \{\check{\lambda}_{j}\}_{j\in \mathcal{N}_{i}})-\bar{\lambda}_{i}\|^{2}\\
\le & 4 \sum_{j\in \mathcal{N}_{i}\cap \mathcal{B}}\tilde{e}_{ij}\sum_{j\in \mathcal{N}_{i}\cap \mathcal{H}}\tilde{e}_{ij}\|\lambda_{i}-\lambda_{j}\|^{2} \notag \\
\le & 4 \sum_{j\in \mathcal{N}_{i}\cap \mathcal{B}}\tilde{e}_{ij}\sum_{j\in \mathcal{N}_{i}\cap \mathcal{H}\cup \{i\}}e_{ij}\|\lambda_{i}-\lambda_{j}\|^{2}  \notag \\
\le & 4 \sum_{j\in \mathcal{N}_{i}\cap \mathcal{B}}\tilde{e}_{ij}[2\|\lambda_{i}-\bar{\lambda}_{i}\|^{2}+2\sum_{j\in \mathcal{N}_{i}\cap \mathcal{H}\cup \{i\}}e_{ij}\|\lambda_{j}-\bar{\lambda}_{i}\|^{2} ] \notag
\end{align}
\begin{align}
\le &  4 \sum_{j\in \mathcal{N}_{i}\cap \mathcal{B}}\tilde{e}_{ij}[2\max_{j\in \mathcal{N}_{i}\cap \mathcal{H}\cup \{i\}}\|\lambda_{j}-\bar{\lambda}_{i}\|^{2}\notag\\
&+2\sum_{j\in \mathcal{N}_{i}\cap \mathcal{H}\cup \{i\}}e_{ij}\|\lambda_{j}-\bar{\lambda}_{i}\|^{2} ] \notag \\
\le & 4 \sum_{j\in \mathcal{N}_{i}\cap \mathcal{B}}\tilde{e}_{ij}[\frac{2}{\min_{j\in \mathcal{N}_{i}\cap \mathcal{H}\cup \{i\}}\tilde{e}_{ij}}\sum_{j\in \mathcal{N}_{i}\cap \mathcal{H}\cup \{i\}}e_{ij}\|\lambda_{j}-\bar{\lambda}_{i}\|^{2}\notag\\
&+2\sum_{j\in \mathcal{N}_{i}\cap \mathcal{H}\cup \{i\}}e_{ij}\|\lambda_{j}-\bar{\lambda}_{i}\|^{2} ] \notag\\
=& \frac{8 \sum_{j\in \mathcal{N}_{i}\cap \mathcal{B}}\tilde{e}_{ij}(1+\min_{j\in \mathcal{N}_{i}\cap \mathcal{H}\cup \{i\}}\tilde{e}_{ij})}{\min_{j\in \mathcal{N}_{i}\cap \mathcal{H}\cup \{i\}}\tilde{e}_{ij}}\cdot  \notag\\
&  \sum_{j\in \mathcal{N}_{i}\cap \mathcal{H}\cup \{i\}}e_{ij}\|\lambda_{j}-\bar{\lambda}_{i}\|^{2} \notag,
\end{align}
which shows that for any benign agent $i$, the robust aggregation rule $SCC(\cdot)$ satisfies Property \ref{d1}. Specifically, we observe that the contraction constant $\tilde{\rho}= \frac{8 \sum_{j\in \mathcal{N}_{i}\cap \mathcal{B}}\tilde{e}_{ij}(1+\min_{j\in \mathcal{N}_{i}\cap \mathcal{H}\cup \{i\}}\tilde{e}_{ij})}{\min_{j\in \mathcal{N}_{i}\cap \mathcal{H}\cup \{i\}}\tilde{e}_{ij}}$. Therefore, combining \eqref{proof-supp-definition1-SCC-1} and Lemma \ref{lemma-AGG-ARC-definition1}, we have
\begin{align}
\label{proof-supp-definition1-SCC-2}
&\|SCC(ARC(\lambda_{i}, \{\check{\lambda}_{j}\}_{j\in \mathcal{N}_{i}}))-\bar{\lambda}_{i}\|^{2}\\
\le& [\frac{8 \sum_{j\in \mathcal{N}_{i}\cap \mathcal{B}}\tilde{e}_{ij}(1+\min_{j\in \mathcal{N}_{i}\cap \mathcal{H}\cup \{i\}}\tilde{e}_{ij})}{\min_{j\in \mathcal{N}_{i}\cap \mathcal{H}\cup \{i\}}\tilde{e}_{ij}} +  \notag\\
& \frac{|\mathcal{N}_{i}\cap \mathcal{B}|\cdot \max_{j\in \mathcal{N}_{i}\cap \mathcal{H}\cup \{i\}} \tilde{e}_{ij}}{1-|\mathcal{N}_{i}\cap \mathcal{B}|\cdot \max_{j\in \mathcal{N}_{i}\cap \mathcal{H}\cup \{i\}} \tilde{e}_{ij}}]\cdot\sum_{j\in \mathcal{N}_{i}\cap \mathcal{H}\cup \{i\}}e_{ij}\|\lambda_{j}-\bar{\lambda}_{i}\|^{2} \notag .
\end{align}
\end{IEEEproof}

\noindent\textbf{Robust Aggregation Rule $IOS(ARC(\cdot))$}
\begin{lemma}
\label{Definition1-IOS}
For any benign agent $i$, suppose in $IOS(\cdot)$ the number of discarding messages  $b_{i}$ as $| \mathcal{N}_{i} \cap \mathcal{B}|$ and in $ARC(\cdot)$ the number of clipping messages  $b_{i}$ as $|\mathcal{N}_{i}\cap \mathcal{B}|$.
Define a neighbor set that includes the neighbors with the largest $b_{i}$ weights, as $\mathcal{N}_{i,b_{i}}:=\underset{\mathcal{N}^{'}\subseteq \mathcal{N}_{i},|\mathcal{N}^{'}|=b_{i}} {\arg\max} \sum_{j\in \mathcal{N}^{'}}\tilde{e}_{ij}.$
When $\sum_{j\in \mathcal{N}_{i,b_{i}}}\tilde{e}_{ij}< \frac{1}{3}$, the associated  weight matrix $E$ is doubly stochastic and its each elements $e_{ij}$ is given by
$$e_{ij}=\left\{\begin{matrix}
\tilde{e}_{ij}+\sum_{j^{'}\in \mathcal{N}_{i} \cap \mathcal{B}}\tilde{e}_{ij^{'}}, &\text{ if } i=j,
 \\
\tilde{e}_{ij},  &\text{ if } i\ne j.
\end{matrix}\right.$$
Then the robust aggregation rule $IOS(ARC(\cdot))$ satisfies Property \ref{d1} with the contraction constant
\begin{align*}
\rho\le \max_{i\in \mathcal{H}}[\frac{(15\sum_{j\in \mathcal{N}_{i,b_{i}}}\tilde{e}_{ij})^{2}}{\min_{j\in \mathcal{N}_{i}\cap \mathcal{H}\cup \{i\}}\tilde{e}_{ij}^{2}(1-3\sum_{j\in \mathcal{N}_{i,b_{i}}}\tilde{e}_{ij})^{2}}\\
+\frac{|\mathcal{N}_{i}\cap \mathcal{B}|\cdot \max_{j\in \mathcal{N}_{i}\cap \mathcal{H}\cup \{i\}} \tilde{e}_{ij}}{1-|\mathcal{N}_{i}\cap \mathcal{B}|\cdot \max_{j\in \mathcal{N}_{i}\cap \mathcal{H}\cup \{i\}} \tilde{e}_{ij}}].
\end{align*}
\end{lemma}
\begin{IEEEproof}
\label{proof-supp-IOS}
First, we prove that the robust aggregation rule $IOS(\cdot)$ satisfies Property \ref{d1} and analyze the corresponding contraction constant $\tilde{ \rho}$. Then, combining the conclusion of the robust aggregation rule $IOS(\cdot)$ and Lemma \ref{lemma-AGG-ARC-definition1}, we prove that the robust aggregation rule $AGG(\cdot)=IOS(ARC(\cdot))$ satisfies Property \ref{d1} and analyze the according contraction constant $\rho$.

Denote the remaining agents after $IOS(\cdot)$ as $U_{i} \subseteq \mathcal{N}_{i} \cup \{i\} $. We analyze $IOS(\cdot)$ from the following two cases.

\textit{Case 1: Benign agent $i$ removes all Byzantine messages.}
According to $(75)$ in \cite{b-Haoxiang-Ye-2023}, we obtain
\begin{align}
\label{proof-Definition1-IOS-case1-2}
&\|IOS(\lambda_{i},\{\check{\lambda}_{j}\}_{j\in \mathcal{N}_{i}})-\bar{\lambda}_{i}\| \\
\le & \frac{\sum_{j\in \mathcal{N}_{i} \cap\mathcal{B}}\tilde{e}_{ij}}{\sum_{j\in \mathcal{N}_{i} \cap\mathcal{H} \cup \{i\}}\tilde{e}_{ij}} \|\lambda_{i}-\bar{\lambda}_{i}\| \notag\\
\le &  \frac{\sum_{j\in \mathcal{N}_{i} \cap\mathcal{B}}\tilde{e}_{ij}}{1-\sum_{j\in \mathcal{N}_{i} \cap\mathcal{B}}\tilde{e}_{ij}} \max_{j \in  \mathcal{N}_{i} \cap\mathcal{H} \cup \{i\}}
\|\lambda_{j}-\bar{\lambda}_{i}\| \notag \\
\le & \frac{\sum_{j\in \mathcal{N}_{i,b_{i}}}\tilde{e}_{ij}}{1-\sum_{j\in \mathcal{N}_{i,b_{i}} }\tilde{e}_{ij}} \max_{j \in  \mathcal{N}_{i} \cap\mathcal{H} \cup \{i\}}
\|\lambda_{j}-\bar{\lambda}_{i}\| \notag\\
\le&  \frac{\sum_{j\in \mathcal{N}_{i,b_{i}}}\tilde{e}_{ij}}{\min_{j\in \mathcal{N}_{i}\cap \mathcal{H}\cup \{i\}}\tilde{e}_{ij}(1-\sum_{j\in \mathcal{N}_{i,b_{i}}}\tilde{e}_{ij})} \sum_{j\in \mathcal{N}_{i}\cap \mathcal{H}\cup \{i\}}e_{ij}\|\lambda_{j}-\bar{\lambda}_{i}\| \notag.
\end{align}
Taking squares for both sides of \eqref{proof-Definition1-IOS-case1-2} yields
\begin{align}
\label{proof-Definition1-IOS-case1-3}
&\|IOS(\lambda_{i},\{\check{\lambda}_{j}\}_{j\in \mathcal{N}_{i}})-\bar{\lambda}_{i}\|^{2} \\
\le&  \frac{(\sum_{j\in \mathcal{N}_{i,b_{i}}}\tilde{e}_{ij})^{2}}{\min_{j\in \mathcal{N}_{i}\cap \mathcal{H}\cup \{i\}}\tilde{e}_{ij}^{2}(1-\sum_{j\in \mathcal{N}_{i,b_{i}}}\tilde{e}_{ij})^{2}}\cdot \notag \\
&  (\sum_{j\in \mathcal{N}_{i}\cap \mathcal{H}\cup \{i\}}e_{ij}\|\lambda_{j}-\bar{\lambda}_{i}\|)^{2} \notag\\
\le & \frac{(\sum_{j\in \mathcal{N}_{i,b_{i}}}\tilde{e}_{ij})^{2}}{\min_{j\in \mathcal{N}_{i}\cap \mathcal{H}\cup \{i\}}\tilde{e}_{ij}^{2}(1-\sum_{j\in \mathcal{N}_{i,b_{i}}}\tilde{e}_{ij})^{2}}\cdot \notag \\
& \sum_{j\in \mathcal{N}_{i}\cap \mathcal{H}\cup \{i\}}e_{ij}\|\lambda_{j}-\bar{\lambda}_{i}\|^{2} \notag,
\end{align}
where the last inequality holds from Jensen's inequality and the row stochasticity of weight matrix $E$.

\textit{Case 2: Benign agent $i$ cannot remove all Byzantine messages.}
According to (17) in \cite{b-Haoxiang-Ye-2023}, we obtain
\begin{align}
\label{proof-Definition1-IOS-case2-1}
&\|IOS(\lambda_{i},\{\check{\lambda}_{j}\}_{j\in \mathcal{N}_{i}})-\bar{\lambda}_{i}\| \\
\le & \frac{15\sum_{j\in \mathcal{N}_{i,b_{i}}}\tilde{e}_{ij}}{1-3\sum_{j\in \mathcal{N}_{i,b_{i}} }\tilde{e}_{ij}} \max_{j \in  \mathcal{N}_{i} \cap\mathcal{H} \cup \{i\}}
\|\lambda_{j}-\bar{\lambda}_{i}\| \notag\\
\le&  \frac{15\sum_{j\in \mathcal{N}_{i,b_{i}}}\tilde{e}_{ij}}{\min_{j\in \mathcal{N}_{i}\cap \mathcal{H}\cup \{i\}}\tilde{e}_{ij}(1-3\sum_{j\in \mathcal{N}_{i,b_{i}}}\tilde{e}_{ij})}\cdot \notag\\
& \sum_{j\in \mathcal{N}_{i}\cap \mathcal{H}\cup \{i\}}e_{ij}\|\lambda_{j}-\bar{\lambda}_{i}\| \notag.
\end{align}
Taking squares for both sides of \eqref{proof-Definition1-IOS-case2-1} yields
\begin{align}
\label{proof-Definition1-IOS-case2-2}
&\|IOS(\lambda_{i},\{\check{\lambda}_{j}\}_{j\in \mathcal{N}_{i}})-\bar{\lambda}_{i}\|^{2} \\
\le & \frac{(15\sum_{j\in \mathcal{N}_{i,b_{i}}}\tilde{e}_{ij})^{2}}{(1-3\sum_{j\in \mathcal{N}_{i,b_{i}} }\tilde{e}_{ij})^{2}} (\sum_{j\in \mathcal{N}_{i}\cap \mathcal{H}\cup \{i\}}e_{ij}
\|\lambda_{j}-\bar{\lambda}_{i}\|)^{2} \notag\\
\le&  \frac{(15\sum_{j\in \mathcal{N}_{i,b_{i}}}\tilde{e}_{ij})^{2}}{\min_{j\in \mathcal{N}_{i}\cap \mathcal{H}\cup \{i\}}\tilde{e}_{ij}^{2}(1-3\sum_{j\in \mathcal{N}_{i,b_{i}}}\tilde{e}_{ij})^{2}}\cdot \notag\\
& \sum_{j\in \mathcal{N}_{i}\cap \mathcal{H}\cup \{i\}}e_{ij}\|\lambda_{j}-\bar{\lambda}_{i}\|^{2}, \notag
\end{align}
where the last inequality holds from Jensen's inequality and the row stochasticity of weight matrix $E$.

Combining \eqref{proof-Definition1-IOS-case1-3} and \eqref{proof-Definition1-IOS-case2-2}, we have
\begin{align}
\label{proof-Definition1-IOS-case2-3}
&\|IOS(\lambda_{i},\{\check{\lambda}_{j}\}_{j\in \mathcal{N}_{i}})-\bar{\lambda}_{i}\|^{2} \\
\le&  \frac{(15\sum_{j\in \mathcal{N}_{i,b_{i}}}\tilde{e}_{ij})^{2}}{\min_{j\in \mathcal{N}_{i}\cap \mathcal{H}\cup \{i\}}\tilde{e}_{ij}^{2}(1-3\sum_{j\in \mathcal{N}_{i,b_{i}}}\tilde{e}_{ij})^{2}}\cdot \notag\\
& \sum_{j\in \mathcal{N}_{i}\cap \mathcal{H}\cup \{i\}}e_{ij}\|\lambda_{j}-\bar{\lambda}_{i}\|^{2}, \notag
\end{align}
which shows that for any benign agent $i$ the robust aggregation rule $IOS(\cdot)$ satisfies Property \ref{d1}. Specifically, we observe that the contraction constant $\tilde{\rho}=\frac{(15\sum_{j\in \mathcal{N}_{i,b_{i}}}\tilde{e}_{ij})^{2}}{\min_{j\in \mathcal{N}_{i}\cap \mathcal{H}\cup \{i\}}\tilde{e}_{ij}^{2}(1-3\sum_{j\in \mathcal{N}_{i,b_{i}}}\tilde{e}_{ij})^{2}}$. Therefore, combining \eqref{proof-Definition1-IOS-case2-3} and Lemma \ref{lemma-AGG-ARC-definition1}, we have
\begin{align}
\label{proof-Definition1-IOS-case2-4}
&\|IOS(ARC(\lambda_{i},\{\check{\lambda}_{j}\}_{j\in \mathcal{N}_{i}}))-\bar{\lambda}_{i}\|^{2} \\
\le&  [\frac{(15\sum_{j\in \mathcal{N}_{i,b_{i}}}\tilde{e}_{ij})^{2}}{\min_{j\in \mathcal{N}_{i}\cap \mathcal{H}\cup \{i\}}\tilde{e}_{ij}^{2}(1-3\sum_{j\in \mathcal{N}_{i,b_{i}}}\tilde{e}_{ij})^{2}}+\notag\\
& \frac{|\mathcal{N}_{i}\cap \mathcal{B}|\cdot \max_{j\in \mathcal{N}_{i}\cap \mathcal{H}\cup \{i\}} \tilde{e}_{ij}}{1-|\mathcal{N}_{i}\cap \mathcal{B}|\cdot \max_{j\in \mathcal{N}_{i}\cap \mathcal{H}\cup \{i\}} \tilde{e}_{ij}}] \cdot \sum_{j\in \mathcal{N}_{i}\cap \mathcal{H}\cup \{i\}}e_{ij}\|\lambda_{j}-\bar{\lambda}_{i}\|^{2}. \notag
\end{align}
\end{IEEEproof}

\noindent\textbf{Supporting Lemmas for Satisfaction of Property \ref{d1}}
\begin{lemma}
\label{lemma-AGG-ARC-definition1}
Consider the robust aggregation rule $AGG(\cdot) := IOS(ARC(\cdot))$, or $SCC(ARC(\cdot))$. If the corresponding base robust aggregation rule namely,  $IOS(\cdot)$, or $SCC(\cdot)$ satisfies Property \ref{d1} with contraction constant $\tilde{\rho}$, then $AGG(\cdot)$ also satisfies Property \ref{d1} with contraction constant $\rho=\tilde{\rho}+\frac{\sum_{j\in S_{i}^{c}} e_{ij}}{1-\sum_{j\in S_{i}^{c}} e_{ij}}$, where $S_{i}^{c}:=\{j\in \mathcal{N}_{i}\cap \mathcal{H} \cup \{i\}, \|\lambda_{j}\| \ge C_{i} \}$.
\begin{IEEEproof}
\label{proof-lemma-ARC-definition1}
Property \ref{d1} is analogous to those used in \cite{b-ZhaoxianWu-2023,b-RunhuaWang-2023,b-Haoxiang-Ye-2023,b-Haoxiang-Ye-2025}, with a key difference: we adopt the value of $\sum_{j\in \mathcal{N}_{i}\cap\mathcal{H}\cup\{i\}}e_{ij}\|\lambda_{j}-\bar{\lambda}_{i}\|^{2}$ as the proximity measure, whereas the cited works use the value of $\max_{j\in \mathcal{N}_{i}\cap\mathcal{H}\cup\{i\}}\|\lambda_{j}-\bar{\lambda}_{i}\|^{2}$. This modification facilitates our following analysis.
For any benign agent $i \in \mathcal{H}$, we denote the set $S_{i}:=\mathcal{N}_{i}\cap \mathcal{H} \cup \{i\}$. Let $S_{i}^{c}:=\{j\in \mathcal{N}_{i}\cap \mathcal{H} \cup \{i\}, \|\lambda_{j}\| \ge C_{i} \}$ be the set of indices of clipped variables in $S_{i}$. For any dual variable $\lambda_{j}, \forall j \in \mathcal{N}_{i}$, denote $y_{j}:=clip_{C_{i}}(\lambda_{j})=\min(1,\frac{C_{i}}{\|\lambda_{j}\|})\lambda_{j}$. For benign agent $i$, $y_{i}:=\lambda_{i}$. The  weighted average of dual variables from benign agent $i$ own and its neighbors is denoted as $\bar{y}_{i}:=\sum_{j\in \mathcal{N}_{i}\cap \mathcal{H} \cup \{i\}} e_{ij}y_{j}$.

When the base robust aggregation rule namely, $IOS(\cdot)$, or $SCC(\cdot)$ satisfies Property \ref{d1} and performs perfectly, i.e., $\tilde{\rho}=0$, we have $\|AGG(\lambda_{i},\{\check{\lambda}_{j}\}_{j\in \mathcal{N}_{i}})-\bar{y}_{i}\|^{2} =0 $, i.e., $AGG(\lambda_{i},\{\check{\lambda}_{j}\}_{j\in \mathcal{N}_{i}})=\bar{y}_{i}$. Thus, we have
\begin{align}
\label{definition1-AGG-ARC-1}
&\|AGG (\lambda_{i},\{\check{\lambda}_{j}\}_{j\in \mathcal{N}_{i}})-\bar{\lambda}_{i}\|^{2}\\
=& \|AGG (\lambda_{i},\{\check{\lambda}_{j}\}_{j\in \mathcal{N}_{i}})-\bar{y}_{i}+\bar{y}_{i}-\bar{\lambda}_{i}\|^{2} \notag\\
=&\|\bar{y}_{i}-\bar{\lambda}_{i}\|^{2}.\notag
\end{align}
Substituting \eqref{eq-Bias-clipping} in Lemma \ref{lemma-Bias-clipping} into \eqref{definition1-AGG-ARC-1}, we have
\begin{align}
\label{definition1-AGG-ARC-2}
&\|AGG (\lambda_{i},\{\check{\lambda}_{j}\}_{j\in \mathcal{N}_{i}})-\bar{\lambda}_{i}\|^{2}\\
\le& \frac{\sum_{j\in S_{i}^{c}}e_{ij}}{1-\sum_{j\in S_{i}^{c}}e_{ij}}\cdot \sum_{j\in S_{i}}e_{ij}\|\lambda_{j}-\bar{\lambda}_{i}\|^{2} \notag.
\end{align}
Hence, when the base robust aggregation rule namely,  $IOS(\cdot)$, or $SCC(\cdot)$ satisfies Property \ref{d1} with contraction constant $\tilde{\rho}=0$, then the robust aggregation rule $AGG(\cdot)= IOS(ARC(\cdot))$, or $SCC(ARC(\cdot))$ satisfies Property \ref{d1} with contraction constant $\rho=\tilde{\rho}+\frac{\sum_{j\in S_{i}^{c}} e_{ij}}{1-\sum_{j\in S_{i}^{c}} e_{ij}}=0+\frac{\sum_{j\in S_{i}^{c}} e_{ij}}{1-\sum_{j\in S_{i}^{c}} e_{ij}}$.

When the base robust aggregation rule namely,  $IOS(\cdot)$, or $SCC(\cdot)$ satisfies Property \ref{d1} with a contraction constant $\tilde{\rho}>0$, using the inequality $\|a+b\|^{2}\le (1+u) \|a\|^{2}+(1+\frac{1}{u})\|b\|^{2} (u>0)$, we obtain
\begin{align}
\label{definition1-AGG-ARC-3}
&\|AGG  (\lambda_{i},\{\check{\lambda}_{j}\}_{j\in \mathcal{N}_{i}})-\bar{\lambda}_{i}\|^{2}\\
=& \|AGG  (\lambda_{i},\{\check{\lambda}_{j}\}_{j\in \mathcal{N}_{i}})-\bar{y}_{i}+\bar{y}_{i}-\bar{\lambda}_{i}\|^{2} \notag\\
\le& (1+u)\|AGG  (\lambda_{i},\{\check{\lambda}_{j}\}_{j\in \mathcal{N}_{i}})-\bar{y}_{i}\|^{2} + (1+\frac{1}{u})\|\bar{y}_{i}-\bar{\lambda}_{i}\|^{2}\notag\\
=& (\tilde{\rho}+\frac{\sum_{j\in S_{i}^{c}}e_{ij}}{1-\sum_{j\in S_{i}^{c}}e_{ij}})[\sum_{j\in S_{i}^{c}}e_{ij}\|y_{j}-\bar{y}_{i}\|^{2}\notag\\
&+\frac{1-\sum_{j\in S_{i}^{c}}e_{ij}}{\sum_{j\in S_{i}^{c}}e_{ij}}\|\bar{y}_{i}-\bar{\lambda}_{i}\|^{2}]\notag,
\end{align}
where the last equality holds by choosing $u=\frac{\sum_{j\in S_{i}^{c}}e_{ij}}{\tilde{\rho}\cdot(1-\sum_{j\in S_{i}^{c}}e_{ij})}$.

When $\|\bar{\lambda}_{i}\|\le C_{i}$, combining \eqref{eq-Variance-reduction-case1} in Lemma \ref{lemma-Variance-reduction} and \eqref{BC-2} in Lemma \ref{lemma-Bias-clipping}, we have
\begin{align}
\label{definition1-AGG-ARC-4}
&\|AGG (\lambda_{i},\{\check{\lambda}_{j}\}_{j\in \mathcal{N}_{i}})-\bar{\lambda}_{i}\|^{2}\\
\le & (\tilde{\rho}+\frac{\sum_{j\in S_{i}^{c}}e_{ij}}{1-\sum_{j\in S_{i}^{c}}e_{ij}})[\sum_{j\in S_{i}^{c}}e_{ij}\|\lambda_{j}-\bar{\lambda}_{i}\|^{2}-\sum_{j\in S_{i}^{c}}e_{ij}(\|\lambda_{j}\|-C_{i})^{2}\notag\\
&+\frac{1-\sum_{j\in S_{i}^{c}}e_{ij}}{\sum_{j\in S_{i}^{c}}e_{ij}}\cdot \frac{\sum_{j\in S_{i}^{c}}e_{ij}}{1-\sum_{j\in S_{i}^{c}}e_{ij}}\sum_{j\in S_{i}^{c}}e_{ij}(\|\lambda_{j}\|-C_{i})^{2} ] \notag\\
=& (\tilde{\rho}+\frac{\sum_{j\in S_{i}^{c}}e_{ij}}{1-\sum_{j\in S_{i}^{c}}e_{ij}})\sum_{j\in S_{i}^{c}}e_{ij}\|\lambda_{j}-\bar{\lambda}_{i}\|^{2} \notag.
\end{align}

When $\|\bar{\lambda}_{i}\|> C_{i}$, combining \eqref{eq-Variance-reduction-case2} in Lemma \ref{lemma-Variance-reduction} and \eqref{BC-2-case2-2} in Lemma \ref{lemma-Bias-clipping}, we have
\begin{align}
\label{definition1-AGG-ARC-5}
&\|AGG (\lambda_{i},\{\check{\lambda}_{j}\}_{j\in \mathcal{N}_{i}})-\bar{\lambda}_{i}\|^{2}\\
\le & (\tilde{\rho}+\frac{\sum_{j\in S_{i}^{c}}e_{ij}}{1-\sum_{j\in S_{i}^{c}}e_{ij}})[\sum_{j\in S_{i}^{c}}e_{ij}\|\lambda_{j}-\bar{\lambda}_{i}\|^{2}\notag\\
&-(1-\sum_{j\in S_{i}^{c}}e_{ij})(\|\bar{\lambda}_{i}\|-C_{i})^{2}-\sum_{j\in S_{i}^{c}}e_{ij}(\|\lambda_{j}\|-\|\bar{\lambda}_{i}\|)^{2} \notag \\
&+ \sum_{j\in S_{i}^{c}}e_{ij}(\|\lambda_{j}\|-\|\bar{\lambda}_{i}\|)^{2}+(1-\sum_{j\in S_{i}^{c}}e_{ij})(\|\bar{\lambda}_{i}\|-C_{i})^{2}] \notag\\
=& (\tilde{\rho}+\frac{\sum_{j\in S_{i}^{c}}e_{ij}}{1-\sum_{j\in S_{i}^{c}}e_{ij}})\sum_{j\in S_{i}^{c}}e_{ij}\|\lambda_{j}-\bar{\lambda}_{i}\|^{2} \notag.
\end{align}

Hence, when the base robust aggregation rule namely, $IOS(\cdot)$, or $SCC(\cdot)$ satisfies Property \ref{d1} with a contraction constant $\tilde{\rho}$, then the robust aggregation rule $AGG(\cdot) := IOS(ARC(\cdot))$, or $SCC(ARC(\cdot))$ satisfies Property \ref{d1} with contraction constant $\rho=\tilde{\rho}+\frac{\sum_{j\in S_{i}^{c}} e_{ij}}{1-\sum_{j\in S_{i}^{c}} e_{ij}}$.
\end{IEEEproof}
\end{lemma}

\begin{lemma}
\label{lemma-Variance-reduction}
If any benign agent $i$ preprocesses its received dual variables using $ARC(\cdot)$, then the following inequalities hold:

\textit{Case 1: If  $\|\bar{\lambda}_{i}\|\le C_{i}$, we have}
\begin{align}
\label{eq-Variance-reduction-case1}
&\sum_{j\in S_{i}}e_{ij}\|y_{j}-\bar{y}_{i}\|^{2}\\
\le& \sum_{j \in S_{i}}e_{ij}\|\lambda_{j}-\bar{\lambda}_{i}\|^{2}-\sum_{j \in S_{i}^{c}}e_{ij}(\|\lambda_{j}\|-C_{i})^{2}.\notag
\end{align}

\textit{Case 2: If  $\|\bar{\lambda}_{i}\|> C_{i}$, we have}
\begin{align}
\label{eq-Variance-reduction-case2}
&\sum_{j\in S_{i}}e_{ij}\|y_{j}-\bar{y}_{i}\|^{2}\\
\le& \sum_{j \in S_{i}}e_{ij}\|\lambda_{j}-\bar{\lambda}_{i}\|^{2}-(1-\sum_{j \in S_{i}^{c}}e_{ij})(\|\bar{\lambda}_{i}\|-C_{i})^{2}\notag
\\& -\sum_{j \in S_{i}^{c}}e_{ij}(\|\lambda_{j}\|-\|\bar{\lambda}_{i}\|)^{2}.\notag
\end{align}
\begin{IEEEproof}
Based on the row stochasticity of matrix $E$ in Property \ref{d1}, we have
\begin{align}
\label{VR-1}
&\sum_{j\in S_{i}}e_{ij}\|y_{j}-\bar{y}_{i}\|^{2}\\
=& \sum_{j \in S_{i}}e_{ij}\|y_{j}-\bar{\lambda}_{i}+\bar{\lambda}_{i}-\bar{y}_{i}\|^{2}\notag\\
=& \sum_{j \in S_{i}}e_{ij} \left [\|y_{j}-\bar{\lambda}_{i}\|^{2}+\|\bar{\lambda}_{i}-\bar{y}_{i}\|^{2}+2  \left \langle y_{j}-\bar{\lambda}_{i},\bar{\lambda}_{i}-\bar{y}_{i} \right \rangle  \right ] \notag\\
=&\sum_{j \in S_{i}}e_{ij}\|y_{j}-\bar{\lambda}_{i}\|^{2}- \|\bar{\lambda}_{i}-\bar{y}_{i}\|^{2} \notag,
\end{align}
Next we analyze the term $\sum_{j \in S_{i}}e_{ij}\|y_{j}-\bar{\lambda}_{i}\|^{2}$ in \eqref{VR-1}. By the definition of $S_{i}^{c}:=\{j\in \mathcal{N}_{i}\cap \mathcal{H} \cup \{i\}, \|\lambda_{j}\| > C_{i} \}$, for $j\in S_{i}\setminus S_{i}^{c}$, $y_{j}=\lambda_{j}$. Thus, we have
\begin{align}
\label{VR-2}
&\sum_{j \in S_{i}}e_{ij}\|y_{j}-\bar{\lambda}_{i}\|^{2}\\
=&\sum_{j \in S_{i}^{c}}e_{ij}\|y_{j}-\bar{\lambda}_{i}\|^{2}+\sum_{j \in S_{i}\setminus S_{i}^{c} }e_{ij}\|\lambda_{j}-\bar{\lambda}_{i}\|^{2} \notag\\
=& \sum_{j \in S_{i}}e_{ij}\|\lambda_{j}-\bar{\lambda}_{i}\|^{2} +\sum_{j \in S_{i}^{c}}e_{ij}\|y_{j}-\bar{\lambda}_{i}\|^{2}-\sum_{j \in S_{i}^{c}}e_{ij}\|\lambda_{j}-\bar{\lambda}_{i}\|^{2} \notag\\
=&\sum_{j \in S_{i}}e_{ij}\|\lambda_{j}-\bar{\lambda}_{i}\|^{2}+\sum_{j \in S_{i}^{c}}e_{ij}[\|y_{j}-\bar{\lambda}_{i}\|^{2}-\|\lambda_{j}-\bar{\lambda}_{i}\|^{2}]  \notag.
\end{align}

Now we analyze the term $\|y_{j}-\bar{\lambda}_{i}\|^{2}-\|\lambda_{j}-\bar{\lambda}_{i}\|^{2}$ in \eqref{VR-2}. For $j\in S_{i}^{c}$, $y_{j}=\frac{C_{i}}{\|\lambda_{j}\|}\lambda_{j}$, we have
\begin{align}
\label{VR-3}
&\|y_{j}-\bar{\lambda}_{i}\|^{2}-\|\lambda_{j}-\bar{\lambda}_{i}\|^{2}\\
=& \|y_{j}\|^{2}+\|\bar{\lambda}_{i}\|^{2}-2\left \langle y_{j},\bar{\lambda}_{i} \right \rangle-\|\lambda_{j}\|^{2}-\|\bar{\lambda}_{i}\|^{2}+2\left \langle \lambda_{j},\bar{\lambda}_{i}  \right \rangle \notag\\
=& C_{i}^{2}-\|\lambda_{j}\|^{2}-2\left \langle \frac{C_{i}}{\|\lambda_{j}\|}\lambda_{j},\bar{\lambda}_{i} \right \rangle+2\left \langle \lambda_{j},\bar{\lambda}_{i}  \right \rangle\notag\\
=& C_{i}^{2}-\|\lambda_{j}\|^{2}+2 (1-\frac{C_{i}}{\|\lambda_{j}\|}) \left \langle \lambda_{j},\bar{\lambda}_{i}  \right \rangle\notag\\
=& -(\|\lambda_{j}\|-C_{i})(\|\lambda_{j}\|+C_{i})+2(\|\lambda_{j}\|-C_{i})\frac{\left \langle \lambda_{j},\bar{\lambda}_{i}  \right \rangle}{\|\lambda_{j}\|}\notag\\
=& (\|\lambda_{j}\|-C_{i})(\frac{2 \left \langle \lambda_{j},\bar{\lambda}_{i}  \right \rangle}{\|\lambda_{j}\|}-\|\lambda_{j}\|-C_{i})\notag\\
\le & (\|\lambda_{j}\|-C_{i})(\frac{2 \|\lambda_{j} \|\|\bar{\lambda}_{i} \|}{\|\lambda_{j}\|}-\|\lambda_{j}\|-C_{i})\notag\\
=& (\|\lambda_{j}\|-C_{i})(2\|\bar{\lambda}_{i}\|-\|\lambda_{j}\|-C_{i})\notag ,
\end{align}
where the last inequality holds since the fact $\|\lambda_{j}\|-C_{i}> 0, j \in S_{i}^{c}$ and the inequality $\left \langle \lambda_{j},\bar{\lambda}_{i} \right \rangle \le \|\lambda_{j}\|\|\bar{\lambda}_{i}\|$.

Substituting \eqref{VR-3} into \eqref{VR-2}, we have
\begin{align}
\label{VR-4}
&\sum_{j \in S_{i}}e_{ij}\|y_{j}-\bar{\lambda}_{i}\|^{2}\\
=&\sum_{j \in S_{i}}e_{ij}\|\lambda_{j}-\bar{\lambda}_{i}\|^{2}+\sum_{j \in S_{i}^{c}}e_{ij}(\|\lambda_{j}\|-C_{i})(2\|\bar{\lambda}_{i}\|-\|\lambda_{j}\|-C_{i}).  \notag
\end{align}

Substituting \eqref{VR-4} into \eqref{VR-1}, we have
\begin{align}
\label{VR-5}
&\sum_{j\in S_{i}}e_{ij}\|y_{j}-\bar{y}_{i}\|^{2} = \sum_{j \in S_{i}}e_{ij}\|\lambda_{j}-\bar{\lambda}_{i}\|^{2}\\
+&\sum_{j \in S_{i}^{c}}e_{ij}(\|\lambda_{j}\|-C_{i})(2\|\bar{\lambda}_{i}\|-\|\lambda_{j}\|-C_{i})  - \|\bar{\lambda}_{i}-\bar{y}_{i}\|^{2}. \notag
\end{align}

We proceed to analyze \eqref{VR-5} for these two cases: $\|\bar{\lambda}_{i}\|\le C_{i}$ and $\|\bar{\lambda}_{i}\|> C_{i}$.

\textit{Case 1: When $\|\bar{\lambda}_{i}\|\le C_{i}$, we have}
\begin{align}
\label{VR-5-case1}
&\sum_{j\in S_{i}}e_{ij}\|y_{j}-\bar{y}_{i}\|^{2}\\
\le & \sum_{j \in S_{i}}e_{ij}\|\lambda_{j}-\bar{\lambda}_{i}\|^{2}+\sum_{j \in S_{i}^{c}}e_{ij}(\|\lambda_{j}\|-C_{i})(C_{i}-\|\lambda_{j}\|)   \notag \\
& - \|\bar{\lambda}_{i}-\bar{y}_{i}\|^{2} \notag\\
\le & \sum_{j \in S_{i}}e_{ij}\|\lambda_{j}-\bar{\lambda}_{i}\|^{2}-\sum_{j \in S_{i}^{c}}e_{ij}(\|\lambda_{j}\|-C_{i})^{2}. \notag
\end{align}

\textit{Case 2: When $\|\bar{\lambda}_{i}\|> C_{i}$, we have}
\begin{align}
\label{VR-5-case2}
&\sum_{j\in S_{i}}e_{ij}\|y_{j}-\bar{y}_{i}\|^{2}\\
=& \sum_{j \in S_{i}}e_{ij}\|\lambda_{j}-\bar{\lambda}_{i}\|^{2}+\sum_{j \in S_{i}^{c}}e_{ij}[(\|\bar{\lambda}_{i}\|-C_{i})^{2}-(\|\lambda_{j}\|-\|\bar{\lambda}_{i}\|)^{2}]   \notag \\
& - \|\bar{\lambda}_{i}-\bar{y}_{i}\|^{2} \notag\\
\le& \sum_{j \in S_{i}}e_{ij}\|\lambda_{j}-\bar{\lambda}_{i}\|^{2}+(\sum_{j \in S_{i}^{c}}e_{ij}-1)(\|\bar{\lambda}_{i}\|-C_{i})^{2}\notag\\
&-\sum_{j \in S_{i}^{c}}e_{ij}(\|\lambda_{j}\|-\|\bar{\lambda}_{i}\|)^{2}  \notag\\
=& \sum_{j \in S_{i}}e_{ij}\|\lambda_{j}-\bar{\lambda}_{i}\|^{2}-(1-\sum_{j \in S_{i}^{c}}e_{ij})(\|\bar{\lambda}_{i}\|-C_{i})^{2}\notag\\
&-\sum_{j \in S_{i}^{c}}e_{ij}(\|\lambda_{j}\|-\|\bar{\lambda}_{i}\|)^{2}.\notag
\end{align}
To derive the first inequality, we use the inequality $(\|\bar{\lambda}_{i}\|-C_{i})^{2}\le (\|\bar{\lambda}_{i}\|-\|\bar{y}_{i}\|)^{2}\le \|\bar{\lambda}_{i}-\bar{y}_{i}\|^{2}$ which holds based on facts $\|\bar{\lambda}_{i}\|>C_{i}$ and $\|\bar{y}_{i}\|\le C_{i}$.
\end{IEEEproof}
\end{lemma}

\begin{lemma}
\label{lemma-Bias-clipping}
If any agent $i$ preprocesses its received dual variables using $ARC(\cdot)$, then the following inequality holds:
\begin{align}
\label{eq-Bias-clipping}
\|\bar{\lambda}_{i}-\bar{y}_{i}\|^{2}\le \frac{\sum_{j\in S_{i}^{c}}e_{ij}}{1-\sum_{j\in S_{i}^{c}}e_{ij}}\cdot \sum_{j\in S_{i}}e_{ij}\|\lambda_{j}-\bar{\lambda}_{i}\|^{2}.
\end{align}
\begin{IEEEproof}
Based on the fact $\lambda_{j}=y_{j}, \forall j \in S_{i}\setminus S_{i}^{c}$, we have
\begin{align}
\label{BC-1}
&\|\bar{\lambda}_{i}-\bar{y}_{i}\|^{2}\\
=&\|\sum_{j\in S_{i}}e_{ij}(\lambda_{j}-y_{j})\|^{2} \notag\\
=&\|\sum_{j\in S_{i}^{c}}e_{ij}(\lambda_{j}-y_{j})+\sum_{j\in S_{i}\setminus S_{i}^{c} }e_{ij}(\lambda_{j}-y_{j})\|^{2} \notag\\
=& \|\sum_{j\in S_{i}^{c}}e_{ij}(\lambda_{j}-y_{j})\|^{2} \notag,
\end{align}
Based on the definition $y_{j}=\frac{C_{i}}{\|\lambda_{j}\|}\cdot \lambda_{j}, \forall j \in S_{i}^{c} $ and Jensen's inequality, we get
\begin{align}
\label{BC-2}
\|\bar{\lambda}_{i}-\bar{y}_{i}\|^{2}\le &\sum_{j\in S_{i}^{c}}e_{ij} \cdot \sum_{j\in S_{i}^{c}}e_{ij} (\|\lambda_{j}\|-C_{i})^{2}.
\end{align}

Next we analyze \eqref{BC-2} for two cases: $\|\bar{\lambda}_{i}\|\le C_{i}$ and $\|\bar{\lambda}_{i}\|> C_{i}$.

\textit{Case 1: When $\|\bar{\lambda}_{i}\|\le C_{i}$, $\|\lambda_{j}\|-C_{i}\le \|\lambda_{j}\|-\|\bar{\lambda}_{i}\|$ holds, we have}
\begin{align}
\label{BC-2-case1}
&\|\bar{\lambda}_{i}-\bar{y}_{i}\|^{2}\\
\le & \sum_{j\in S_{i}^{c}}e_{ij} \cdot \sum_{j\in S_{i}^{c}}e_{ij} (\|\lambda_{j}\|-\|\bar{\lambda}_{i}\|)^{2}\notag\\
\le & \sum_{j\in S_{i}^{c}}e_{ij} \cdot \sum_{j\in S_{i}^{c}}e_{ij}\|\lambda_{j}-\bar{\lambda}_{i}\|^{2} \notag\\
\le & \frac{\sum_{j\in S_{i}^{c}}e_{ij}}{1-\sum_{j\in S_{i}^{c}}e_{ij}}\cdot \sum_{j\in S_{i}^{c}}e_{ij}\|\lambda_{j}-\bar{\lambda}_{i}\|^{2}  \notag \\
\le & \frac{\sum_{j\in S_{i}^{c}}e_{ij}}{1-\sum_{j\in S_{i}^{c}}e_{ij}}\cdot \sum_{j\in S_{i}}e_{ij}\|\lambda_{j}-\bar{\lambda}_{i}\|^{2}.  \notag
\end{align}

\textit{Case 2: When $\|\bar{\lambda}_{i}\|> C_{i}$, we have}
\begin{align}
\label{BC-2-case2-1}
&\|\bar{\lambda}_{i}-\bar{y}_{i}\|^{2}\\
\le & \sum_{j\in S_{i}^{c}}e_{ij} \cdot \sum_{j\in S_{i}^{c}}e_{ij} (\|\lambda_{j}\|-\|\bar{\lambda}_{i}\|+\|\bar{\lambda}_{i}\|-C_{i})^{2}\notag\\
\le &  \sum_{j\in S_{i}^{c}}e_{ij} \cdot \sum_{j\in S_{i}^{c}}e_{ij} [(1+u)(\|\lambda_{j}\|-\|\bar{\lambda}_{i}\|)^{2}+(1+\frac{1}{u})(\|\bar{\lambda}_{i}\|-C_{i})^{2}].\notag
\end{align}
To derive the last inequalities, we use $\|a+b\|^{2}\le (1+u) \|a\|^{2}+(1+\frac{1}{u})\|b\|^{2} (u>0)$.

Substituting $u= \frac{\sum_{j\in S_{i}^{c}}e_{ij}}{1-\sum_{j\in S_{i}^{c}}e_{ij}}$ into \eqref{BC-2-case2-1} and rearranging the terms, we obtain
\begin{align}
\label{BC-2-case2-2}
&\|\bar{\lambda}_{i}-\bar{y}_{i}\|^{2}\\
\le & \frac{\sum_{j\in S_{i}^{c}}e_{ij}}{1-\sum_{j\in S_{i}^{c}}e_{ij}}\cdot \sum_{j\in S_{i}^{c}}e_{ij} (\|\lambda_{j}\|-\|\bar{\lambda}_{i}\|)^{2}+\sum_{j\in S_{i}^{c}}e_{ij} (\|\bar{\lambda}_{i}\|-C_{i})^{2}  \notag\\
=& \frac{\sum_{j\in S_{i}^{c}}e_{ij}}{1-\sum_{j\in S_{i}^{c}}e_{ij}}\cdot \sum_{j\in S_{i}^{c}}e_{ij} (\|\lambda_{j}\|-\|\bar{\lambda}_{i}\|)^{2}\notag\\
&+  \frac{\sum_{j\in S_{i}^{c}}e_{ij}}{1-\sum_{j\in S_{i}^{c}}e_{ij}} \cdot (1-\sum_{j\in S_{i}^{c}}e_{ij})(\|\bar{\lambda}_{i}\|-C_{i})^{2} \notag\\
=&  \frac{\sum_{j\in S_{i}^{c}}e_{ij}}{1-\sum_{j\in S_{i}^{c}}e_{ij}}\cdot \sum_{j\in S_{i}^{c}}e_{ij} (\|\lambda_{j}\|-\|\bar{\lambda}_{i}\|)^{2}\notag\\
&+  \frac{\sum_{j\in S_{i}^{c}}e_{ij}}{1-\sum_{j\in S_{i}^{c}}e_{ij}} \cdot \sum_{j\in S_{i} \setminus  S_{i}^{c}}e_{ij}(\|\bar{\lambda}_{i}\|-C_{i})^{2} \notag,
\end{align}
where the last equality holds due to the row stochasticity of matrix $E$, i.e., $\sum_{j\in S_{i}^{c}}e_{ij}+\sum_{j\in S_{i} \setminus  S_{i}^{c}}e_{ij}=1$.

Since $\|\lambda_{j}\|\le C_{i}, \forall j \in S_{i}\setminus S_{i}^{c}$, $\|\bar{\lambda}_{i}\|-C_{i}\le \|\bar{\lambda}_{i}\|-\|\lambda_{j}\|, \forall j \in S_{i}\setminus S_{i}^{c}$ holds. Therefore, we have
\begin{align}
\label{BC-2-case2-2}
&\|\bar{\lambda}_{i}-\bar{y}_{i}\|^{2}\\
\le &  \frac{\sum_{j\in S_{i}^{c}}e_{ij}}{1-\sum_{j\in S_{i}^{c}}e_{ij}}\cdot \sum_{j\in S_{i}^{c}}e_{ij} (\|\lambda_{j}\|-\|\bar{\lambda}_{i}\|)^{2}\notag\\
&+  \frac{\sum_{j\in S_{i}^{c}}e_{ij}}{1-\sum_{j\in S_{i}^{c}}e_{ij}} \cdot \sum_{j\in S_{i} \setminus  S_{i}^{c}}e_{ij}(\|\bar{\lambda}_{i}\|-\|\lambda_{j}\|)^{2} \notag\\
\le & \frac{\sum_{j\in S_{i}^{c}}e_{ij}}{1-\sum_{j\in S_{i}^{c}}e_{ij}}\cdot \sum_{j\in S_{i}}e_{ij} \|\lambda_{j}-\bar{\lambda}_{i}\|^{2}. \notag
\end{align}
\end{IEEEproof}
\end{lemma}

\subsection{Satisfaction of Property \ref{d2} by Robust Aggregation Rules $AGG(\cdot)$ }\label{aa4}

\noindent\textbf{Robust Aggregation Rule $CTM(\cdot)$ in the scalar case}
\begin{lemma}
\label{Definition2-CTM-scalar}
For any benign agent $i$, when the dual variable is one-dimensional, i.e., $d=1$, and $CTM(\cdot)$ discards the smallest $b_i$ and largest $b_i$ received messages, the robust aggregation rule $CTM(\cdot)$ satisfies Property \ref{d2}, i.e.,
$$
\|CTM(\lambda_{i},\{\check{\lambda}_{j}\}_{j\in \mathcal{N}_{i}})\|
\le
\max_{j\in \mathcal{N}_{i}\cap \mathcal{H} \cup \{i\}}\|\lambda_{j}\|.
$$
\end{lemma}

\begin{IEEEproof}
\label{proof-supp-CTM-scalar}
When $d=1$, $CTM(\cdot)$ reduces to the scalar trimmed mean rule. Let
$m_i:=\min_{j\in \mathcal{N}_{i}\cap \mathcal{H}}\lambda_j,\
M_i:=\max_{j\in \mathcal{N}_{i}\cap \mathcal{H}}\lambda_j$.
Denote by $\mathcal{U}_i\subseteq \mathcal{N}_i$ the set of retained neighbors after discarding the smallest $b_i$ and largest $b_i$ received messages. Then,
\begin{align}
\label{proof-supp-definition2-CTM-1}
CTM(\lambda_i,\{\check{\lambda}_j\}_{j\in\mathcal{N}_i})
=\frac{1}{|\mathcal{N}_i|-2b_i+1}
\left(
\lambda_i+\sum_{j\in \mathcal{U}_i}\check{\lambda}_j
\right).
\end{align}

We claim that every retained received message lies in $[m_i,M_i]$. Indeed, if some retained $\check{\lambda}_j$ satisfies $\check{\lambda}_j>M_i$, then it must come from a Byzantine neighbor, since $M_i$ is the largest benign received value. As there are at most $b_i$ Byzantine neighbors, all received messages larger than $M_i$ must be removed by discarding the largest $b_i$ values, which is a contradiction. Hence, $\check{\lambda}_j\le M_i,\ \forall j\in \mathcal{U}_i$.
Similarly, one can show that $\check{\lambda}_j\ge m_i,\ \forall j\in \mathcal{U}_i$.
Therefore,
\begin{align}
\label{proof-supp-definition2-CTM-2}
m_i\le \check{\lambda}_j\le M_i,\qquad \forall j\in \mathcal{U}_i.
\end{align}

By \eqref{proof-supp-definition2-CTM-1} and \eqref{proof-supp-definition2-CTM-2}, all scalars involved in the final average lie in $[\min\{\lambda_i,m_i\},\ \max\{\lambda_i,M_i\}]$.
Hence, the output of $CTM(\cdot)$ also lies in this interval, and thus
\begin{align}
\|CTM(\lambda_i,\{\check{\lambda}_j\}_{j\in\mathcal{N}_i})\|
&\le
\max\{\|\lambda_i\|,\|m_i\|,\|M_i\|\}\\
&=
\max_{j\in \mathcal{N}_{i}\cap \mathcal{H} \cup \{i\}}\|\lambda_{j}\|.\notag
\end{align}
This proves that $CTM(\cdot)$ satisfies Property \ref{d2} in the scalar case.
\end{IEEEproof}

\noindent\textbf{Robust Aggregation Rule $SCC(ARC(\cdot))$}
\begin{lemma}
\label{Definition2-SCC}
For any benign agent $i$, in $SCC(\cdot)$ it chooses a clipping parameter $\tau_{i}$ to clip its received messages and in $ARC(\cdot)$ it clips $b_{i}$ messages. Then, the robust aggregation rule $SCC(ARC(\cdot))$ satisfies Property \ref{d2}, i.e.,
$$\|SCC(ARC(\lambda_{i},\{\check{\lambda}_{j}\}_{j\in \mathcal{N}_{i}}))\|\le \max_{j\in \mathcal{N}_{i}\cap \mathcal{H} \cup \{i\}}\|\lambda_{j}\|.$$
\end{lemma}
\begin{IEEEproof}
\label{proof-supp-SCC}
Based on the scheme of $SCC$, we have
\begin{align}
\label{proof-supp-definition2-SCC-1}
&\|SCC(ARC(\lambda_{i},\{\check{\lambda}_{j}\}_{j\in \mathcal{N}_{i}}))\|\\
=& \| \sum_{j\in \mathcal{N}_{i} \cup \{i\}}\tilde{e}_{ij}[\lambda_{i}+clip(clip_{C_{i}}(\check{\lambda}_{j})-\lambda_{i}, \tau_{i})]\|\notag\\
\le & \sum_{j\in \mathcal{N}_{i} \cup \{i\}}\tilde{e}_{ij} \| \lambda_{i}+clip(clip_{C_{i}}(\check{\lambda}_{j})-\lambda_{i}, \tau_{i})\| \notag,
\end{align}
where the last inequality holds by using Jensen's inequality.

Considering the definition $clip(clip_{C_{i}}(\check{\lambda}_{j})-\lambda_{i}, \tau_{i}):=\min(1,\frac{\tau_{i}}{\|clip_{C_{i}}(\check{\lambda}_{j})-\lambda_{i}\|})\cdot (clip_{C_{i}}(\check{\lambda}_{j})-\lambda_{i})$, for any agent $j\in \mathcal{N}_{i} $, if $\min(1,\frac{\tau_{i}}{\|clip_{C_{i}}(\check{\lambda}_{j})-\lambda_{i}\|})=1$, we have
\begin{align}
 &\| \lambda_{i}+clip(clip_{C_{i}}(\check{\lambda}_{j})-\lambda_{i}, \tau_{i})\| \\
=&\|\lambda_{i}+clip_{C_{i}}(\check{\lambda}_{j})-\lambda_{i}\| \notag \\
=&\|clip_{C_{i}}(\check{\lambda}_{j})\|.\notag
\end{align}
If $\min(1,\frac{\tau_{i}}{\|clip_{C_{i}}(\check{\lambda}_{j})-\lambda_{i}\|})=\frac{\tau_{i}}{\|clip_{C_{i}}(\check{\lambda}_{j})-\lambda_{i}\|}$, we have
\begin{align}
 &\| \lambda_{i}+clip(clip_{C_{i}}(\check{\lambda}_{j})-\lambda_{i}, \tau_{i})\|\\
=&\| \lambda_{i}+ \frac{\tau_{i}}{\|clip_{C_{i}}(\check{\lambda}_{j})-\lambda_{i}\|}(clip_{C_{i}}(\check{\lambda}_{j})-\lambda_{i}) \| \notag\\
=&\|(1-\frac{\tau_{i}}{\|clip_{C_{i}}(\check{\lambda}_{j})-\lambda_{i}\|})\lambda_{i}+\frac{\tau_{i}}{\|clip_{C_{i}}(\check{\lambda}_{j})-\lambda_{i}\|}clip_{C_{i}}(\check{\lambda}_{j})\| \notag\\
\le&\max\{\|\lambda_{i}\|,\max_{j\in \mathcal{N}_{i}}\|clip_{C_{i}}(\check{\lambda}_{j})\|  \}.\notag
\end{align}
Thus, we have
\begin{align}
\label{eq:xx-001}
&\| \lambda_{i}+clip(clip_{C_{i}}(\check{\lambda}_{j})-\lambda_{i}, \tau_{i})\| \\
\le &\max\{\|\lambda_{i}\|,\max_{j\in \mathcal{N}_{i}}\|clip_{C_{i}}(\check{\lambda}_{j})\|  \}. \notag
\end{align}

Combining \eqref{proof-supp-definition2-SCC-1} and \eqref{eq:xx-001}, we have
\begin{align}
\label{proof-supp-definition2-SCC-2}
&\|SCC(ARC(\lambda_{i},\{\check{\lambda}_{j}\}_{j\in \mathcal{N}_{i}}))\|\\
\le & \sum_{j\in \mathcal{N}_{i} \cup \{i\}}\tilde{e}_{ij} \cdot \max\{\|\lambda_{i}\|,\max_{j\in \mathcal{N}_{i}}\|clip_{C_{i}}(\check{\lambda}_{j})\|  \} \notag\\
\le & \max\{\|\lambda_{i}\|,\max_{j\in\mathcal{N}_{i}}\|clip_{C_{i}}(\check{\lambda}_{j})\|  \}\notag \\
\le & \max_{j\in \mathcal{N}_{i}\cap \mathcal{H}  \cup \{i\}}\|\lambda_{j}\| \notag,
\end{align}
in which the second inequality holds true due to the fact $\sum_{j\in \mathcal{N}_{i} \cup \{i\}}\tilde{e}_{ij}\le 1$. The last inequality holds since $ARC(\cdot)$ guarantees that the norm of any clipped dual variable in the set $\{clip_{C_{i}}(\check{\lambda}_{j})\}_{j\in \mathcal{N}_{i}}$ must be smaller than the maximal norm of all benign dual variables in  $\{\lambda_{j}\}_{j\in \mathcal{N}_{i}\cap \mathcal{H}}$, i.e., $\|clip_{C_{i}}(\check{\lambda}_{j})\|\le \max_{j\in \mathcal{N}_{i}\cap \mathcal{H}  }\|\lambda_{j}\|, \forall j \in \mathcal{N}_{i}$.
\end{IEEEproof}

\noindent\textbf{Robust Aggregation Rule $IOS(ARC(\cdot))$}
\begin{lemma}
\label{Definition2-IOS}
For any benign agent $i$, in $IOS(\cdot)$ it discards $b_{i}$ messages and in $ARC(\cdot)$ it clips $b_{i}$ messages. Then, the robust aggregation rule $IOS(ARC(\cdot))$ satisfies Property \ref{d2}, i.e.,
$$\|IOS(ARC(\lambda_{i},\{\check{\lambda}_{j}\}_{j\in \mathcal{N}_{i}}))\|\le \max_{j\in \mathcal{N}_{i}\cap \mathcal{H}  \cup \{i\}}\|\lambda_{j}\|.$$
\end{lemma}
\begin{IEEEproof}
\label{proof-supp-IOS}
Denote the remaining agents after $IOS(\cdot)$ as $U_{i}\subset \mathcal{N}_{i} \cup \{i\}$. Based on the schemes of $IOS(\cdot)$ and $ARC(\cdot)$, we have
\begin{align}
\label{proof-supp-definition2-IOS-1}
&\|IOS(ARC(\lambda_{i},\{\check{\lambda}_{j}\}_{j\in \mathcal{N}_{i}}))\| \\
=&\|\frac{1}{\sum_{j \in U_{i}}\tilde{e}_{ij}}\cdot [\tilde{e}_{ii}\lambda_{i}+\sum_{j\in U_{i}\setminus \{i\}}\tilde{e}_{ij}\cdot
clip_{C_{i}}(\check{\lambda}_{j})]\|\notag\\
\le & \frac{1}{\sum_{j \in U_{i}}\tilde{e}_{ij}} \cdot [\tilde{e}_{ii}\|\lambda_{i}\| + \sum_{j\in U_{i}\setminus \{i\} }\tilde{e}_{ij} \|clip_{C_{i}}(\check{\lambda}_{j})\|] \notag\\
\le & \max\{\|\lambda_{i}\|,\max_{j\in U_{i}\setminus \{i\} }\|clip_{C_{i}}(\check{\lambda}_{j})\|  \}\notag \\
\le & \max_{j\in \mathcal{N}_{i}\cap \mathcal{H}  \cup \{i\}}\|\lambda_{j}\| \notag,
\end{align}
where the last inequality holds since $ARC(\cdot)$ guarantees that the norm of any clipped dual variable in the set $\{clip_{C_{i}}(\check{\lambda}_{j})\}_{j\in \mathcal{N}_{i}}$ must be smaller than the maximal norm of all benign dual variables in $\{\lambda_{j}\}_{j\in \mathcal{N}_{i}\cap \mathcal{H}}$, i.e., $\|clip_{C_{i}}(\check{\lambda}_{j})\|\le \max_{j\in \mathcal{N}_{i}\cap \mathcal{H}  }\|\lambda_{j}\|, \forall j \in \mathcal{N}_{i}$.
\end{IEEEproof}

\section{Proof of Theorem \ref{t1}}\label{B}
\begin{IEEEproof}
\label{proof-t1}
For notational convenience, we define a function $L_{i}^{t}(P):=\left \langle P-P_{i}^{t},\nabla C_{i}^{t}(P_{i}^{t}) +\frac{\lambda_{i}^{t}}{M}  \right \rangle +\frac{1}{2 \alpha}\|P-P_{i}^{t}\|^{2}$. Therefore, the update of primal variables $P_{i}^{t+1}$ in Algorithm \ref{alg1} can be rewritten as $P_{i}^{t+1}=\arg\min_{P\in \Omega_{i}} L_{i}^{t}(P)$. Given the definition of $L_{i}^{t}(P)$, we have $\nabla^{2}L_{i}^{t}(P)=\frac{1}{\alpha}>0$. Therefore, the function $L_{i}^{t}(P)$ is $\frac{1}{\alpha}$-strongly convex. According to the definition of a strongly convex function, we obtain
\begin{align}
\label{proof-t1-1}
L_{i}^{t}(\widetilde{P}_{i}^{t*})
\ge& L_{i}^{t}(P_{i}^{t+1})+\left \langle \nabla L_{i}^{t}(P_{i}^{t+1})  , \widetilde{P}_{i}^{t*}-P_{i}^{t+1}  \right \rangle \\
&+ \frac{1}{2\alpha}\|\widetilde{P}_{i}^{t*}-P_{i}^{t+1}\|^{2}.\notag
\end{align}
Since $P_{i}^{t+1}=\arg\min_{P\in \Omega_{i}} L_{i}^{t}(P)$, we obtain the optimality condition $\left \langle \nabla L_{i}^{t}(P_{i}^{t+1})  , \widetilde{P}_{i}^{t*}-P_{i}^{t+1}  \right \rangle \ge 0$. Hence, we have
\begin{align}
\label{proof-t1-2}
L_{i}^{t}(\widetilde{P}_{i}^{t*})\ge L_{i}^{t}(P_{i}^{t+1})+ \frac{1}{2\alpha}\|\widetilde{P}_{i}^{t*}-P_{i}^{t+1}\|^{2}.
\end{align}
From the definition $L_{i}^{t}(P):=\left \langle P-P_{i}^{t},\nabla C_{i}^{t}(P_{i}^{t}) + \frac{\lambda_{i}^{t} }{M} \right \rangle +\frac{1}{2 \alpha}\|P-P_{i}^{t}\|^{2}$, we can rewrite \eqref{proof-t1-2} as
\begin{align}
\label{proof-t1-3}
&\left \langle \widetilde{P}_{i}^{t*}-P_{i}^{t},\nabla C_{i}^{t}(P_{i}^{t}) + \frac{\lambda_{i}^{t}}{M} \right \rangle +\frac{1}{2 \alpha}\|\widetilde{P}_{i}^{t*}-P_{i}^{t}\|^{2}\ge  \\
&\left \langle P_{i}^{t+1}-P_{i}^{t},\nabla C_{i}^{t}(P_{i}^{t}) +\frac{\lambda_{i}^{t}}{M} \right \rangle +\frac{1}{2 \alpha}\|P_{i}^{t+1}-P_{i}^{t}\|^{2}\notag\\
&+ \frac{1}{2\alpha}\|\widetilde{P}_{i}^{t*}-P_{i}^{t+1}\|^{2}.\notag
\end{align}
Adding $C_{i}^{t}(P_{i}^{t})$ to both sides of \eqref{proof-t1-3} and rearranging the terms, we obtain
\begin{align}
\label{proof-t1-4}
&C_{i}^{t}(P_{i}^{t})+\left \langle P_{i}^{t+1}-P_{i}^{t},\nabla C_{i}^{t}(P_{i}^{t}) + \frac{\lambda_{i}^{t}}{M} \right \rangle +\frac{1}{2 \alpha}\|P_{i}^{t+1}-P_{i}^{t}\|^{2}\notag\\
\le &C_{i}^{t}(P_{i}^{t})+\left \langle \widetilde{P}_{i}^{t*}-P_{i}^{t},\nabla C_{i}^{t}(P_{i}^{t}) + \frac{\lambda_{i}^{t}}{M} \right \rangle \notag\\
&+\frac{1}{2 \alpha}(\|\widetilde{P}_{i}^{t*}-P_{i}^{t}\|^{2}- \|\widetilde{P}_{i}^{t*}-P_{i}^{t+1}\|^{2})\notag \\
\le & C_{i}^{t}(\widetilde{P}_{i}^{t*})+\left \langle \widetilde{P}_{i}^{t*}-P_{i}^{t},\frac{\lambda_{i}^{t}}{M} \right \rangle \notag\\
&+\frac{1}{2 \alpha}(\|\widetilde{P}_{i}^{t*}-P_{i}^{t}\|^{2}- \|\widetilde{P}_{i}^{t*}-P_{i}^{t+1}\|^{2}),
\end{align}
where the last inequality holds because cost function $C_{i}^{t}(\cdot)$ is convex, i.e., $C_{i}^{t}(P_{i}^{t})+\left \langle \widetilde{P}_{i}^{t*}-P_{i}^{t},\nabla C_{i}^{t}(P_{i}^{t}) \right \rangle\le C_{i}^{t}(\widetilde{P}_{i}^{t*})$. Rearranging \eqref{proof-t1-4}, we have
\begin{align}
\label{proof-t1-5}
&C_{i}^{t}(P_{i}^{t})-C_{i}^{t}(\widetilde{P}_{i}^{t*})\\
\le& \underbrace{\left \langle \widetilde{P}_{i}^{t*}-P_{i}^{t+1}, \frac{\lambda_{i}^{t}}{M} \right \rangle}_{A_{1}} \underbrace{-\left \langle P_{i}^{t+1}-P_{i}^{t},\nabla C_{i}^{t}(P_{i}^{t}) \right \rangle}_{A_{2}} \notag\\
&+\underbrace{\frac{1}{2 \alpha}(\|\widetilde{P}_{i}^{t*}-P_{i}^{t}\|^{2}- \|\widetilde{P}_{i}^{t*}-P_{i}^{t+1}\|^{2})}_{A_{3}}-\frac{1}{2 \alpha}\|P_{i}^{t+1}-P_{i}^{t}\|^{2}.\notag
\end{align}
Next, we analyze $A_{1}$, $A_{2}$ and $A_{3}$ in turn.

\noindent\textbf{Bounding $A_{1}$:}
According to the definition $\widetilde{G}_{i}^{t}(P_{i})=\frac{1}{M} P_{i}-\frac{1}{M}D^{t}$, we obtain
\begin{align}
\label{proof-t1-6}
A_{1}=&\left \langle \widetilde{P}_{i}^{t*}-P_{i}^{t+1},\frac{\lambda_{i}^{t}}{M}\right \rangle\\
=&\left \langle \widetilde{G}_{i}^{t}(\widetilde{P}_{i}^{t*})-\widetilde{G}_{i}^{t}(P_{i}^{t+1}), \lambda_{i}^{t} \right \rangle\notag\\
=& \left \langle \lambda_{i}^{t} ,\widetilde{G}_{i}^{t}(\widetilde{P}_{i}^{t*}) \right \rangle-\left \langle \lambda_{i}^{t} ,\widetilde{G}_{i}^{t}(P_{i}^{t+1}) \right \rangle\notag\\
&+\left \langle \bar{\lambda}^{t} ,\widetilde{G}_{i}^{t}(P_{i}^{t+1}) \right \rangle-\left \langle \bar{\lambda}^{t} ,\widetilde{G}_{i}^{t}(P_{i}^{t+1}) \right \rangle\notag\\
=& \left \langle \lambda_{i}^{t} ,\widetilde{G}_{i}^{t}(\widetilde{P}_{i}^{t*}) \right \rangle+\left \langle \bar{\lambda}^{t}-\lambda_{i}^{t} ,\widetilde{G}_{i}^{t}(P_{i}^{t+1}) \right \rangle\notag\\
&-\left \langle \bar{\lambda}^{t} ,\widetilde{G}_{i}^{t}(P_{i}^{t+1}) \right \rangle+ \left \langle \bar{\lambda}^{t} ,\widetilde{G}_{i}^{t}(\widetilde{P}_{i}^{t*}) \right \rangle-\left \langle \bar{\lambda}^{t} ,\widetilde{G}_{i}^{t}(\widetilde{P}_{i}^{t*}) \right \rangle\notag \\
=& \left \langle \lambda_{i}^{t}-\bar{\lambda}^{t} ,\widetilde{G}_{i}^{t}(\widetilde{P}_{i}^{t*}) \right \rangle+\left \langle \bar{\lambda}^{t}-\lambda_{i}^{t} ,\widetilde{G}_{i}^{t}(P_{i}^{t+1}) \right \rangle\notag\\
&+ \left \langle \bar{\lambda}^{t} ,\widetilde{G}_{i}^{t}(\widetilde{P}_{i}^{t*}) \right \rangle-\left \langle \bar{\lambda}^{t} ,\widetilde{G}_{i}^{t}(P_{i}^{t+1}) \right \rangle. \notag
\end{align}

\noindent\textbf{Bounding $A_{2}$:} Under Assumption \ref{a1}, we obtain
\begin{align}
\label{proof-t1-7}
A_{2} \le & \|P_{i}^{t+1}-P_{i}^{t}\|\nabla C_{i}^{t}(P_{i}^{t})\|\\
\le & \frac{u_{1}}{2}\cdot \|P_{i}^{t+1}-P_{i}^{t}\|^{2}+\frac{1}{2u_{1}}\cdot \|\nabla C_{i}^{t}(P_{i}^{t})\|^{2}\notag\\
\le & \frac{u_{1}}{2}\cdot \|P_{i}^{t+1}-P_{i}^{t}\|^{2}+\frac{\varphi^{2}}{2u_{1}},\notag
\end{align}
where $u_{1}>0$ is any positive constant. To derive the second inequality, we use $2\left \langle a  , b  \right \rangle \le u\|a\|^{2}+\frac{1}{u}\|b\|^{2}$ for any $u>0$.

\textbf{Bounding $A_{3}$:}
Similar to the derivation of \eqref{proof-t2-8} Under Assumption \ref{a1}, we obtain
\begin{align}
\label{proof-t1-8}
A_{3}=& \frac{1}{2 \alpha}(\|\widetilde{P}_{i}^{t*}-P_{i}^{t}\|^{2}- \|\widetilde{P}_{i}^{t*}-P_{i}^{t+1}\|^{2})\\
\le & \frac{R}{\alpha} \|\widetilde{P}_{i}^{t*}-\widetilde{P}_{i}^{t-1*}\|+\frac{1}{2\alpha}(\|P_{i}^{t}-\widetilde{P}_{i}^{t-1*}\|^{2} -\|\widetilde{P}_{i}^{t*}-P_{i}^{t+1}\|^{2}).\notag
\end{align}
Substituting \eqref{proof-t1-6}, \eqref{proof-t1-7} and \eqref{proof-t1-8} into \eqref{proof-t1-5} and rearranging
the terms, we have
\begin{align}
\label{proof-t1-9}
&C_{i}^{t}(P_{i}^{t})-C_{i}^{t}(\widetilde{P}_{i}^{t*})\\
\le &(\frac{u_{1}}{2}-\frac{1}{2\alpha})\|P_{i}^{t+1}-P_{i}^{t}\|^{2}+\frac{R}{\alpha} \|\widetilde{P}_{i}^{t*}-\widetilde{P}_{i}^{t-1*}\|\notag\\
&+\frac{1}{2\alpha}(\|P_{i}^{t}-\widetilde{P}_{i}^{t-1*}\|^{2} -\|\widetilde{P}_{i}^{t*}-P_{i}^{t+1}\|^{2})+\left \langle \bar{\lambda}^{t} ,\widetilde{G}_{i}^{t}(\widetilde{P}_{i}^{t*}) \right \rangle\notag\\
&-\left \langle \bar{\lambda}^{t} ,\widetilde{G}_{i}^{t}(P_{i}^{t+1}) \right \rangle+\left \langle \lambda_{i}^{t}-\bar{\lambda}^{t} ,\widetilde{G}_{i}^{t}(\widetilde{P}_{i}^{t*}) \right \rangle\notag\\
&+\left \langle \bar{\lambda}^{t}-\lambda_{i}^{t} ,\widetilde{G}_{i}^{t}(P_{i}^{t+1}) \right \rangle+\frac{\varphi^{2}}{2u_{1}}.\notag
\end{align}
Since $\widetilde{\bm{P}}^{t*}:=[\widetilde{P}_{1}^{t*},\cdots,\widetilde{P}_{M}^{t*}]$ is the optimal solution of problem \eqref{online-decentralized-economic-dispatch-problem} at time period $t$, we have $\sum_{i=1}^{M}\widetilde{G}_{i}^{t}(\widetilde{P}_{i}^{t*})=0$. Summing over $i \in  \mathcal{M}$ on both sides of \eqref{proof-t1-9}, we have
\begin{align}
\label{proof-t1-10}
&\sum_{i\in \mathcal{M}}C_{i}^{t}(P_{i}^{t})-\sum_{i\in \mathcal{M}}C_{i}^{t}(\widetilde{P}_{i}^{t*})\\
\le &(\frac{u_{1}}{2}-\frac{1}{2\alpha})\sum_{i\in \mathcal{M}}\|P_{i}^{t+1}-P_{i}^{t}\|^{2}+\frac{R}{\alpha} \sum_{i\in \mathcal{M}}\|\widetilde{P}_{i}^{t*}-\widetilde{P}_{i}^{t-1*}\|\notag\\
&+\frac{1}{2\alpha}\sum_{i\in \mathcal{M}}(\|P_{i}^{t}-\widetilde{P}_{i}^{t-1*}\|^{2} -\|\widetilde{P}_{i}^{t*}-P_{i}^{t+1}\|^{2})\notag\\
&-\sum_{i\in \mathcal{M}}\left \langle \bar{\lambda}^{t} ,\widetilde{G}_{i}^{t}(P_{i}^{t+1}) \right \rangle +\underbrace{\sum_{i\in \mathcal{M}} \left \langle \lambda_{i}^{t}-\bar{\lambda}^{t} ,\widetilde{G}_{i}^{t}(\widetilde{P}_{i}^{t*}) \right \rangle}_{A_{4}} \notag\\
&+\underbrace{\sum_{i\in \mathcal{M}} \left \langle \bar{\lambda}^{t}-\lambda_{i}^{t} ,\widetilde{G}_{i}^{t}(P_{i}^{t+1}) \right \rangle}_{A_{5}}+\frac{\varphi^{2}M}{2u_{1}}.\notag
\end{align}
Next, we analyze $A_{4}$ and $A_{5}$ in turn.
Based on Assumption \ref{a1}, Lemma \ref{lemma8} and the fact $\sum_{i\in \mathcal{M}} \|\lambda_{i}^{t}-\bar{\lambda}^{t}\|\le \sqrt{M}\cdot\|\widetilde{\Lambda}^{t}-\frac{1}{M}\widetilde{\bm{1}}\widetilde{\bm{1}}^{\top}\widetilde{\Lambda}^{t}
\|_{F}$, we obtain
\begin{align}
\label{proof-t1-11}
A_{4}=&\sum_{i\in \mathcal{M}}\left \langle \lambda_{i}^{t}-\bar{\lambda}^{t} ,\widetilde{G}_{i}^{t}(\widetilde{P}_{i}^{t*}) \right \rangle\\
\le & \sum_{i\in \mathcal{M}} \|\lambda_{i}^{t}-\bar{\lambda}^{t}\|\|\widetilde{G}_{i}^{t}(\widetilde{P}_{i}^{t*})\|\notag \\
\le & \frac{2\widetilde{\psi}^{2}M\beta}{\widetilde{\epsilon}\sqrt{\widetilde{\epsilon}}}.\notag
\end{align}
Similar to the derivation of \eqref{proof-t1-11}, we obtain
\begin{align}
\label{proof-t1-12}
A_{5}=&\sum_{i\in \mathcal{M}} \left \langle \bar{\lambda}^{t}-\lambda_{i}^{t} ,\widetilde{G}_{i}^{t}(P_{i}^{t+1}) \right \rangle \\
\le& \frac{2\widetilde{\psi}^{2}M\beta}{\widetilde{\epsilon}\sqrt{\widetilde{\epsilon}}}\notag.
\end{align}
Substituting \eqref{proof-t1-11} and \eqref{proof-t1-12} into \eqref{proof-t1-10}, we have
\begin{align}
\label{proof-t1-13}
&\sum_{i\in \mathcal{M}}C_{i}^{t}(P_{i}^{t})-\sum_{i\in \mathcal{M}}C_{i}^{t}(\widetilde{P}_{i}^{t*})\\
\le &(\frac{u_{1}}{2}-\frac{1}{2\alpha})\sum_{i\in \mathcal{M}}\|P_{i}^{t+1}-P_{i}^{t}\|^{2}+\frac{R}{\alpha} \sum_{i\in \mathcal{M}}\|\widetilde{P}_{i}^{t*}-\widetilde{P}_{i}^{t-1*}\|\notag\\
&+\frac{1}{2\alpha}\sum_{i\in \mathcal{M}}(\|P_{i}^{t}-\widetilde{P}_{i}^{t-1*}\|^{2} -\|\widetilde{P}_{i}^{t*}-P_{i}^{t+1}\|^{2})\notag\\
&-\sum_{i\in \mathcal{M}}\left \langle \bar{\lambda}^{t} ,\widetilde{G}_{i}^{t}(P_{i}^{t+1}) \right \rangle + \frac{4\widetilde{\psi}^{2}M\beta}{\widetilde{\epsilon}\sqrt{\widetilde{\epsilon}}}+\frac{\varphi^{2}M}{2u_{1}}.\notag
\end{align}

Combining \eqref{proof-t1-13} and Lemma \ref{lemma9}, we have
\begin{align}
\label{proof-t1-14}
&\frac{\widetilde{\Delta}^{t}}{2\beta}+\sum_{i\in \mathcal{M}}C_{i}^{t}(P_{i}^{t})-\sum_{i\in \mathcal{M}}C_{i}^{t}(\widetilde{P}_{i}^{t*})\\\le &(\frac{u_{1}}{2}-\frac{1}{2\alpha})\sum_{i\in \mathcal{M}}\|P_{i}^{t+1}-P_{i}^{t}\|^{2}+\frac{R}{\alpha} \sum_{i\in \mathcal{M}}\|\widetilde{P}_{i}^{t*}-\widetilde{P}_{i}^{t-1*}\|\notag\\
&+\frac{1}{2\alpha}\sum_{i\in \mathcal{M}}(\|P_{i}^{t}-\widetilde{P}_{i}^{t-1*}\|^{2} -\|\widetilde{P}_{i}^{t*}-P_{i}^{t+1}\|^{2})\notag\\
&+\underbrace{\sum_{i\in \mathcal{M}}\left \langle \bar{\lambda}^{t}  , \widetilde{G}_{i}^{t}(P_{i}^{t})  \right \rangle-\sum_{i\in \mathcal{M}}\left \langle \bar{\lambda}^{t} ,\widetilde{G}_{i}^{t}(P_{i}^{t+1}) \right \rangle}_{A_{6}}\notag\\
&-\sum_{i\in \mathcal{M}}\left \langle \lambda  , \widetilde{G}_{i}^{t}(P_{i}^{t})  \right \rangle+ \frac{4\widetilde{\psi}^{2}M\beta}{\widetilde{\epsilon}\sqrt{\widetilde{\epsilon}}}+2\widetilde{\psi}^{2}M\beta\notag \\
&+\frac{\varphi^{2}M}{2u_{1}}+\frac{M\theta}{2}\|\lambda\|^{2}.\notag
\end{align}
Next we analyze the term $A_{6}$.

\noindent\textbf{Bounding $A_{6}$:}
According to the definition $\widetilde{G}_{i}^{t}(P_{i})=\frac{1}{M} P_{i}-\frac{1}{M}D^{t}$, Assumption \ref{a1}, Lemma \ref{lemma7} and Lemma \ref{lemma8}, we have
\begin{align}
\label{proof-t1-15}
A_{6}=&\sum_{i\in \mathcal{M}}\left \langle \bar{\lambda}^{t}  , \widetilde{G}_{i}^{t}(P_{i}^{t})  \right \rangle-\sum_{i\in \mathcal{M}}\left \langle \bar{\lambda}^{t} ,\widetilde{G}_{i}^{t}(P_{i}^{t+1}) \right \rangle\\
=&\sum_{i\in \mathcal{M}} \left \langle\bar{\lambda}^{t}  ,  \widetilde{G}_{i}^{t}(P_{i}^{t})- \widetilde{G}_{i}^{t}(P_{i}^{t+1})  \right \rangle \notag\\
=& \frac{1}{M}\sum_{i\in \mathcal{M}} \left \langle\bar{\lambda}^{t}  ,  P_{i}^{t}- P_{i}^{t+1}  \right \rangle \notag \notag \\
\le & \frac{u_{2}}{2M}\sum_{i\in \mathcal{M}} \|P_{i}^{t+1}- P_{i}^{t}\|^{2}+\frac{1}{2u_{2}} \|\bar{\lambda}^{t}\|^{2}\notag\\
\le & \frac{u_{2}}{2M}\sum_{i\in \mathcal{M}} \|P_{i}^{t+1}- P_{i}^{t}\|^{2}+\frac{1}{2u_{2}}\cdot\frac{\widetilde{\psi}^{2} }{\theta^{2}},\notag
\end{align}
where $u_{2}>0$ is any positive constant. Letting $u_{2}=\frac{M}{2\alpha}$, we can rewrite \eqref{proof-t1-15} as
\begin{align}
\label{proof-t1-16}
A_{6}=&\sum_{i\in \mathcal{M}}\left \langle \lambda_{i}^{t}  , \widetilde{G}_{i}^{t}(P_{i}^{t})  \right \rangle-\sum_{i\in \mathcal{M}}\left \langle \bar{\lambda}^{t} ,\widetilde{G}_{i}^{t}(P_{i}^{t+1}) \right \rangle\\
\le & \frac{1}{4\alpha}\sum_{i\in \mathcal{M}} \|P_{i}^{t+1}- P_{i}^{t}\|^{2}+\frac{\widetilde{\psi}^{2}\alpha}{\theta^{2}M}.\notag
\end{align}
Substituting \eqref{proof-t1-16} into \eqref{proof-t1-14} and rearranging the terms, we have
\begin{align}
\label{proof-t1-17}
&\frac{\widetilde{\Delta}^{t}}{2\beta}+\sum_{i\in \mathcal{M}}C_{i}^{t}(P_{i}^{t})-\sum_{i\in \mathcal{M}}C_{i}^{t}(\widetilde{P}_{i}^{t*})\\\le &(\frac{u_{1}}{2}+\frac{1}{4\alpha}-\frac{1}{2\alpha})\sum_{i\in \mathcal{M}}\|P_{i}^{t+1}-P_{i}^{t}\|^{2}+\frac{R}{\alpha} \sum_{i\in \mathcal{M}}\|\widetilde{P}_{i}^{t*}-\widetilde{P}_{i}^{t-1*}\|\notag\\
&+\frac{1}{2\alpha}\sum_{i\in \mathcal{M}}(\|P_{i}^{t}-\widetilde{P}_{i}^{t-1*}\|^{2} -\|\widetilde{P}_{i}^{t*}-P_{i}^{t+1}\|^{2})\notag\\
&+ (\frac{4\widetilde{\psi}^{2}M}{\widetilde{\epsilon}\sqrt{\widetilde{\epsilon}}}+2\widetilde{\psi}^{2}M)\cdot\beta+\frac{\varphi^{2}M}{2u_{1}}+\frac{\widetilde{\psi}^{2}\alpha}{\theta^{2}M}\notag\\
&-\sum_{i\in \mathcal{M}}\left \langle \lambda  , \widetilde{G}_{i}^{t}(P_{i}^{t})  \right \rangle+\frac{M\theta}{2}\|\lambda\|^{2}\notag\\
=&\frac{R}{\alpha} \sum_{i\in \mathcal{M}}\|\widetilde{P}_{i}^{t*}-\widetilde{P}_{i}^{t-1*}\|\notag\\
&+\frac{1}{2\alpha}\sum_{i\in \mathcal{M}}(\|P_{i}^{t}-\widetilde{P}_{i}^{t-1*}\|^{2} -\|\widetilde{P}_{i}^{t*}-P_{i}^{t+1}\|^{2})\notag\\
&+\varphi^{2}M\cdot \alpha+ (\frac{4\widetilde{\psi}^{2}M}{\widetilde{\epsilon}\sqrt{\widetilde{\epsilon}}}+2\widetilde{\psi}^{2}M)\cdot\beta+\frac{\widetilde{\psi}^{2}}{M}\cdot\frac{\alpha}{\theta^{2}}\notag\\
&-\sum_{i\in \mathcal{M}}\left \langle \lambda  , \widetilde{G}_{i}^{t}(P_{i}^{t})  \right \rangle+\frac{M\theta}{2}\|\lambda\|^{2},\notag
\end{align}
where the last equality holds by setting $u_{1}=\frac{1}{2\alpha}$.  Summing over $t \in  [1,T]$ on both sides of \eqref{proof-t1-17}, we have
\begin{align}
\label{proof-t1-18}
Reg^{T}_{\mathcal{M}}\le &\frac{R}{\alpha}\sum_{t=1}^{T} \sum_{i\in \mathcal{M}}\|\widetilde{P}_{i}^{t*}-\widetilde{P}_{i}^{t-1*}\|\\
&+\underbrace{\frac{1}{2\alpha}\sum_{t=1}^{T}\sum_{i\in \mathcal{M}}(\|P_{i}^{t}-\widetilde{P}_{i}^{t-1*}\|^{2} -\|\widetilde{P}_{i}^{t*}-P_{i}^{t+1}\|^{2})}_{A_{7}}\notag\\
&+\varphi^{2}M\cdot \alpha T+ (\frac{4\widetilde{\psi}^{2}M}{\widetilde{\epsilon}\sqrt{\widetilde{\epsilon}}}+2\widetilde{\psi}^{2}M)\cdot\beta T+\frac{\widetilde{\psi}^{2}}{M}\cdot\frac{\alpha T}{\theta^{2}}\notag\\
&-\sum_{t=1}^{T}\sum_{i\in \mathcal{M}}\left \langle \lambda  , \widetilde{G}_{i}^{t}(P_{i}^{t})  \right \rangle+\frac{M\theta T}{2}\|\lambda\|^{2}\underbrace{-\sum_{t=1}^{T}\frac{\widetilde{\Delta}^{t}}{2\beta}}_{A_{8}}.\notag
\end{align}
Next, we analyze the terms $A_{7}$ and $A_{8}$ in turn.

\noindent\textbf{Bounding $A_{7}$:} Similar to the derivation of \eqref{proof-t2-19}, we have
\begin{align}
\label{proof-t1-19}
A_{7}=&\frac{1}{2\alpha}\sum_{t=1}^{T}\sum_{i\in \mathcal{M}}(\|P_{i}^{t}-\widetilde{P}_{i}^{t-1*}\|^{2} -\|\widetilde{P}_{i}^{t*}-P_{i}^{t+1}\|^{2})\\
\le & \frac{1}{2\alpha} \sum_{i\in \mathcal{M}}\|P_{i}^{1}-\widetilde{P}_{i}^{0*}\|^{2}.\notag
\end{align}

\noindent\textbf{Bounding $A_{8}$:}
Similar to the derivation of \eqref{proof-t2-20}, according to the definition $\widetilde{\Delta}^{t}:=M \|\bar{\lambda}^{t+1}-\lambda\|^{2}-M \|\bar{\lambda}^{t}-\lambda\|^{2}$, we have
\begin{align}
\label{proof-t1-20}
A_{6}=-\sum_{t=1}^{T}\frac{\widetilde{\Delta}^{t}}{2\beta}\le  \frac{M}{2\beta}\|\lambda\|^{2},
\end{align}

Substituting \eqref{proof-t1-19} and \eqref{proof-t1-20} into \eqref{proof-t1-18} and rearranging the terms, we have
\begin{align}
\label{proof-t1-21}
\hspace{-1em} &Reg^{T}_{\mathcal{M}}+\sum_{t=1}^{T}\sum_{i\in \mathcal{M}}\left \langle \lambda  , \widetilde{G}_{i}^{t}(P_{i}^{t})  \right \rangle-(\frac{M\theta T}{2}+\frac{M}{2\beta})\|\lambda\|^{2}\\
\hspace{-1em} \le &\frac{R}{\alpha}\sum_{t=1}^{T} \sum_{i\in \mathcal{M}}\|\widetilde{P}_{i}^{t*}-\widetilde{P}_{i}^{t-1*}\|+\frac{1}{2\alpha} \sum_{i\in \mathcal{M}}\|P_{i}^{1}-\widetilde{P}_{i}^{0*}\|^{2} \notag\\
\hspace{-1em} &+\varphi^{2}M\cdot \alpha T+ (\frac{4\widetilde{\psi}^{2}M}{\widetilde{\epsilon}\sqrt{\widetilde{\epsilon}}}+2\widetilde{\psi}^{2}M)\cdot\beta T+\frac{\widetilde{\psi}^{2}}{M}\cdot\frac{\alpha T}{\theta^{2}}.\notag
\end{align}

i) Substituting $\lambda=0$ into \eqref{proof-t1-21} and rearranging the terms, we have
\begin{align}
\label{proof-t1-22}
&Reg^{T}_{\mathcal{M}}\\
\le &\frac{R}{\alpha}\sum_{t=1}^{T} \sum_{i\in \mathcal{M}}\|\widetilde{P}_{i}^{t*}-\widetilde{P}_{i}^{t-1*}\|+\frac{1}{2\alpha} \sum_{i\in \mathcal{M}}\|P_{i}^{1}-\widetilde{P}_{i}^{0*}\|^{2}\notag\\
&+\varphi^{2}M\cdot \alpha T+ (\frac{4\widetilde{\psi}^{2}M}{\widetilde{\epsilon}\sqrt{\widetilde{\epsilon}}}+2\widetilde{\psi}^{2}M)\cdot\beta T+\frac{\widetilde{\psi}^{2}}{M}\cdot\frac{\alpha T}{\theta^{2}}\notag.
\end{align}

ii) Substituting $\lambda=\frac{\sum_{t=1}^{T}\sum_{i\in \mathcal{M}}\widetilde{G}_{i}^{t}(P_{i}^{t})}{2(\frac{M\theta T}{2}+\frac{M}{2 \beta})}$ into \eqref{proof-t1-21} and rearranging the terms, we have
\begin{align}
\label{proof-t1-23}
&\|\sum_{t=1}^{T}\sum_{i\in \mathcal{M}}\widetilde{G}_{i}^{t}(P_{i}^{t})\|^{2}\\
\le &[2M\theta T+\frac{2M}{\beta}] \cdot [\frac{R}{\alpha}\sum_{t=1}^{T} \sum_{i\in \mathcal{M}}\|\widetilde{P}_{i}^{t*}-\widetilde{P}_{i}^{t-1*}\|\notag\\
&+\frac{1}{2\alpha} \sum_{i\in \mathcal{M}}\|P_{i}^{1}-\widetilde{P}_{i}^{0*}\|^{2}+\varphi^{2}M\cdot \alpha T+\frac{\widetilde{\psi}^{2}}{M}\cdot\frac{\alpha T}{\theta^{2}}\notag\\
&+ (\frac{4\widetilde{\psi}^{2}M}{\widetilde{\epsilon}\sqrt{\widetilde{\epsilon}}}+2\widetilde{\psi}^{2}M)\cdot\beta T+2MF\cdot T].\notag
\end{align}
To derive the above inequality, we use the fact $|Reg^{T}_{\mathcal{M}}|\le 2	MF\cdot T$ which holds based on Assumption \ref{a1}.
\end{IEEEproof}

\noindent\textbf{Supporting Lemmas for Proof of Theorem \ref{t1}}
\begin{lemma}
\label{lemma7}
Under Assumptions \ref{a1} and \ref{a3}, for any agent $i\in \mathcal{M}$ and $t\in [0,\cdots,T]$, $\lambda_{i}^{t+\frac{1}{2}}$ and $\lambda_{i}^{t+1}$ generated by Algorithm \ref{alg1} satisfy
\begin{align}
\label{l7}
\|\lambda_{i}^{t+\frac{1}{2}}\|\le \frac{\widetilde{\psi}}{\theta}, \quad \|\lambda_{i}^{t+1}\|\le \frac{\widetilde{\psi}}{\theta}.
\end{align}
\begin{IEEEproof}
\label{proof-l7}
Combining the initialization $P_{i}^{0}=\lambda_{i}^{0}=D^{0}=0$ and the updates of $\lambda_{i}^{t+\frac{1}{2}}$ and $\lambda_{i}^{t+1}$ in Algorithm \ref{alg1}, we have $\|\lambda_{i}^{0+\frac{1}{2}}\|=0\le \frac{\widetilde{\psi}}{\theta}$ and $\|\lambda_{i}^{0+1}\|=0\le \frac{\widetilde{\psi}}{\theta}$. Therefore, when $t=0$, the propositions $\|\lambda_{i}^{t+\frac{1}{2}}\|\le \frac{\widetilde{\psi}}{\theta}$ and $\|\lambda_{i}^{t+1}\|\le \frac{\widetilde{\psi}}{\theta}$ hold.

Next, we prove the conclusion by mathematical induction. Suppose that when $t=t^{'}$, the propositions  $\|\lambda_{i}^{t^{'}+\frac{1}{2}}\|\le \frac{\widetilde{\psi}}{\theta}$ and $\|\lambda_{i}^{t^{'}+1}\|\le \frac{\widetilde{\psi}}{\theta}$ hold. We analyze when $t=t^{'}+1$, whether $\|\lambda_{i}^{t^{'}+1+\frac{1}{2}}\|\le \frac{\widetilde{\psi}}{\theta}$ and $\|\lambda_{i}^{t^{'}+1+1}\|\le \frac{\widetilde{\psi}}{\theta}$ hold. According to the update of $\lambda_{i}^{t+\frac{1}{2}}$ in Algorithm \ref{alg1}, we have
\begin{align}
\label{proof-l7-1}
\|\lambda_{i}^{t^{'}+1+\frac{1}{2}}\|&=\|\lambda_{i}^{t^{'}+1}+\beta\cdot(\widetilde{G}_{i}^{t^{'}+1}(P_{i}^{t^{'}+1})-\theta \lambda_{i}^{t^{'}+1})\|\notag\\
&\le (1-\beta \theta)\|\lambda_{i}^{t^{'}+1}\|+\beta \|\widetilde{G}_{i}^{t^{'}+1}(P_{i}^{t^{'}+1})\|\notag\\
&\le (1-\beta \theta)\cdot \frac{\widetilde{\psi}}{\theta}+\beta\widetilde{\psi} \notag \\
&=\frac{\widetilde{\psi}}{\theta},
\end{align}
where the second inequality holds based on $\|\lambda_{i}^{t^{'}+1}\|\le \frac{\widetilde{\psi}}{\theta}$ and  Assumption \ref{a1}. According to the update of $\lambda_{i}^{t+1}$ in Algorithm \ref{alg1}, we have
\begin{align}
\label{proof-l7-2}
&\|\lambda_{i}^{t^{'}+1+1}\| \\
=&\|\sum_{j\in \mathcal{N}_{i}\cup\{i\}}\widetilde{e}_{ij}\lambda_{j}^{t^{'}+1+\frac{1}{2}}\|\notag\\
\le& \sum_{j\in \mathcal{N}_{i}\cup\{i\}}\widetilde{e}_{ij}\|\lambda_{j}^{t^{'}+1+\frac{1}{2}}\|\notag\\
\le& \sum_{j\in \mathcal{N}_{i}\cup\{i\}}\widetilde{e}_{ij}\cdot \frac{\widetilde{\psi}}{\theta}\notag\\
=& \frac{\widetilde{\psi}}{\theta}, \notag
\end{align}
where the second inequality holds according to \eqref{proof-l7-1}. To derive the last equality, we use Assumption \ref{a3} which shows $\sum_{j\in \mathcal{N}_{i}\cup\{i\}}\widetilde{e}_{ij}=1$. Hence, when $t=t^{'}+1$, $\|\lambda_{i}^{t^{'}+1+\frac{1}{2}}\|\le \frac{\widetilde{\psi}}{\theta}$ and $\|\lambda_{i}^{t^{'}+1+1}\|\le \frac{\widetilde{\psi}}{\theta}$ hold.
\end{IEEEproof}
\end{lemma}

\begin{lemma}
\label{lemma8}
Define a matrix $\widetilde{\Lambda}^{t+1}=[\cdots,\bm{\lambda}_{i}^{t+1},\cdots]\in \mathbb{R}^{M\times d} $ that collects the dual variables $\bm{\lambda}_{i}^{t+1}$ of all agents $i\in \mathcal{M}$ generated by Algorithm \ref{alg1}. Under Assumptions \ref{a1} and \ref{a3}, we have
\begin{align}
\label{l8}
\|\widetilde{\Lambda}^{t+1}-\frac{1}{M}\widetilde{\bm{1}}\widetilde{\bm{1}}^{\top}\widetilde{\Lambda}^{t+1}
\|^{2}_{F}\le \frac{4\beta^{2}\widetilde{\psi}^{2} M }{\widetilde{\epsilon}^{3}},
\end{align}
where $\widetilde{\epsilon}:=1-\widetilde{\kappa}$.
\begin{IEEEproof}
Define $\widetilde{G}^{t}(\widetilde{P}^{t})=[\cdots,\widetilde{G}_{i}^{t}(P_{i}^{t}),\cdots]\in \mathbb{R}^{M \times d}$ to collect the local constraints $\widetilde{G}_{i}^{t}(P_{i}^{t})$ of all agents $i\in \mathcal{M}$. With these notations, we can rewrite the updates of  $\bm{\lambda}_{i}^{t+1}$ and $\bm{\lambda}_{i}^{t+\frac{1}{2}}$ in Algorithm \ref{alg1} in a compact form of
\begin{align}
\label{lemma8-proof-1}
\widetilde{\Lambda}^{t+\frac{1}{2}}=\widetilde{\Lambda}^{t}+\beta \cdot (\widetilde{G}^{t}(\widetilde{P}^{t})-\theta \widetilde{\Lambda}^{t}),
\end{align}
\begin{align}
\label{lemma8-proof-2}
\widetilde{\Lambda}^{t+1}=\widetilde{E}\widetilde{\Lambda}^{t+\frac{1}{2}}.
\end{align}
Combining \eqref{lemma8-proof-1} and \eqref{lemma8-proof-2}, and also using the fact that $\widetilde{E}$ is doubly stochastic by Assumption \ref{a3}, we have
\begin{align}
\label{lemma8-proof-3}
&\|\widetilde{\Lambda}^{t+1}-\frac{1}{M}\widetilde{\bm{1}}\widetilde{\bm{1}}^{\top}\widetilde{\Lambda}^{t+1}
\|^{2}_{F}\\
=&\|\widetilde{E}(\widetilde{\Lambda}^{t}+\beta \cdot (\widetilde{G}^{t}(\widetilde{P}^{t})-\theta \widetilde{\Lambda}^{t}))\notag\\
&-\frac{1}{M}\widetilde{\bm{1}}\widetilde{\bm{1}}^{\top}\widetilde{E}(\widetilde{\Lambda}^{t}+\beta \cdot (\widetilde{G}^{t}(\widetilde{P}^{t})-\theta \widetilde{\Lambda}^{t}))\|_{F}^{2}\notag
\end{align}
\begin{align}
=&\|(1-\beta\theta)\widetilde{E}\widetilde{\Lambda}^{t}-(1-\beta\theta)\frac{1}{M}\widetilde{\bm{1}}\widetilde{\bm{1}}^{\top}\widetilde{\Lambda}^{t}\notag\\
&+\beta \widetilde{E} \widetilde{G}^{t}(\widetilde{P}^{t})-\beta \frac{1}{M}\widetilde{\bm{1}}\widetilde{\bm{1}}^{\top}    \widetilde{E}\widetilde{G}^{t}(\widetilde{P}^{t}) \|_{F}^{2}\notag\\
\le&\frac{(1-\beta\theta)^{2}}{1-u}\|\widetilde{E}\widetilde{\Lambda}^{t}-\frac{1}{M}\widetilde{\bm{1}}\widetilde{\bm{1}}^{\top}\widetilde{\Lambda}^{t}\|_{F}^{2}\notag\\
&+\frac{\beta^{2}}{u}\|\widetilde{E}\widetilde{G}^{t}(\widetilde{P}^{t})-\frac{1}{M}\widetilde{\bm{1}}\widetilde{\bm{1}}^{\top}\widetilde{E}\widetilde{G}^{t}(\widetilde{P}^{t})\|_{F}^{2} \notag\\
=&\frac{(1-\beta\theta)^{2}}{1-u}\|(\widetilde{E}-\frac{1}{M}\widetilde{\bm{1}}\widetilde{\bm{1}}^{\top})(\widetilde{\Lambda}^{t}-\frac{1}{M}\widetilde{\bm{1}}\widetilde{\bm{1}}^{\top}\widetilde{\Lambda}^{t})\|_{F}^{2}\notag\\
&+\frac{\beta^{2}}{u}\|(\widetilde{E}-\frac{1}{M}\widetilde{\bm{1}}\widetilde{\bm{1}}^{\top})(\widetilde{G}^{t}(\widetilde{P}^{t})-\frac{1}{M}\widetilde{\bm{1}}\widetilde{\bm{1}}^{\top}\widetilde{G}^{t}(\widetilde{P}^{t}))\|_{F}^{2} \notag\\
\le& \frac{(1-\beta\theta)^{2}}{1-u}\|\widetilde{E}-\frac{1}{M}\widetilde{\bm{1}}\widetilde{\bm{1}}^{\top}\|^{2}\|\widetilde{\Lambda}^{t}-\frac{1}{M}\widetilde{\bm{1}}\widetilde{\bm{1}}^{\top}\widetilde{\Lambda}^{t}\|_{F}^{2}\notag\\
&+\frac{\beta^{2}}{u}\|\widetilde{E}-\frac{1}{M}\widetilde{\bm{1}}\widetilde{\bm{1}}^{\top}\|^{2}\|\widetilde{G}^{t}(\widetilde{P}^{t})-\frac{1}{M}\widetilde{\bm{1}}\widetilde{\bm{1}}^{\top}\widetilde{G}^{t}(\widetilde{P}^{t})\|_{F}^{2}\notag\\
\le &  \frac{1}{1-u}\|\widetilde{E}-\frac{1}{M}\widetilde{\bm{1}}\widetilde{\bm{1}}^{\top}\|^{2}\|\widetilde{\Lambda}^{t}-\frac{1}{M}\widetilde{\bm{1}}\widetilde{\bm{1}}^{\top}\widetilde{\Lambda}^{t}\|_{F}^{2}\notag\\
&+\frac{\beta^{2}}{u}\|\widetilde{E}-\frac{1}{M}\widetilde{\bm{1}}\widetilde{\bm{1}}^{\top}\|^{2}\|\widetilde{G}^{t}(\widetilde{P}^{t})-\frac{1}{M}\widetilde{\bm{1}}\widetilde{\bm{1}}^{\top}\widetilde{G}^{t}(\widetilde{P}^{t})\|_{F}^{2},\notag
\end{align}
where $u\in (0,1)$ is any positive constant. To derive the second inequality, we use the fact that $\|AB\|_{F}^{2}\le \|A\|^{2}\|B\|_{F}^{2}$.
By Assumption \ref{a3}, $\widetilde{\kappa} := \|\widetilde{E}-\frac{1}{M}\widetilde{\bm{1}}\widetilde{\bm{1}}^{\top}\|^{2} < 1$. Thus, we have
\begin{align}
\label{lemma8-proof-4}
&\|\widetilde{\Lambda}^{t+1}-\frac{1}{M}\widetilde{\bm{1}}\widetilde{\bm{1}}^{\top}\widetilde{\Lambda}^{t+1}
\|^{2}_{F}\\
\le& \frac{\widetilde{\kappa}}{1-u}\|\widetilde{\Lambda}^{t}-\frac{1}{M}\widetilde{\bm{1}}\widetilde{\bm{1}}^{\top}\widetilde{\Lambda}^{t}\|_{F}^{2}\notag\\
&+\frac{\beta^{2}\widetilde{\kappa}}{u}\|\widetilde{G}^{t}(\widetilde{P}^{t})-\frac{1}{M}\widetilde{\bm{1}}\widetilde{\bm{1}}^{\top}\widetilde{G}^{t}(\widetilde{P}^{t})\|_{F}^{2}.\notag
\end{align}
We bound the term $\frac{\beta^{2}\widetilde{\kappa}}{u}\|\widetilde{G}^{t}(\widetilde{P}^{t})-\frac{1}{M}\widetilde{\bm{1}}\widetilde{\bm{1}}^{\top}\widetilde{G}^{t}(\widetilde{P}^{t})\|_{F}^{2}$ at the right-hand side of \eqref{lemma8-proof-4} as
\begin{align}
\label{lemma8-proof-5}
&\frac{\beta^{2}\widetilde{\kappa}}{u}\|\widetilde{G}^{t}(\widetilde{P}^{t})-\frac{1}{M}\widetilde{\bm{1}}\widetilde{\bm{1}}^{\top}\widetilde{G}^{t}(\widetilde{P}^{t})\|_{F}^{2}\\
=&\frac{\beta^{2}\widetilde{\kappa}}{u}\sum_{i\in \mathcal{M}}\|\widetilde{G}_{i}^{t}(P_{i}^{t})-\frac{1}{M}\sum_{i\in \mathcal{H}}\widetilde{G}_{i}^{t}(P_{i}^{t})\|^{2}\notag\\
\le & \frac{2\beta^{2}\widetilde{\kappa}}{u}\sum_{i\in \mathcal{M}}\|\widetilde{G}_{i}^{t}(P_{i}^{t})\|^{2}+\frac{2\beta^{2}\widetilde{\kappa}M}{u}\|\frac{1}{M}\sum_{i\in \mathcal{H}}\widetilde{G}_{i}^{t}(P_{i}^{t})\|^{2}\notag\\
\le & \frac{4\beta^{2}\widetilde{\psi}^{2}\widetilde{\kappa}M}{u},\notag
\end{align}
where the last inequality holds because of Assumption \ref{a1}.

Substituting \eqref{lemma8-proof-5} into \eqref{lemma8-proof-4}, we obtain
\begin{align}
\label{lemma8-proof-6}
&\|\widetilde{\Lambda}^{t+1}-\frac{1}{M}\widetilde{\bm{1}}\widetilde{\bm{1}}^{\top}\widetilde{\Lambda}^{t+1}
\|^{2}_{F}\\
\le& \frac{\widetilde{\kappa}}{1-u}\|\widetilde{\Lambda}^{t}-\frac{1}{M}\widetilde{\bm{1}}\widetilde{\bm{1}}^{\top}\widetilde{\Lambda}^{t}\|_{F}^{2}+\frac{4\beta^{2}\widetilde{\psi}^{2}\widetilde{\kappa}M}{u}\notag\\
=& (1-\widetilde{\epsilon})\cdot\frac{1}{1-u}\|\widetilde{\Lambda}^{t}-\frac{1}{M}\widetilde{\bm{1}}\widetilde{\bm{1}}^{\top}\widetilde{\Lambda}^{t}\|_{F}^{2}+(1-\widetilde{\epsilon})\cdot \frac{4\beta^{2}\widetilde{\psi}^{2}M}{u},\notag
\end{align}
where $\widetilde{\epsilon}:=1-\widetilde{\kappa}$.

Set $u=\frac{\widetilde{\epsilon}}{1+\widetilde{\epsilon}}$. Therefore, we have $\frac{1}{1-u}= 1+\widetilde{\epsilon}$. In consequence, \eqref{lemma8-proof-6} can be rewritten as
\begin{align}
\label{lemma8-proof-7}
&\|\widetilde{\Lambda}^{t+1}-\frac{1}{M}\widetilde{\bm{1}}\widetilde{\bm{1}}^{\top}\widetilde{\Lambda}^{t+1}
\|^{2}_{F}\\
\le& (1-\widetilde{\epsilon}^{2})\|\widetilde{\Lambda}^{t}-\frac{1}{M}\widetilde{\bm{1}}\widetilde{\bm{1}}^{\top}\widetilde{\Lambda}^{t}\|_{F}^{2}+\frac{4\beta^{2}\widetilde{\psi}^{2}M}{\widetilde{\epsilon}}.\notag
\end{align}
We write \eqref{lemma8-proof-7} recursively to yield
\begin{align}
\label{lemma8-proof-8}
&\|\widetilde{\Lambda}^{t+1}-\frac{1}{M}\widetilde{\bm{1}}\widetilde{\bm{1}}^{\top}\widetilde{\Lambda}^{t+1}
\|^{2}_{F}\\
\le& (1-\widetilde{\epsilon}^{2})^{t+1}\|\widetilde{\Lambda}^{0}-\frac{1}{M}\widetilde{\bm{1}}\widetilde{\bm{1}}^{\top}\widetilde{\Lambda}^{0}\|_{F}^{2}+\sum_{l=0}^{t}(1-\widetilde{\epsilon}^{2})^{t-l}\cdot\frac{4\beta^{2}\widetilde{\psi}^{2}M}{\widetilde{\epsilon}}\notag.
\end{align}
With the same initialization $\bm{\lambda}_{i}^{0}$ for all agents $i\in \mathcal{M}$, we can rewrite \eqref{lemma8-proof-8} as
\begin{align}
\label{lemma8-proof-9}
\|\widetilde{\Lambda}^{t+1}-\frac{1}{M}\widetilde{\bm{1}}\widetilde{\bm{1}}^{\top}\widetilde{\Lambda}^{t+1}
\|^{2}_{F} \le & \sum_{l=0}^{t}(1-\widetilde{\epsilon}^{2})^{t-l}\cdot\frac{4\beta^{2}\widetilde{\psi}^{2}M}{\widetilde{\epsilon}}\notag\\
\le & \frac{4\beta^{2}\widetilde{\psi}^{2}M}{\widetilde{\epsilon}^{3}}.
\end{align}
\end{IEEEproof}
\end{lemma}

\begin{lemma}
\label{lemma9}
For any agent $i\in \mathcal{M}$ and $t\in [0,\cdots,T]$, consider $\lambda_{i}^{t+1}$ generated by Algorithm \ref{alg1}. Under Assumptions \ref{a1} and \ref{a3}, we have
\begin{align}
\label{l9}
\frac{\widetilde{\Delta}^{t}}{2\beta} \le & \sum_{i\in \mathcal{M}}\left \langle \bar{\lambda}^{t}  , \widetilde{G}_{i}^{t}(P_{i}^{t})  \right \rangle -\sum_{i\in \mathcal{M}}\left \langle \lambda  , \widetilde{G}_{i}^{t}(P_{i}^{t})  \right \rangle\\
&+2\widetilde{\psi}^{2}M\beta+\frac{M\theta}{2}\|\lambda\|^{2}.\notag
\end{align}
where $\widetilde{\Delta}^{t}:=\sum_{i \in \mathcal{M}}\|\lambda_{i}^{t+1}-\lambda\|^{2}-(1-\beta\theta)\sum_{i \in \mathcal{M}}\|\lambda_{i}^{t}-\lambda\|^{2}$ and $\lambda\in \mathbb{R}^{d}$ is an arbitrary vector.
\begin{IEEEproof}
\label{proof-l9}
Combining the updates of $\lambda_{i}^{t+\frac{1}{2}}$ and $\lambda_{i}^{t+1}$ in Algorithm \ref{alg1}, we have
\begin{align}
\label{proof-l9-1}
&M\|\bar{\lambda}^{t+1}-\lambda\|^{2}\\
=&M\|\frac{1}{M}\sum_{i\in \mathcal{M}}[\sum_{j\in \mathcal{N}_{i}\cup\{i\}}\widetilde{e}_{ij}\lambda_{j}^{t}+\beta \sum_{j\in \mathcal{N}_{i}\cup\{i\}}\widetilde{e}_{ij}(\widetilde{G}_{j}^{t}(P_{j}^{t})-\theta \lambda_{j}^{t})]-\lambda\|^{2}\notag\\
=& M\|\bar{\lambda}^{t}-\lambda+\frac{\beta}{M}\sum_{i\in\mathcal{M}}(\widetilde{G}_{i}^{t}(P_{i}^{t})-\theta \lambda_{i}^{t})  \|^{2}\notag\\
=&M \|\bar{\lambda}^{t}-\lambda\|^{2}+\sum_{i\in \mathcal{M}}\beta^{2}\|\widetilde{G}_{i}^{t}(P_{i}^{t})-\theta \lambda_{i}^{t}\|^{2}\notag\\
&+2M \beta\left \langle \bar{\lambda}^{t}-\lambda  ,\frac{1}{M}\sum_{i\in \mathcal{M}}(\widetilde{G}_{i}^{t}(P_{i}^{t})-\theta \lambda_{i}^{t}) \right \rangle  \notag\\
=& M \|\bar{\lambda}^{t}-\lambda\|^{2}+ 2\beta \sum_{i\in \mathcal{M}} \left \langle \bar{\lambda}^{t}  , \widetilde{G}_{i}^{t}(P_{i}^{t})  \right \rangle -  2\beta \sum_{i\in \mathcal{M}} \left \langle \lambda  , \widetilde{G}_{i}^{t}(P_{i}^{t})  \right \rangle\notag\\
&+\underbrace{\sum_{i\in \mathcal{M}}\beta^{2}\|\widetilde{G}_{i}^{t}(P_{i}^{t})-\theta \lambda_{i}^{t}\|^{2}}_{A_{1}}\underbrace{-2M\beta \theta \left \langle \bar{\lambda}^{t}-\lambda  , \bar{\lambda}^{t}  \right \rangle }_{A_{2}}, \notag
\end{align}
where the second equality holds since the weight matrix $\widetilde{E}:=[\widetilde{e}_{ij}]$ is column stochastic which is shown in Assumption \ref{a3}.
Next, we analyze the terms $A_{1}$ and $A_{2}$ in turn.

\noindent\textbf{Bounding $A_{1}$:}
Based on inequality $\|a+b\|^{2}\le 2\|a\|^{2}+2\|b\|^{2}$, we have
\begin{align}
\label{proof-l9-2}
A_{1}=&\sum_{i\in \mathcal{M}}\beta^{2}\|\widetilde{G}_{i}^{t}(P_{i}^{t})-\theta \lambda_{i}^{t}\|^{2}\\
\le & \sum_{i\in \mathcal{M}}2\beta^{2}\|\widetilde{G}_{i}^{t}(P_{i}^{t})\|^{2}+\sum_{i\in \mathcal{M}}2\beta^{2}\theta^{2}\| \lambda_{i}^{t}\|^{2} \le  4\beta^{2}\widetilde{\psi}^{2}M,\notag
\end{align}
where the last inequality holds, since the conclusions in Assumption \ref{a1} and Lemma
\ref{lemma7}.

\noindent\textbf{Bounding $A_{2}$:}
Based on the inequality $-2\left \langle a-b  , a  \right \rangle \le \|b\|^{2}-\|a-b\|^{2}$, we obtain
\begin{align}
\label{proof-l9-3}
A_{2}=&-2M \beta\theta\left \langle \bar{\lambda}^{t}-\lambda  ,  \bar{\lambda}^{t} \right \rangle \\
\le & M\beta\theta[\|\lambda\|^{2}-\|\bar{\lambda}^{t}-\lambda\|^{2}]\notag\\
= & \beta\theta M\|\lambda\|^{2}-\beta\theta M\|\bar{\lambda}^{t}-\lambda\|^{2}\notag\\
\le & \beta\theta M\|\lambda\|^{2}.\notag
\end{align}
Substituting \eqref{proof-l9-2} and \eqref{proof-l9-3} into \eqref{proof-l9-1} and rearranging
the terms, we have
\begin{align}
\label{proof-l9-4}
&M\|\bar{\lambda}^{t+1}-\lambda\|^{2}-M\|\bar{\lambda}^{t}-\lambda\|^{2}\\
\le& \sum_{i\in \mathcal{M}}2\beta\left \langle \bar{\lambda}^{t} , \widetilde{G}_{i}^{t}(P_{i}^{t}) \right \rangle -\sum_{i\in \mathcal{M}}2\beta\left \langle \lambda , \widetilde{G}_{i}^{t}(P_{i}^{t}) \right \rangle
\notag\\
&+ 4\beta^{2}\widetilde{\psi}^{2}M+\beta\theta M\|\lambda\|^{2}.\notag
\end{align}
Defining $\widetilde{\Delta}^{t}:=M\|\bar{\lambda}^{t+1}-\lambda\|^{2}-M\|\bar{\lambda}^{t}-\lambda\|^{2}$ and dividing both sides of \eqref{proof-l9-4} by $2\beta$, we have
\begin{align}
\label{proof-l9-5}
\frac{\widetilde{\Delta}^{t}}{2\beta} \le & \sum_{i\in \mathcal{M}}\left \langle \bar{\lambda}^{t}  , \widetilde{G}_{i}^{t}(P_{i}^{t})  \right \rangle -\sum_{i\in \mathcal{M}}\left \langle \lambda  , \widetilde{G}_{i}^{t}(P_{i}^{t})  \right \rangle\\
&+2\widetilde{\psi}^{2}M\beta+\frac{M\theta}{2}\|\lambda\|^{2}.\notag
\end{align}
\end{IEEEproof}
\end{lemma}
\end{appendices}

\balance

\bibliographystyle{IEEEtran}
\bibliography{references}

\end{document}